\documentclass[11pt]{article}
\usepackage{placeins,graphicx,diagbox }

\usepackage{wrapfig} 
\usepackage{times} 
\usepackage{amsmath,tikz,pgfplots,mathtools,float}
\usepackage{amsmath}
\usepackage{subfigure}
\usepackage{multirow}
\usepackage{amssymb,latexsym}  
\usepackage{dsfont} 

\usepackage{amsmath,algorithm} 
\usepackage{amssymb}  
\usepackage{color}
\usepackage{cite}
\newcounter{algo}

\newtheorem{theorem}{Theorem}
\newtheorem{lemma}{Lemma}

\newtheorem{assumption}{Assumption}
\newtheorem{proposition}{Proposition}
\newtheorem{corollary}{Corollary}
\newtheorem{remark}{Remark}

\def\uvs#1{{{\color{black}#1}}}

\newcommand{\red}[1]{{\color{black}#1}}

\newcommand{\wh}{\widehat}

\newcommand{\Real}{\mathbb{R}}

\newcommand{\pmat}[1]{\begin{pmatrix} #1 \end{pmatrix}}

\newcommand{\epsilonbar}{{\bar \epsilon}}

 \newcommand{\remove}[1]{}

\def\Real{\mathbb{R}}

\def\argmax{\mathop{\rm argmax}}

\usepackage{verbatim}
\usepackage{ulem}
\author{Jinlong Lei  \and Uday V. Shanbhag \and Jong-Shi Pang \and Suvrajeet Sen \thanks{Lei and Shanbhag
	(corresponding) are with the
	Department of Industrial and Manufacturing Engineering, Pennsylvania
		State University, University Park, PA 16802, USA while Pang and
		Sen are with the Department of Industrial and Systems Engineering at
		University of Southern California in Los Angeles {\tt\small
		jxl800,udaybag@psu.edu, jongship,s.sen@usc.edu}. This research was partially
supported by the U.S. National Science Foundation grant CMMI-1538605 and
CAREER award CMMI-1246887 (Shanbhag).  Additionally, a preliminary version of this work appeared in the Proceedings of the IEEE Conference on Decision and Control~\cite{shanbhag16inexact}.} }

\date{}
\title{On Synchronous, Asynchronous, and Randomized Best-Response schemes for
	computing equilibria in Stochastic Nash games
}
\begin{document}
\maketitle

\begin{abstract}
In this paper, we consider a  stochastic Nash game in which each player
minimizes a parameterized expectation-valued convex objective
	function. In deterministic
regimes,  proximal best-response (BR) schemes have been shown to be convergent under a
suitable spectral property associated with the proximal  BR
map. However, a direct application of this scheme to stochastic settings
requires obtaining exact solutions to stochastic optimization problems
at each iteration. Instead, we propose an {\em inexact} generalization
of this scheme in which an inexact solution to the  BR problem
is computed in an expected-value sense  via a stochastic approximation (SA)
scheme. On the basis of this framework, we present three inexact
 BR schemes: (i) First, we propose a synchronous inexact
 BR scheme where all players simultaneously update their
strategies; (ii)  Second, we extend this to a randomized setting
where a subset of players is randomly chosen to   update their  strategies
while the other players keep their   strategies invariant; (iii)
 Third, we propose an asynchronous  scheme, where each player
chooses its update frequency while using outdated
rival-specific data in updating its strategy. Under a suitable
contractive property on the proximal  BR map, we proceed to derive a.s.
convergence of the iterates to the Nash equilibrium (NE)  for (i) and (ii) and mean-convergence for
(i)--(iii). In addition, we show that for (i)--(iii), the generated
iterates converge to the unique equilibrium in mean at a
linear rate with a prescribed constant rather than a sub-linear rate.
Finally, we establish the overall iteration complexity of the scheme in
terms of projected stochastic gradient (SG) steps for computing an
$\epsilon-$Nash equilibrium and note that in all settings, the iteration
complexity is ${\cal O}(1/\epsilon^{2(1+c) + \delta})$ where $c = 0$ in
the context of (i) and represents the positive cost of randomization in	(ii) and asynchronicity and delay in (iii). Notably, in the
synchronous regime, we achieve a near-optimal rate from the standpoint of
solving stochastic convex optimization problems by SA schemes. The schemes are
further extended to settings where players solve two-stage stochastic Nash games with linear and quadratic recourse.
 Finally,  preliminary numerics developed on  a multi-portfolio investment
 problem  and a two-stage capacity expansion game  support  the
	 rate   and complexity statements.
 \end{abstract}
\maketitle

\maketitle

%


%
%
%

\section{Introduction}
Nash games represent an important subclass of noncooperative
games~\cite{fudenberg91game,basar99dynamic} and are rooted in the
seminal work by Nash in 1950~\cite{nash50equilibrium}. In the Nash
equilibrium problem (NEP), there is a finite set of players,  where each
player aims at minimizing its own payoff function over a  player-specific strategy set, given the rivals'
strategies.  Nash's eponymous solution concept requires that at an
	equilibrium, no player  can improve its payoff by
unilaterally deviating from its equilibrium strategy.  Over the last
several decades, there has been a surge of interest in  utilizing  Nash
games to model a range of problems in control theory and decision-making
with applications in communication networks, signal processing,
electricity markets,		(cf.~\cite{basar07control,yin09nash2,pang2010design}).
In recent years, stochastic Nash equilibrium models have	found particular relevance in power markets 			\cite{ravat11characterization,shanbhag11complementarity,kannan2013addressing,kannan2011strategic,ehrenmann2011generation,gurkan09approximations}.  Motivated by those applications,  we consider   an $N$-player stochastic
Nash game, where each player  involves solving  a stochastic
constrained optimization problem parameterized by the rivals'
strategies.  We aim to  analyze  a  breadth  of
	distributed  BR schemes for computing a Nash equilibrium
	in regimes complicated by uncertainty, delay, asynchronicity, and
	randomized update rules.

Computation of Nash equilibria has been a compelling concern for the
last several decades. As such game-theoretic models have gained wider
usage in networked regimes, the need for distributed algorithms has
become paramount. In particular, such algorithms require that the
schemes respect the privacy concerns of users and are implementable over
networked regimes. For instance, in flow control and routing problems in
communication networks, gradient-based schemes have proved particularly
useful~\cite{yin09nash2,basar07control,kannan12distributed,pang2010design}.
Such schemes generally  impose   a suitable monotonicity property on the map to
ensure global convergence of iterates.
Gradient-based schemes are characterized by ease of implementation and
lower complexity  in terms of each player step.  Yet, such schemes are
not ``fully rational''  since selfish players  {may deviate} the  update schemes unless they are forced by some authority to  follow such gradient-based schemes.
 In  BR schemes, players select their  BR, given the
current strategies of its rivals' (cf.~\cite{fudenberg98theory,basar99dynamic}). While there has been some
effort to extend such schemes to engineered settings (cf.~\cite{scutari2009mimo}),
 in which the  BR of each player can be expressed in a closed
 form.  Proximal  BR schemes appear to have  been  first discussed
	 by Facchinei and Pang in 2009~\cite{FPang09}.  In subsequent work,
	 asynchronous extensions to these proximal schemes were
		 applied to a class of optimization problems~\cite{scutari2013decomposition} while in
		 \cite{scutari2014real},  a framework that  relied on
		  solving a  sequence of NEPs   was proposed  to solve a class of monotone
		 Nash games.  In regimes with
			 expectation-valued payoffs, there has been far less computational
				 research. Sampled and smoothed counterparts of the
				 iterative regularization schemes
				 presented in~\cite{kannan12distributed} were considered
				 in~\cite{KNShanbhag13} and
				 ~\cite{yousefian15selftuned}. However, much of this
				 work considered settings where the integrands of the
				 expectation were differentiable, ruling out the
				 incorporation of two-stage recourse. Further, much of
				 the past work in computing Nash equilibria has focused
				 on risk-neutral problems with the sole exception
				 being~\cite{kannan2013addressing}, where gradient
				 methods were combined with cutting-plane schemes. In recent work,
			 the last three authors have considered Nash games
				 characterized by risk-aversion as well as two-stage
				 recourse and have examined a smoothing-based sampling
				 scheme that is more aligned with sample average
				 approximation (SAA) techniques~\cite{pang2017two}.
				 \uvs{\em Yet, the question of computing a
		 Nash equilibrium via natively distributed  BR schemes in
		 stochastic regimes in settings complicated by delay,
		 asynchronicity, and two-stage recourse remains a compelling
			 open question}, motivating the current research.
		
   {\bf  Challenge and
 motivation:} A natural question is why does
it remain challenging to extend   BR schemes to stochastic
regimes?  When contending with Nash games where each player's payoff function is
expectation-valued,  computing a  BR requires solving a
stochastic optimization exactly or accurately. Unfortunately, unless this
 BR problem is convex and the expectation as well as its gradients
are available as closed-form expressions, extensions of  BR
schemes to the stochastic regime remain impractical since the
 BR problem is essentially  a stochastic optimization problem
and requires Monte-Carlo sampling schemes~\cite{shapiro09lectures};
in effect, this leads to a two-loop scheme in which the upper loop
represents the   BR iterations while the lower loop captures
the solution of the   BR problem. Practical implementations of
such schemes that provide asymptotically accurate solutions at each step remain
unavailable and at best, this avenue provides a rather coarse
approximation of the equilibrium.

 {\bf A single-loop approach:} Instead, we consider
{developing} single-loop {\em inexact}  BR schemes that are
practically implementable; here
the   BR is computed {\em  inexactly}  in an expected-value sense but require that the
inexactness sequence be driven to zero. This is achieved by utilizing a
SA scheme that utilizes an increasing  number of
 projected SG steps. Note that in comparison with a
two-loop scheme, in this setting, the ``inside loop'' requires a known
number of iterations, thereby making this scheme a single-loop scheme. While this is a
relatively simple scheme, the resulting convergence of iterates is by
no means guaranteed and  claiming optimality of the rate is far from
immediate. In fact, convergence of   BR schemes holds
under diverse settings requiring either that the proximal-response map
(defined later) admits a suitable contractive property or  that player objectives admit a potential function
	(cf.~\cite{fudenberg98theory,FPang09}). Our work  draws
	inspiration from the work by Facchinei and Pang~\cite{FPang09} in which the contractive
	property of the proximal BR map is employed as a tool to ensure convergence of the
	iterates produced by the   BR scheme.\\
	
  {\bf Contributions:} Given the inherent challenges in solving this problem in finite time, we consider
several  {\em inexact} variants  of this scheme and make the following
contributions that are summarized in Table \ref{TAB-contribute}:
\paragraph{ (i) \em Synchronous schemes:}  We  propose  a synchronous  inexact proximal   BR scheme to find  the Nash equilibrium,
 and  prove  that the  generated iterates  converge  almost surely and in mean    to the unique Nash equilibrium
   when  the proximal   BR map  is contractive  in 2-norm.   Furthermore, when the  inexactness
is dropped at a suitable rate, the rate of convergence  of the
	iterates is provably {\bf linear} or geometric rather than {\bf sub-linear}. Based on this linear rate of convergence for the
 BR iterates and  by assuming that an inexact solution is computed via a
SA,   we  derive    the  overall iteration complexity  in projected SG   steps for  computing an
{$\epsilon$-NE$_2$ (see definition \eqref{def-epsilon-NE1})}
and show that  this    overall complexity  is  ${\cal O}( 	{N}^{1+\delta/2}/\epsilon^{2+\delta})$
where $\delta>0$. {In addition, the lower bound  of the  derived  iteration complexity satisfies $\Omega(N/\epsilon^2)$.}
 Furthermore, for a specific selection of the algorithmic parameters,
 the  overall complexity  is shown to be exactly of ${\cal
	 O}(N/\epsilon^2)$, which  is optimal for the resolution of stochastic
	 convex optimization  by standard SA  schemes.

\paragraph{(ii) \em Randomized extension:} Subsequently, a randomized  inexact proximal   BR
(BR) scheme motivated by  \cite{nesterov2012efficiency,combettes2015stochastic} is proposed,   in which
a subset of  players  is  randomly chosen    to update  at each iteration. 
To be specific, player $i$ is chosen with some positive probability $p_i$ at each major iteration.
This randomized protocol allows players to initiate an update according to a local clock, e.g, the Poisson clock,
and locally   choose the inexactness  sequence. With  a suitably  selected   inexactness  sequence,
 the estimates are shown to converge to the Nash equilibrium almost
 surely and  in mean at a prescribed  linear rate  {when  the proximal   BR map  \red{is contractive}  in 2-norm.
The expectation of the total number of projected SG steps to
  compute an $\epsilon$-NE$_2$  is shown to be
   ${\cal O}\left (( \sqrt{N}/\epsilon)^{2   \ln(\tilde{\eta}_0^{-1})/\ln(\tilde{\eta}^{-1})+\delta}\right)$	for some $\delta>0$,
    where $\tilde{\eta}$   and $\tilde{\eta}_0$ are defined in \eqref{asy-bound-inexct-sequnece} and   \eqref{rand-def-eta}, respectively.}
    Further, as noted in Remark \ref{rand-remark}, this bound is worse  than
that shown for the synchronous algorithm, a consequence of the  randomized index selection  scheme  that
accommodates   much flexibility  into the   update scheme.

\paragraph{(iii) \em Asynchronous and delay-tolerant schemes:}
Synchronous algorithms  require players to update their strategies
simultaneously, but this is often difficult to mandate in
networked settings with a large collection of noncooperative  players.
Additionally, players may often not have  access to their rivals' latest strategies.
 Motivated by the  asynchronous   algorithm specified in \cite{bertsekas1989parallel},
     we propose an  asynchronous  inexact proximal   BR
	 scheme,  where each player determines its update time  while using possibly
outdated rival information. Yet, we assume that   each player  updates at least once  in  any time interval of length $B_1,$  and
that  the communication delays are  uniformly  bounded  by $B_2$.
When the proximal   BR map   is  contractive   in {$\infty$-norm},   the  iterates are proved  to converge   in mean  to the
unique equilibrium at a linear rate for   an appropriately  chosen
inexactness sequence.  Furthermore,   we  derive  the overall iteration
complexity of   the projected SG steps  for computing {an $\epsilon$-NE$_{\infty}$ (see definition \eqref{def-epsilon-NE2})}
 and show that the   complexity    bound  is of ${\cal O}\Big((1/\epsilon)^{2B_1\left(1+\left \lceil
			  \frac{B_2}{B_1}\right \rceil\right)+\delta}\Big)$ for some $\delta>0.$
  Specially, if  the players update in a cyclic fashion,  then   the   complexity
  bound  improves to   ${\cal O}\Big((1/\epsilon)^{2\left(1+ \left \lceil
			  \frac{B_2}{N}\right\rceil+\delta\right)}\Big)$.
\paragraph{(iv) \em Incorporating private two-stage recourse:}  To show that all the aforementioned     avenues can be extended to accommodate recourse-based objectives arising from two-stage stochastic  programming~\cite{dantzig1955linear},
we  consider  an extended stochastic NEP where  each player solves a
\uvs{two-stage  stochastic program with recourse.}  We separately
investigate the case of linear and quadratic recourse with a particular
accent on computing subgradients of the random recourse function. 
On the basis of the existence and boundedness of the stochastic subgradient,    the
 proposed synchronous, randomized,   and asynchronous   inexact proximal
  BR algorithms are still applicable and the  formulated
 convergence  results  hold as well. 
\paragraph{(v) \em Preliminary numerics:} Finally,  our numerical studies
on  competitive  portfolio selection  problems  and  two-stage competitive capacity expansion problems suggest  that the
empirical behavior corresponds well with the rate statements and the
complexity bounds.
 \begin{table} [!htb]
 \centering
\newcommand{\tabincell}[2]{\begin{tabular}{@{}#1@{}}#2\end{tabular}}
\scriptsize
 \centering
     \begin{tabular}{|c|c|c|c|}
        \hline
      Update scheme  &  Asymptotic convergence & Rate of convergence  & Iteration complexity    \\ \hline
      \tabincell{c} { Synchronous  (Algorithm \ref{inexact-sbr-cont})  \\ {(using $\|.\|_2$ norm) }} &   \tabincell{c} {  a.s. convergence (Proposition \ref{prp-as} )
     \\ {convergence in mean} (Proposition  \ref{prp-mean-convergence})} & geometric  (Proposition \ref{prp-linear}) &
 \tabincell{c} {     {$\epsilon$-NE$_2$:~} \\Theorem \ref{prop-iter-comp} (b): ${\cal O}\left(( \sqrt{N}/\epsilon)^{2+\delta}\right)$ \\ Corollary \ref{cor-bound}:  ${\cal O}( {N}/\epsilon^2)$ }
     \\ \hline   \tabincell{c} { Randomized (Algorithm  \ref{rand-inexact-sbr-cont})  \\ {(using $\|.\|_2$ norm) }}  &  \tabincell{c} { a.s. convergence (Lemma \ref{lem-as-asy})
      \\   convergence in mean   (Lemma \ref{rand-geometric})}
      & geometric (Lemma \ref{rand-geometric})
  &    \tabincell{c} {  {$\epsilon$-NE$_2$:~} \\
  Theorem \ref{rand-thm-complexity}: ${\cal O}\left(( \sqrt{N}/\epsilon)^{2 \ln(\tilde{\eta}_0^{-1})/\ln(\tilde{\eta}^{-1})+\delta}\right)$} \\ \hline
 \tabincell{c} {   Asynchronous  (Algorithm  \ref{asy-inexact-sbr-cont})  \\ (using $\|.\|_{\infty}$ norm) }&     convergence in mean   (Lemma \ref{asy-lem-gemo}) &
    geometric  (Lemma \ref{asy-lem-gemo}) &  \tabincell{c} {     {$\epsilon$-NE$_{\infty}$:~} \\Theorem
		\ref{asy-upper-bound}(a): ${\cal
			O}\Big((1/\epsilon)^{2B_1\left(1+\lceil
						\frac{B_2}{B_1}\rceil\right)+\delta}\Big)$ \\
			Theorem \ref{asy-upper-bound}(b): ${\cal
				O}\Big((1/\epsilon)^{2\left(1+\lceil
							\frac{B_2}{N}\rceil\right)+\delta}\Big)$ }  \\ \hline
      \end{tabular}
      \vskip 2mm
            \centering \caption{ Summary of Contributions   \label{TAB-contribute}}
\end{table}

  {\bf Organization:} The remainder of this paper is
organized as follows:  In Section \ref{sec:formulation}, we   formulate
the stochastic Nash game, provide some basic assumptions  as well as
some background on  the proximal  BR map.
 In Section~\ref{sec:II}, we introduce a synchronous  inexact  proximal
  BR scheme,  derive rate statements, and analyze
the overall iteration complexity when an inexact solution to the  BR problem is solved via a SA scheme.
A randomized   and  an asynchronous  inexact proximal   BR
scheme are proposed  in Sections \ref{sec:random}  and \ref{sec:asy}, respectively,   where
the   rate  and complexity statements established.  Finally, all
	the proposed avenues are extended to accommodate recourse-based objectives in Section \ref{sec:recourse}.
  We present some numerical results in  Section  \ref{sec:numerics} and   conclude  in Section~\ref{sec:conclusion} with a brief summary of our main findings.

  {\bf Notations:} When referring to a vector $x$, it is assumed to
be a column vector while $x^T$ denotes its transpose. Generally, $\|x\|$ {denotes}
the Euclidean vector norm, i.e., $\|x\|=\sqrt{x^Tx}$, while other norms
will be specified appropriately (such as the $1$-norm or the
$\infty$-norm).  We use $\Pi_X[x]$ to denote the Euclidean projection of
a vector $x$ on a set $X$, i.e., $\Pi_X[x]=\min_{y \in X}\|x-y\|$.   We  abbreviate ``almost surely''  by \textit{a.s.}  and  use $\mathbb{E}[{z}]$ to denote the expectation of a random variable~$z$.
 For of a square matrix $A$, we denote  by   $ \rho(A)$  the spectral radius and by  $\lambda_{\min} (A)$   the smallest eigenvalue of  $\frac{A+A^T}{2}$.
 {We  denote by  $\otimes$  and  $\mathbf{I}_m $ the   Kronecker  product and $m \times m$ identity matrix, respectively.}

\section{Problem Formulation and Preliminaries} \label{sec:formulation}
\subsection{Problem Statement}
There exists  a  set of $N $ players indexed by  $i$ \uvs{ where $i
	\in\mathcal{N} \triangleq \{1, \cdots,N\}.$}
For any $i \in \cal{N}$, the $i$th player  has a strategy set  $X_i \subseteq \mathbb{R}^{n_i}$  and a payoff function $ f_i(x_i,x_{-i})$
depending on its own strategy  $x_i$ and  on  the vector of rivals' strategies     $x_{-i}\triangleq \{x_j\}_{j \neq i}   $.
 Suppose $n \triangleq \sum_{i=1}^N n_i,$ $X \triangleq \prod_{i} X_i$ and $X_{-i}\triangleq   \prod_{j \neq i} X_j$.  Let us consider  a  stochastic setting of the
 NEP denoted by (SNash),   in which the  objective  of player $i$,  given rivals'  strategies $x_{-i}$,  is to   solve the  following constrained stochastic  program
\begin{align}
\tag{SNash$_i(x_{-i})$} \min_{x_i \in X_i}  \quad \ f_i(x_i,x_{-i})\triangleq \mathbb{E}\left[
\psi_i(x_i,x_{-i};\xi(\omega)) \right],
\end{align}
where  $\psi_i(\cdot): X \times \mathbb{R}^d \to \mathbb{R}$ is a scalar-valued  function, and the expectation is taken with respect to the
random vector  $\xi: { \Omega} \to \Real^d$  defined on the probability space $({ \Omega}, {\cal F},
\mathbb{P})$.   An NE  of  the stochastic Nash game  (SNash) $x^*=\{x_i^*\}_{i=1}^N$  satisfies the following:
\begin{align}
	x_i^* \mbox{ solves } \quad (\mbox{SNash$_i(x_{-i}^*)$}), ~~\forall  i  \in \mathcal{ N}.  \nonumber
\end{align}
    In other words,  $x^*$ is an NE if no player  can improve the payoff by unilaterally deviating from the  equilibrium strategy $x_i^*$.  For notational simplicity,  $\xi$ is used to denote $\xi(\omega)$ throughout the paper.  In this work,
 we assume that for any $i \in \mathcal{N}$,  $f_i(x_i,x_{-i})$ and $ \psi_i(x_i, x_{-i};\xi )$ is smooth  in  $x_i$  for every given $x_{-i}$ and $\xi$. Further,  there exists a stochastic oracle such that  for any $i \in \mathcal{N}$
 and  every given $x,\xi$ returns a  sample  $\nabla_{x_i} \psi_i(x_i, x_{-i};\xi )$, which   is  an unbiased estimator of $ \nabla_{x_i} f_i(x_i,x_{-i}).$  We impose the following conditions  on the  stochastic Nash game.
\begin{assumption}~\label{assump-play-prob}  Let the following hold.\\
(a) $X_i $ is  a  closed,  compact, and    convex set;\\
(b) $f_i(x_i,x_{-i})$ is convex and  twice continuously  differentiable  in $x_i$ over an open set containing  $X_i$  for every
$x_{-i} \in X_{-i}$;\\
(c) For all $x_{-i} \in X_{-i}  $ and  any $\omega \in  \Omega$, $  \psi_i(x_i, x_{-i};\xi(\omega))$ is differentiable  in
$x_i$   over an open set containing    $X_i$ such that  $ \nabla_{x_i} f_i(x_i,x_{-i})=  \mathbb{E}[\nabla_{x_i} \psi_i(x_i, x_{-i};\xi )]$;\\
(d) For any $i \in \cal{N}$ and  all $x \in X $, there exists a constant $M_i>0$ such that
 $\mathbb{E}[\|\nabla_{x_i} \psi_i(x_i,x_{-i};\xi)\|^2] \leq  M_i^2.$
\end{assumption}

\subsection{ Background on proximal  BR maps}
 In this paper,  we consider the class of stochastic Nash games in which
the proximal  BR map (which is defined subsequently) admits a  contractive property~\cite{FPang09}. Suppose $\wh{x}(y)$ is defined as follows:
\begin{align}\label{eq-prox-resp}
\wh{x}(y) \triangleq 
{\operatornamewithlimits{\mbox{argmin}}_{x \in X}}
 \left[ \sum_{i=1}^N   \mathbb{E}\left[
\psi_i(x_i,y_{-i};\xi ) \right]+ {\mu \over 2} \|x-y\|^2 \right].
\end{align}
It is easily seen that the objective function is now separable in $x_i$ and \eqref{eq-prox-resp}
reduces to a set of player-specific proximal   BR problems,
{in which the $i$th player's} problem is given by the following:
\begin{align}\label{eq-br}
\wh{x}_{i}(y) \triangleq 
{\operatornamewithlimits{\mbox{argmin}}_{x_i \in X_i}}
\left [  \mathbb{E}\left[ \psi_i(x_i,y_{-i};\xi ) \right]+ {\mu \over 2} \|x_i-y_i\|^2 \right].
\end{align}
Analogous to the avenue adopted in ~\cite{FPang09},   we may define the $N \times N$ real matrix $\Gamma=[\gamma_{ij}]_{i,j=1}^N:$ \begin{align}
&  \Gamma \triangleq \pmat{\frac{\mu}{\mu+\zeta_{1,\min}} &
	\frac{\zeta_{12,\max} }{\mu+\zeta_{1,\min}} & \hdots &
		\frac{\zeta_{1N,\max}}{\mu+\zeta_{1,\min}}\\
\frac{\zeta_{21,\max}}{\mu+\zeta_{2,\min}} &
	\frac{\mu}{\mu+\zeta_{2,\min}} & \hdots &
		\frac{\zeta_{2N,\max}}{\mu+\zeta_{2,\min}}\\
	\vdots & & \ddots & \\
\frac{\zeta_{N1,\max}}{\mu+\zeta_{N,\min}} &
	\frac{\zeta_{N2,\max}}{\mu+\zeta_{N,\min}} & \hdots &
		\frac{\mu}{\mu+\zeta_{N,\min}}} \label{matrix-hessian}
\end{align}
with
 \begin{equation}\label{minmax-twice-diff}
\zeta_{i,\min} \triangleq \inf_{x \in X} \lambda_{\min}  \left (\nabla^2_{x_i} f_i(x) \right)  \mbox{ ~and ~} \zeta_{ij,\max} \triangleq \sup_{x \in X}  \| \nabla^2_{x_ix_j} f_i(x) \|  ~~      \forall j \neq i.
\end{equation}
Then we obtain the following inequality:
\begin{align}
  \pmat{\|\wh{x}_1(y') -\wh{x}_1(y)\| \\
				\vdots \\
\|\wh{x}_N(y') - \wh{x}_N(y)\|}  \leq   \Gamma   \pmat{\|y_1'- y_1\| \\
				\vdots \\
\|y'_N - y_N\|}  .\label{cont-prox-best-resp}
\end{align}
If the spectral radius    $ \rho(\Gamma)< 1$,  then  there exist a   scalar $a\in (0,1)$ and a monotonic norm $|~\| \bullet \| ~|$   such that
\begin{align}
\left |~\left\| \pmat{\|\wh{x}_1(y') -\wh{x}_1(y)\| \\
				\vdots \\
\|\wh{x}_N(y') - \wh{x}_N(y)\|} \right\|~\right|\leq a
\left |~\left\|  \pmat{\|y_1'- y_1\| \\
				\vdots \\ \|y'_N - y_N\|} \right\|~\right| .  \label{contraction}
\end{align}
Note that sufficient conditions for the contractive property of the  BR   map $ \wh{x}(\bullet)$
are provided in~\cite[Proposition 12.17]{FPang09}.

\begin{remark}~
Let $D_i \in \mathbb{R}^{n_i \times n_i}~\forall i \in \cal{N}$ be any
arbitrary chosen   nonsingular  matrices.  Define   $$\zeta_{i,\min} \triangleq
\inf_{x \in X}~ \lambda_{\min}  \left (D_i^T \nabla^2_{x_i} f_i(x) D_i\right)
,~\zeta_{ij,\max} \triangleq \sup_{x \in X}~\| D_i^T\nabla^2_{x_ix_j} f_i(x)
D_j\|~~\forall j \neq i. $$ It is easily seen that  definition
\eqref{minmax-twice-diff} is only a special case  with  $D_i\in \mathbb{R}^{n_i
\times n_i}$  being an identity matrix. By making simple modifications to the
proof given  in \cite[Section 12.6.1]{FPang09} we  are able to obtain the
inequality \eqref{cont-prox-best-resp} with $\Gamma$ defined by
\eqref{matrix-hessian}.  As discussed   in \cite{scutari2014real},
matrices $D_i~\forall i \in \cal{N}$  provide an additional  degree of freedom
in deriving   the contraction property of  the proximal  BR map
$\wh{x}(\bullet)$.
\end{remark}

 Recall that   the proximal  BR $\wh{x}_i(y)$ requires solving a
stochastic optimization problem  defined in \eqref{eq-br}. Stochastic
optimization problems have been studied extensively over the last several decades through Monte-Carlo sampling
schemes  such as SAA and SA.  SAA provides a foundation for relating the
estimators of an expectation-valued problem obtained by solving the
deterministic sample-average approximation.   Asymptotic convergence and error bounds for the SAA
estimators  have been extensively investigated in  \cite{shapiro09lectures}.
SA  was first considered by Robbins and Monro \cite{robbins1951stochastic} for
seeking roots of a regression function with noisy observations. Such  schemes
have found wide applications in stochastic problems such as convex
optimization,  variational inequality  problems, and systems and control problems.  Here, we   apply the SA scheme  to estimate the optimal solution of the
stochastic  BR problem~\eqref{eq-br}.

\section{A  Synchronous  Inexact   Proximal  BR Scheme}\label{sec:II}
In this section, we present a  synchronous  inexact   proximal  BR scheme   for which both asymptotic convergence analysis, rate statements, and overall iteration complexity results are provided.  
\subsection{Description of algorithm}
Recall that the exact  proximal  BR $\wh{x}_{i}( y_k)$ is   defined as the  following:
\begin{align}
 \wh{x}_{i}(y_k)& \triangleq 
{\operatornamewithlimits{\mbox{argmin}}_{x_i \in X_i}}
 \big[  f_i(x_i,y_{-i,k})   + {\mu \over 2} \|x_i-y_{i,k}\|^2 \big]. \label{def-wh_x_inft}
\end{align}
Since $ f_i(x_i,y_{-i,k}) $ is the expectation of $ \psi_i(x_i,y_{-i,k};\xi ) $  with respect to $\xi$,
a closed-form expression of  $\wh{x}_{i}(y_k)$  is unavailable.  Instead  we propose an inexact proximal   BR scheme (Algorithm~\ref{inexact-sbr-cont})   computing  an {\em inexact}  BR  that satisfies
  \eqref{inexact-sub} rather than computing the exact  BR,  which is available  only in an asymptotic sense.
   It is worth pointing out that     $\mathcal{F}_k$   in equation \eqref{inexact-sub}
denotes  the $\sigma$-field of the entire information  used by  the algorithm up to  (and including)  the update of  $x_k$.
We will  further define  $\mathcal{F}_k$ in Section \ref{sec:syn:complexity} and
 specify the deterministic sequence  $\{\alpha_{i,k}\}_{k\geq 0}$  in  Algorithm \ref{inexact-sbr-cont}
when we proceed to   investigate  the convergence properties.

\begin{algorithm}[H]
\caption{{Inexact proximal  BR scheme for stochastic Nash
	games}} \label{inexact-sbr-cont}
 Set $k:=0$. Let  $ y_{i,0}=x_{i,0} \in X_i$,   and  $\{\alpha_{i,k}\}_{k \geq 0}$
be a given  deterministic  sequence for $i = 1, \hdots, N$.
\begin{enumerate}
\item[(1)] For $i = 1, \hdots, N$,  let $x_{i,k+1}$ satisfy the following
 \begin{align} \label{inexact-sub}
x_{i,k+1}  \in  \left\{  z\in X_i  : \mathbb{E}\left[\|z -
\wh{x}_{i}(y_k)\|^2\big |\mathcal{F}_k\right]   \leq \alpha_{i,k}^2~~a.s. \right\}.
\end{align}
\item[(2)] For $i = 1, \hdots, N$,  $y_{i,k+1} :=x_{i,k+1}$;
\item [(3)] $k:=k+1$; If $k < K$, return to (1); else STOP.
\end{enumerate}
\end{algorithm}
\subsection{Convergence analysis}

We now proceed to analyze this scheme in greater detail and initiate our
discussion by noting the convergence of the exact  BR scheme.
  \begin{proposition}[~\cite{FPang09}]~
Consider the stochastic Nash game in which the $i$th player solves
(SNash$_i(x_{-i})$), given $x_{-i} \in X_{-i}$.  Let the sequence of exact responses be denoted by $\{y_{k}^{\infty} \}_{k=1}^{\infty}$.  Assume that the matrix $\Gamma$  defined by  \eqref{matrix-hessian}  satisfies $
\rho(\Gamma)< 1$. Then the following hold:
\begin{enumerate}
\item[(i)] The contraction given by \eqref{contraction} holds for the proximal  BR map.
\item[(ii)] For $i = 1, \cdots, N,$ the  BR iterates $y_{i,k}^{\infty}  \to x_i^*$ as $k \to \infty$
\end{enumerate}
\end{proposition}
Naturally, one  may  derive the uniqueness of the Nash equilibrium from
the existence of this contractive property (cf.~\cite{FPang09}). We prove that the sequence
$\{x_k\}_{k \geq 0}$ generated by Algorithm~\ref{inexact-sbr-cont} converges almost surely and
in expectation to the unique Nash equilibrium of the stochastic Nash game
in Propositions \ref{prp-as} and \ref{prp-mean-convergence},  respectively.   The results rely  on comparing
the inexact responses with the  exact  BR map, which is known to be contractive by
 the following  sufficient  condition:
\begin{assumption}~\label{assp-gamma}
The $\Gamma$ matrix defined in \eqref{matrix-hessian} satisfies  $\|\Gamma\|< 1$.
\end{assumption}
   Define  $a  \triangleq  \|\Gamma\|.$   Then  $a  \in(0,1)$ by Assumption \ref{assp-gamma}, and hence
   by   \eqref{cont-prox-best-resp} we obtain the following:
\begin{align}
\left\| \pmat{\|\wh{x}_1(y') -\wh{x}_1(y)\| \\
				\vdots \\
\|\wh{x}_N(y') - \wh{x}_N(y)\|} \right\|\leq a \left\| \pmat{\|y_1'- y_1\| \\
				\vdots \\
\|y'_N - y_N\|} \right\|\quad\forall y,y' \in X . \label{cont-prox-resp}
\end{align}

We  utilize the following lemma for  random and  deterministic recursions to prove the almost sure  convergence and
 convergence in mean, respectively.
	\begin{lemma}\label{supp-conv}
	(a) \cite[Lemma 2.2.10]{Polyak87}  Let $\{v_k\}_{k \geq 0}$ be a sequence of random variables, $v_k\geq 0$, $\mathbb{E}[v_0]<\infty$, such that
$\mathbb{E}[v_{k+1}| v_0,\cdots,v_k] \leq q_kv_k+\zeta_k$, a.s. where $0 \leq q_k < 1,~
    \zeta_k \geq 0, ~ \sum_{k=0}^{\infty} (1-q_k) = \infty,  ~\sum_{k=0}^{\infty}  \zeta_k< \infty, \; {\rm  and ~}
   \lim_{k\to \infty} \frac{ \zeta_k}{1-q_k} =0.$
   Then $v_k \to 0 ~a.s.$

  \noindent  (b)\cite[Lemma 2.2.3]{Polyak87} Let $ u_{k+1} \leq q_k u_k +  \zeta_k,\, 0 \leq q_k < 1, \,
    \zeta_k \geq 0$, where $ \sum_{k=0}^{\infty} (1-q_k) = \infty \; {\rm  and ~}
   \lim_{k\to \infty} \frac{ \zeta_k}{1-q_k} =0.$ Then $   \limsup\limits_{k \to
		\infty} u_k \leq 0$. In particular,  if $u_k  \geq 0$, then $\lim\limits_{k \to \infty} u_k = 0.$
\end{lemma}
\begin{proposition}[{\bf Almost Sure Convergence}] \label{prp-as}~
Let the sequence $\{x_k\}_{k \geq 0}$ be generated
by Algorithm~\ref{inexact-sbr-cont}. Suppose Assumptions \ref{assump-play-prob} and   \ref{assp-gamma} hold,  and that $ \alpha_{i,k} \geq 0 $ with $\sum_{k=0}^\infty \alpha_{i,k}<\infty$ for  any $i \in \mathcal{N} $. Then for  any $i \in \mathcal{N} $,
  $$ \lim_{ k \to\infty}  x_{i,k} =x_i^* \quad a.s.$$
\end{proposition}
\noindent {\bf Proof.}
We begin by noting that  $x_i^* = \wh{x}_i (x^*)$  and examine $\|x_{i,k+1}- x_i^*\|$:
\begin{align}
  \quad \|x_{i,k+1}- x_i^*\|
			& \leq  \|x_{i,k+1} -\wh{x}_i (y_k)\|
			+\| \wh{x}_{i}(y_k) -
			\wh{x}_i (x^*)\|, \label{ineq1}
\end{align}
where the inequality follows from the triangle inequality.
 Furthermore,  by  $y_k=x_k$ and the contractive property \eqref{cont-prox-resp} of $\wh{x}(\bullet)$, we have  the following bound:
\begin{align}\label{proximal-map-bound}
\left\|  \pmat{
	 \|\wh{x} _{1}(y_k) -
	 \wh{x}_{1}(x^*)\|\\
				\vdots \\
	 \|\wh{x} _{N}(y_k) -
	 \wh{x}_{N}(x^*)\|} \right\|   & = \left\|  \pmat{
	 \|\wh{x} _{1}(x_k) -
	 \wh{x}_{1}(x^*)\|\\
				\vdots \\
	 \|\wh{x} _{N}(x_k) -
	 \wh{x}_{N}(x^*)\|} \right\|  \leq a  \left\|  \pmat{
	 \|x_{1,k}-x_1^*\|\\
				\vdots \\
	 \|x_{N,k}- x_N^*\|} \right\|.
\end{align}	
Consequently,  from \eqref{ineq1} by    the triangle inequality,  we obtain the following inequality:
\begin{align}\label{ineaxt-contractive}
& v_{k+1} \triangleq \left\|   \pmat{ \|x_{1,k+1} -{x}_{1}^*\| \\
				\vdots \\
 \|x_{N,k+1}  - {x}_{N}^*\|} \right \|   \leq
a \left\|  \pmat{
	 \|x_{1,k}-x_1^*\|\\
				\vdots \\
	 \|x_{N,k}- x_N^*\|} \right \| + \left\|  \pmat{
	 \|x_{1,k+1} -
			 \wh{x}_{1}(y_k)\|\\
				\vdots \\
	 \|x_{N,k+1} -
			 \wh{x}_{N}(y_k)\|} \right\|.
\end{align}

  Denote  by $\alpha_{\max,k} \triangleq \max_{i} \alpha_{i,k}$.
 Then by   \eqref{inexact-sub}  and conditional Jensen's inequality,  we have  that for all $k \geq 0$:
\begin{align}\label{inequ2}
  \mathbb{E} \left[    \left\|  \pmat{
	 \|x_{1,k+1} -
			 \wh{x}_{1}(y_k)\|\\
				\vdots \\
	 \|x_{N,k+1} -
			 \wh{x}_{N}(y_k)\|} \right\| \Bigg| \mathcal{F}_k\right]  \notag & \leq  \sqrt{ \mathbb{E} \left [\sum_{i=1}^N \|x_{i,k+1} - \wh{x}_i (y_k)\|^2\big|\mathcal{F}_k \right]}  \\ & \leq  \sqrt{ \sum_{i=1}^N \alpha_{i,k}^2  } \leq \sqrt{N}\alpha_{\max,k} \quad a.s.
\end{align}By taking expectations  conditioned on  $\mathcal{F}_k$  with respect to both sides of \eqref{ineaxt-contractive},
 from  \eqref{inequ2}  we have that
\begin{equation}
\begin{split} \mathbb{E}[v_{k+1}|\mathcal{F}_k]& \leq a v_k +  \sqrt{N}\alpha_{\max,k} \quad a.s. \label{def-uk}
\end{split}
\end{equation}
 Since $\sum_k(1-a) = \infty$ and  $  \lim\limits_{k\to \infty} \frac{\alpha_{\max,k}}{(1-a)}=0$,  by   Lemma~\ref{supp-conv}(a)
  and the summability of  $\alpha_{\max,k} $, we conclude that $ v_k \to 0 ~a.s.$
Consequently,   the result follows by $v_k$ defined   in \eqref{ineaxt-contractive}.
 \hfill $\Box$

{\begin{proposition} \label{prp-mean-convergence}~
Let the sequence $\{x_k\}_{k \geq 0}$ be generated
by Algorithm~\ref{inexact-sbr-cont}.  Suppose  Assumptions \ref{assump-play-prob} and   \ref{assp-gamma} hold, and that $ \alpha_{i,k} \geq 0$ with  $\lim\limits_{k\to \infty}  \alpha_{i,k}=0$ for  any $i \in \mathcal{N} $.
 Then we have the following  assertions:

 (a) ({\bf convergence of the variance of $x_k$}) $\lim\limits_{k\to \infty} \mathbb{V}\mbox{ar}(x_k)=0.$

 (b) ({\bf convergence in mean})
 $\lim\limits_{k\to \infty}  \mathbb{E}[\|x_{i,k} -x_i^*\|] =0~~\forall i \in \mathcal{N}.$
\end{proposition}
\noindent {\bf Proof.}
  (a) Note that the variance of $x_k$ may be bounded as follows:
\begin{equation}\label{variance-inequ0}
\begin{split}
\mathbb{V}{ar}(x_k)&=\mathbb{E}\Big[\big\|\mathbb{E}[x_k]-x_k \big\|^2\Big]
= \mathbb{E}\Big[\big\|\mathbb{E}[x_k-x^*]-\left(x_k-x^*\right) \big\|^2\Big]
\leq \mathbb{E}\big[\| x_k-x^*  \|^2\big].
\end{split}
\end{equation}
By  Assumption \ref{assump-play-prob}(a),  we define the diameter of the $X_i$   as follows:
 \begin{equation}\label{diameter}
D_{X_i}\triangleq\sup\{ d(x_i,x_i'): x_i,x_i' \in X_i \}<\infty.
\end{equation}
By  \eqref{ineaxt-contractive}, we have the following:
\begin{equation}
\begin{split}
 \|x_{ k+1}- x ^*\|^2 &   \leq  \sum_{i=1}^N\|x_{i,k+1} -\wh{x}_i (y_k)\|^2
			 +a^2  \|   x_{ k}-x^*\|^2
			 +  2 a \left\|  \pmat{
	 \|x_{1,k}-x_1^*\|\\
				\vdots \\
	 \|x_{N,k}- x_N^*\|} \right \|   \left\|  \pmat{
	 \|x_{1,k+1} -
			 \wh{x}_{1}(y_k)\|\\
				\vdots \\
	 \|x_{N,k+1} -
			 \wh{x}_{N}(y_k)\|} \right\|. \notag
\end{split}
 \end{equation}
 	  Thus, by  taking expectations   conditioned on $\mathcal{F}_k $, from  \eqref{inexact-sub}, \eqref{inequ2}, and     \eqref{diameter} we have that
\begin{equation}\label{expect-xestimate}
\begin{split} \mathbb{E} \left[ \|x_{ k+1}- x ^*\|^2  \big | \mathcal{F}_k\right]& \leq a^2 \|x_k-x^*\|^2
+\sum_{i=1}^N  \alpha_{i,k}^2 +2\sqrt{ \sum_{i=1}^N D_{X_i}^2  } \sqrt{ \sum_{i=1}^N \alpha_{i,k}^2  }\quad a.s.
\end{split}
\end{equation}
Since  $ \alpha_{i,k}$ is deterministic and $\lim\limits_{k\to \infty}  \alpha_{i,k}=0$,  by taking expectations on   both sides
of  \eqref{expect-xestimate}  and by  applying Lemma~\ref{supp-conv}(b)  we have that   $\lim\limits_{k \to \infty} \mathbb{E}[ \|x_{ k}- x ^*\|^2] = 0.$  Thus, the result  (a) follows by the expression \eqref{variance-inequ0} .
}

{
(b) By Jensen's inequality, we have that
$   \mathbb{E} \left[    \left\|   x_{i,k } - x_i^*\right\| \right]   \leq
\sqrt{ \mathbb{E} \left[    \left\|   x_{i,k } - x_i^*\right\|^2 \right] }, $
which implies   result  (b) by invoking the fact that  $\lim\limits_{k \to \infty} \mathbb{E}[ \|x_{ i,k }- x ^*\|^2|] = 0~\forall i\in \mathcal{N}.$
 \hfill $\Box$
 }

\subsection{Rate of convergence}
We begin with a supporting lemma that bounds a
sub-linear rate by a prescribed linear rate.  This has been  proven in \cite{ahmadi16analysis} but a simpler proof is provided here.
\begin{lemma}\label{linear-bd} ~
Given a function $zc^z$ with  $0<c<1$. Then for all $z \geq 0$, we have that
$$ zc^z \red{\leq} Dq^z, \mbox{ where }
 c < q < 1 \mbox{ and } D \red{\geq}\frac{1}{\ln(q/c)^e}.$$
\end{lemma}
\noindent {\bf Proof.}
 Since $zc^z  \red{\leq}  Dq^z$ for all $z\geq 0$. If $p \triangleq c/q < 1$, it follows that
$ \max_{z \geq 0} z p^z  \red{\leq}  D$. Then a maximizer of $zp^z$ satisfies the following:
$p^z + z p^z \ln(p) = 0$ implying that $z = -{1}/{\ln(p)} > 0.$
($z$ can be seen to be a maximizer by taking second derivatives).
It follows that $$\max_{z \geq 0} zp^z =  -\frac{1}{\ln(p)} p^{\frac{1}{-\ln(p)}} =  -\frac{1}{\ln(p)} \left(\frac{1}{p}\right)^{\frac{1}{\ln(p)}} = \frac{1}{\ln(1/p)} \red{ e^{-1}} = \frac{1}{\ln(q/c)^e}. \qquad \hfill \Box $$

\begin{proposition}[{\bf Geometric rate of convergence}] \label{prp-linear}~
Consider the synchronous  inexact proximal    BR scheme (Algorithm \ref{inexact-sbr-cont})
 where $\mathbb{E}[\|x_{i,0}-x_i^*\|] \leq C $ and $\alpha_{i,k} =\eta^{k+1} $ for some  $\eta \in (0,1)$.
Suppose   Assumptions \ref{assump-play-prob} and   \ref{assp-gamma} hold.  {Define $c  \triangleq \max\{a,\eta\}$, and
\begin{equation}
\begin{split}
u_k  \triangleq \mathbb{E} \left[\left\|   \pmat{ \|x_{1,k}  -{x}_{1}^*\| \\
				\vdots \\
 \|x_{N,k}   - {x}_{N}^*\|} \right \| \right ] . \label{def-uk}
\end{split}
\end{equation}
Then the following holds for  any $q \in (c,1)$ and  $D \triangleq 1/\ln((q/c)^e)$
\begin{align} \label{geom-err} u_{k} &  \leq  \sqrt{N}(C+D)q^k\quad \forall k \geq 0.
\end{align}}
\end{proposition}

\noindent {\bf Proof.}
{By taking expectations on both sides of  \eqref{ineaxt-contractive},  from \eqref{inequ2} we obtain the following   recursion:
\begin{equation}\label{err-iter}
u_{k+1} \leq a u_k + \sqrt{N} \alpha_{\max,k} \quad \forall k \geq 0.
\end{equation}
Since  $\alpha_{\max,k}= \eta^{k+1}$,  by \eqref{err-iter} we obtain the following sequence of inequalities:
\begin{equation}\label{seq-bd-uk}
\begin{split}
	u_{k} & \leq a u_{k-1} + \sqrt{N} \eta^{k}    \leq a^k u_0 + \sqrt{N} \sum_{j=1}^{k} a^{k-j} \eta^j  \leq   c^k u_0+ \sqrt{N} \left( \sum_{j=1}^k c^k \right)
		 = (u_0+ \sqrt{N} k) c^k.
		 \end{split}
\end{equation} By Lemma~\ref{linear-bd},
 we have that   $k c^k  \leq  D q^k  $ for any $q \in (c,1)$ and   $D \geq 1/\ln((q/c)^e)$. Then  by \eqref{seq-bd-uk}, we get
\begin{align*}
 (u_0+ \sqrt{N} k) c^k & \leq  (u_0 + \sqrt{N}D) q^k.
\end{align*}}
Since  $u_0 \leq \sqrt{N} C$ by $\mathbb{E}[\|x_{i,0}-x_i^*\|] \leq C~~\forall i \in \mathcal{ N}$,  it follows  that
$ u_{k} \leq \sqrt{N}(C+D) q^k. $
\hfill $\Box$

{ \begin{remark} ~  This  result shows that if the inexactness sequence is driven to zero at a suitable rate, there is {\bf no}
degradation in the overall rate of convergence; in effect, the rate stays
linear and can be precisely specified. Notably, the impact of the
number of players on the rate appears only in terms of the constant.
 While there exists prior research  where  a  geometric rate of convergence is retained by solving
subproblems accurately enough  (cf. \cite{friedlander2012hybrid,schmidt2011convergence}),   this result appears to be the first such result in the context of inexact proximal best-response schemes
for stochastic Nash games.
\end{remark}}

\subsection{Overall iteration complexity analysis} \label{sec:syn:complexity}
We now utilize the results from the prior subsection to develop overall
iteration complexity bounds. Recall that under the prescribed setting,
 the regularized variant, as prescribed by \eqref{def-wh_x_inft},  is always strongly convex with parameter $\mu$.
Given the strong convexity of each player's objective function, SA   provides an avenue for obtaining  an inexact solution.
In such a scheme, at the major iteration $k$,  given
$y_k$ the $i$th player takes the following sequence of {\em stochastic} gradient steps from $t = 1, \hdots, j_{i,k}$,  where $j_{i,k}$ steps suffice in obtaining an inexact	solution of accuracy $\alpha_{i,k}$:
\begin{align} \tag{SA$_{i,k}$}
z_{i,t+1} & := \Pi_{X_i} \left[ z_{i,t} -  		 \gamma_{t} \underbrace{ \left(\nabla_{x_i} \label{sa-inner}
 \psi_i(z_{i,t},y_{-i,k};\xi_{i,k}^t) + \mu (z_{i,t}-y_{i,k})\right) }_{\triangleq G_i(z_{i,t},y_k,\xi_{i,k}^t)}\right],
\end{align}
where  $z_{i,1}=x_{i,k}$, and $ \gamma_{t}=1/\mu(t+1)$. Note that $ \nabla_{x_i} \psi_i(z_{i,t},y_{-i,k};\xi_{i,k}^t)$ denotes the \red{sampled}  gradient  of $f_i(\cdot)$ at point  $(z_{i,t},x_{-i,k})$ while $ \gamma_{t}$ is  a non-summable but square summable sequence.

Let the SA scheme  \eqref{sa-inner}  from   $t = 1, \hdots, j_{i,k}$ be  employed for obtaining $x_{i,k+1}$.
Define $ \xi_{i,k}\triangleq (  \xi_{i,k}^1,\cdots,  \xi_{i,k}^{j_{i,k}}),$ {
  $ \xi_{i,k}^{[t]}=( \xi_{i,k}^1,\cdots, \xi_{i,k}^t)$},
and  $\mathcal{F}_k \triangleq \sigma\{x_0, \xi_{i,l}, i \in \mathcal{N}, 0 \leq l \leq k-1\}.$
  Then by Algorithm \ref{inexact-sbr-cont},  we see that $x_k$  is adapted to $\mathcal{F}_k$, and
     $\wh{x}_{i} (y_{k} )$ is also adapted to $\mathcal{F}_k$  since $\wh{x}_{i} (y_{k} )=\wh{x}_{i} (x_{k} )$ is an  optimal solution  to  the parameterized  problem \eqref{def-wh_x_inft}.   \red{We derive an error bound for \eqref{sa-inner} as specified by the following lemma.}

{\begin{lemma}\label{lem-rate-sa} ~Let Assumption~\ref{assump-play-prob} hold.
Consider the synchronous inexact proximal BR scheme given by   Algorithm \ref{inexact-sbr-cont}.
Assume  that  \red{for any $t=1,\cdots, j_{i,k},$} $  \mathbb{E}\left[  \psi_i(z_{i,t},y_{-i,k};\xi_{i,k}^t)  \big | \sigma\{\mathcal{F}_k, \xi_{i,k}^{[t-1]}\}\right]  =\nabla_{x_i}
f_i(z_{i,t},y_{-i,k})   ~a.s.,$ and   $  \mathbb{E}\left[ \|  \psi_i(z_{i,t},y_{-i,k};\xi_{i,k}^t) \|^2 \big | \sigma\{\mathcal{F}_k, \xi_{i,k}^{[t-1]}\}\right]  =
\| \nabla_{x_i}f_i(z_{i,t},y_{-i,k})  \|^2 ~a.s.$   Define  $Q_i \triangleq  \frac{2 M_i^2}{ \mu^2}+2D_{X_i}^2.$ Then
\red{the following holds for any $t=1,\cdots, j_{i,k}:$}
 $$ \mathbb{E}\left[\|z_{i,t}-\wh{x}_i(y_k)\|^2\big | \mathcal{F}_k\right] \leq   {Q}_i/{(t+1)}\quad a.s.  $$
\end{lemma}
  \noindent{\bf Proof.}
Set $A_t= \|z_{i,t}-\wh{x}_i(y_k)\|^2~{\rm and}~a_t=\mathbb{E}\left[ A_t\big | \mathcal{F}_k \right].$
Since $X_i$ is a convex set and $\wh{x}_i(y_k) \in X_i$,  from algorithm  \eqref{sa-inner} by the nonexpansive property of the  projection operator  we have
\begin{equation}\label{quard-dist}
  \begin{split}
    A_{t+1}  \leq A_t+   \gamma_{t}^2\|  G_i(z_{i,t},y_k, \xi_{i,k}^t)\|^2-2 \gamma_{t}  \left(z_{i,t}-\wh{x}_i(y_k)\right)^T G_i(z_{i,t},y_k, \xi_{i,k}^t).
\end{split}
\end{equation}
Since  $x_k$  is adapted to $\mathcal{F}_k$, by \eqref{sa-inner} it is seen that $z_{i,t}$ is adapted to $\sigma\{\mathcal{F}_k, \xi_{i,k}^{[t-1]}\}$.
  Then   by invoking $  \mathbb{E}\left[  \psi_i(z_{i,t},y_{-i,k};\xi_{i,k}^t)  \big | \sigma\{\mathcal{F}_k, \xi_{i,k}^{[t-1]}\}\right]  =\nabla_{x_i}
f_i(z_{i,t},y_{-i,k})   ~a.s.,$ and  the tower property of conditional expectation we have that for any $t =1,\cdots,j_{i,k}$:
\begin{align}
  & \mathbb{E}\left[ \left(z_{i,t}-\wh{x}_i(y_k)\right)^T G_i(z_{i,t},y_k, \xi_{i,k}^t)  \big| \mathcal{F}_k\right]  =
   \mathbb{E}\left[   \mathbb{E}\left[  \left(z_{i,t}-\wh{x}_i(y_k)\right)^T G_i(z_{i,t},y_k, \xi_{i,k}^t)\big | \sigma\{\mathcal{F}_k, \xi_{i,k}^{[t-1]}\}\right] \Big| \mathcal{F}_k\right]  \notag \\&   =
   \mathbb{E}\left[  \left(z_{i,t}-\wh{x}_i(y_k)\right)^T\left(\nabla_{x_i}f_i(z_{i,t},y_{-i,k}) + \mu (z_{i,t}-y_{i,k})\right)\Big| \mathcal{F}_k\right] ~~a.s. \label{cond-expect0}
\end{align}

Note that $\tilde{f}_i(x_i)=f_i(x_i,y_{-i,k}) + \frac{\mu}{2} \|x_i-y_{i,k}\|^2$ is strongly convex with parameter $\mu$ and $\wh{x}_i(y_k)$ denotes the optimal solution of $\min_{x_i \in X_i}\tilde{f}_i(x_i)$. Then  by $z_{i,t} \in X_i$,
we  may conclude that for any $t =1,\cdots,j_{i,k}$:
\begin{align*}
    &\left(z_{i,t}-\wh{x}_i(y_k)\right)^T \nabla \tilde{f}_i(\wh{x}_i(y_k))  \geq 0,  \\
    {~~\rm and~~}& \left(z_{i,t}-\wh{x}_i(y_k)\right)^T \left(\nabla  \tilde{f}_i(z_{i,t})- \nabla \tilde{f}_i(\wh{x}_i(y_k))\right)  \geq  \mu\|z_{i,t}-\wh{x}_i(y_k)\|^2= \mu A_t.
\end{align*}
Consequently, by adding the two inequalities we obtain the following:
$$   \left(z_{i,t}-\wh{x}_i(y_k)\right)^T  \left ( \nabla_{x_i}
f_i(z_{i,t},y_{-i,k}) + \mu (z_{i,t}-y_{i,k})\right)=\left(z_{i,t}-\wh{x}_i(y_k)\right)^T \nabla  \tilde{f}_i(z_{i,t})
 \geq  \mu A_t.$$
After taking conditional expectations and invoking  \eqref{cond-expect0}, we obtain
\begin{equation}\label{cond-expect2}
  \begin{split}
  \mathbb{E}\left[ \left(z_{i,t}-\wh{x}_i(y_k)\right)^T G_i(z_{i,t},y_k, \xi_{i,k}^t)  \big| \mathcal{F}_k\right]  \geq
  \mu \mathbb{E}\left[ A_t \big| \mathcal{F}_k\right]=    \mu a_t~~a.s .
\end{split}
\end{equation}

Since   $   \|G_i(z_{i,t},y_k, \xi_{i,k}^t)\|^2 \leq 2  \|\nabla_{x_i}  \psi_i(z_{i,t},x_{-i,k}; \xi_{i,k}^t) \|^2 + 2 \| \mu (z_{i,t}-y_{i,k})\|^2 ,$
 by  Assumption~\ref{assump-play-prob}(d) and \eqref{diameter},  and by
  invoking  $  \mathbb{E}\left[ \|  \psi_i(z_{i,t},y_{-i,k};\xi_{i,k}^t) \|^2 \big | \sigma\{\mathcal{F}_k, \xi_{i,k}^{[t-1]}\}\right]  =
\| \nabla_{x_i}f_i(z_{i,t},y_{-i,k})  \|^2 ~a.s.,$
 we have that
  $$   \mathbb{E}\left[\|  G_i(z_{i,t},y_k, \xi_{i,k}^t)\|^2 \big | \sigma\{\mathcal{F}_k, \xi_{i,k}^{[t-1]}\}\right]    \leq 2 M_i^2 + 2 \| \mu (z_{i,t}-y_{i,k})\|^2\leq 2( M_i^2+\mu^2D_{X_i}^2)~~a.s .$$
   Then  by the tower property of conditional expectations  we obtain the following:
   \begin{equation}\label{cond-expect1}
\mathbb{E}\left[ \|  G_i(z_{i,t},y_k, \xi_{i,k}^t)\|^2 \big| \mathcal{F}_k\right] \leq   2( M_i^2+\mu^2D_{X_i}^2)~~a.s.
\end{equation}
By taking  expectations conditioned   on   $\mathcal{F}_k$  on both sides of  \eqref{quard-dist},  from
 \eqref{cond-expect2} and  \eqref{cond-expect1} we have
$$a_{t+1}\leq (1-2\mu \gamma_{t} ) a_t+ 2 ( M_i^2+\mu^2D_{X_i}^2) \gamma_{t}^2 ~~a.s .$$
By substituting $ \gamma_{t}=\frac{1}{\mu (t+1)}$   we obtain
 $a_{t+1}\leq (1- \frac{2}{t+1} ) a_t+\frac{Q_i}{(t+1)^2}~~a.s., $ where $Q_i \triangleq  \frac{2 M_i^2}{ \mu^2}+2D_{X_i}^2.$

 We   inductively  prove that    $a_t \leq \frac{Q_i}{(t+1)}~a.s. $
  By  \eqref{diameter} we have that  $A_1 \leq D_{X_i}^2,$ and hence $a_1 \leq {Q_i \over 2}$.
 Suppose that   $a_t \leq \frac{Q_i}{(t+1)}~a.s. $ holds for some $t \geq 1$.  It remains to show that
 $a_{t+1} \leq \frac{Q_i}{(t+2)}~a.s. $ Note that
 \begin{align*}
  a_{t+1}&\leq   \left(1- \frac{2}{t+1} \right) \frac{Q_i}{(t+1)}+\frac{Q_i}{(t+1)^2}
  =Q_i\frac{t}{(t+1)^2} \leq Q_i\frac{1}{t+2}~~a.s.,
\end{align*}
where the last  inequality holds by the fact that  $ t(t+2)  \leq (t+1)^2.$
Consequently, we conclude that  $a_t \leq \frac{Q_i}{(t+1)}~a.s. $, and hence the lemma follows.
 \hfill $\Box$
}

 To address strategies  that are an approximate Nash equilibrium, we use the concept of $\epsilon-$Nash equilibrium.
  A random strategy profile $x: { \Omega} \to \Real^{n}$ is an $\epsilon-$NE$_2$ if
 \begin{align}\label{def-epsilon-NE1}
    &  \mathbb{E}  \left\|   \pmat{ \|x_{1} -{x}_{1}^*\| \\
				\vdots \\
 \|x_{N}  - {x}_{N}^*\|} \right \|  \leq \epsilon.
\end{align}
Next, we derive a bound on the overall iteration
complexity for  the synchronous  inexact BR scheme.

\begin{theorem}\label{prop-iter-comp}~
Let   Algorithm \ref{inexact-sbr-cont}  be applied to the $N-$player stochastic Nash game (SNash),
 where $\mathbb{E}[\|x_{i,0}-x_i^*\|] \leq C $ and $\alpha_{i,k}=\eta^{k+1}$ for some $ \eta \in (0,1)$.
    Assume that  an inexact solution to the proximal   BR problem given by \eqref{inexact-sub} is
computed via  \eqref{sa-inner}.   Suppose  Assumptions \ref{assump-play-prob} and   \ref{assp-gamma} hold.
Then the following hold.

(a) {\bf (Overall iteration complexity)}
 Define $c\triangleq  \max\{a,\eta\}$, let $ q \in (c,1)$,  and $D= 1/\ln((q/c)^e)$.
Then the number of projected stochastic gradient steps for player $i$ to compute an  $\epsilon-$NE$_2$ is
no greater than
\begin{align}\label{comp-bound}
 \ell_i(\eta) \triangleq \frac{Q_i}{  \eta^{4}\ln(1/\eta^2)}
			\left(\frac{{\sqrt{N}(C+D)}}{\epsilon}\right)^{\frac{\ln(1/\eta^2)}{\ln(1/ q)}}+ \left \lceil
	\frac{\ln\left(\sqrt{N}(C+D)/\epsilon\right)}{\ln(1/ q)}  \right \rceil.
\end{align}

(b) {\bf (Bounds on complexity)}  The expression \eqref{comp-bound}  satisfies   $\ell_i(\eta)  ={\Omega}(1/\epsilon^2)$.
 By setting $\eta = a$ and  $D= {1+2\delta^{-1}\over e \ln(a^{-1})},$
 $\ell_i(a)$ satisfies the following   for any   given  $\delta > 0$:
 \begin{equation}\label{upper-bound}
\ell_i(a) \leq \frac{Q_i}{a^{4}\ln(1/a^2)}
			\left(\frac{{\sqrt{N}(C+D)}}{\epsilon}\right)^{2+ \delta}+  \left \lceil
	\ln\left(\frac{\sqrt{N}(C+D) }{\epsilon}\right) \frac{2+\delta}{ 2 \ln(a^{-1})} \right \rceil.
\end{equation}
\end{theorem}
\noindent {\bf Proof.}
(a) Suppose  $j_{i,k}=\left \lceil \frac{Q_i}{ \eta^{2(k+1)} }\right\rceil$ steps of  the SA scheme (SA$_{i,k}$)
are taken at major iteration $k$ 	 to obtain $x_{i,k+1}$, i.e.,  $x_{i,k+1}=z_{i,j_{i,k}}.$
  Then by Lemma \ref{lem-rate-sa} and $j_{i,k} \geq \frac{Q_i}{ \eta^{2(k+1)}}$, we have that
	 $$ \mathbb{E}\left[\|x_{i,k+1}-\wh{x}_i(y_k)\|^2\big | \mathcal{F}_k\right]
	  = \mathbb{E} \left[\|z_{i,j_{i,k}}-\wh{x}_i(y_k)\|^2\big | \mathcal{F}_k\right] \leq
	\frac{ Q_i }{j_{i,k}} \leq\eta^{2(k+1)}=\alpha_{i,k}^2 \quad \forall i \in \mathcal{N}.$$
Based on \eqref{geom-err}, we have that
 \begin{align}  u_{k+1} \leq \sqrt{N}(C+D) q^{k+1} \leq \epsilon
		  & \Rightarrow \ q^{k+1} \leq
\frac{\epsilon}{\sqrt{N}(C+D)} \triangleq \epsilonbar \label{def-epsbar}   \Rightarrow \ k \geq \frac{
\ln( 1/\epsilonbar)}{\ln(1/q)} -1.  				\end{align}
By recalling that for $\beta > 1$, the following holds:
\begin{align}\label{exponent-sum}
\sum_{k=1}^K \beta^k &  \leq \int_1^{K+1} \beta^x dx \leq
\frac{\beta^{K+1}} {\ln(\beta)},
\end{align}
{which allows us to bound the overall iteration complexity of player $i$ as  follows:
\begin{align*}   \sum_{k=0}^{\left\lceil
	\frac{\ln(1/\epsilonbar)}{\ln(1/ q)}\right\rceil-1}
  j_{i,k}  & \leq   \sum_{k=0}^{\left\lceil
	\frac{\ln(1/\epsilonbar)}{\ln(1/ q)}\right\rceil-1}  \left(1+\frac{Q_i}{\eta^{2(k+1)}}\right)= \sum_{k=1}^{\left\lceil
	\frac{\ln(1/\epsilonbar)}{\ln(1/ q)}\right\rceil}
  \frac{Q_i}{\eta^{2k}}+ \left\lceil
	\frac{\ln(1/\epsilonbar)}{\ln(1/ q)}\right \rceil	  \leq  \frac{Q_i }{\ln(1/\eta^2)}
			 \eta^{ -2\left( {\ln(1/\epsilonbar) \over \ln(1/ q)} +2\right)} + \left \lceil
	\frac{\ln(1/\epsilonbar)}{\ln(1/ q)}\right\rceil .
\end{align*}
Note that
 \begin{equation}\label{change-exp}
\eta^{ -2 \frac{\ln(1/\epsilonbar)}{\ln(1/ q)}  }
	  = \left(e^{\ln(\eta^{-2})}\right)^{\frac{\ln(1/\epsilonbar)}{\ln(1/ q)}}
			  =  e^{\ln(1/\epsilonbar))^{\frac{\ln(\eta^{-2})}{\ln(1/ q)}}}
			  	=   {(1/\epsilonbar)^{\frac{\ln(1/\eta^2)}{\ln(1/ q)}}}  .
\end{equation}}
Therefore,  the overall iteration complexity of player $i$  is bounded by the following:
 \begin{align*}  \frac{Q_i}{\eta^{4}\ln(1/\eta^2)}
			\left(\frac{1}{\epsilonbar}\right)^{\frac{\ln(1/\eta^2)}{\ln(1/ q)}}+ \left \lceil
	\frac{\ln(1/\epsilonbar)}{\ln(1/ q)}  \right \rceil .
\end{align*}
Then the result (a) follows  by  the definition of $\epsilonbar$ in \eqref{def-epsbar}.

\noindent  {\bf (b).} {The proof of the first  assertion  $\ell_i(\eta)  ={\Omega}(1/\epsilon^2)$
can be found in \cite[Proposition 4(b)]{shanbhag16inexact}.}

  Define $ q= c e^{\delta_0/2}$, where $\delta_0=\frac{\delta \ln(a^{-1})}{1+\delta/2}$.
 Since $c=a$ by $ \eta=a$, we obtain the following:
\begin{align*}
\ln(1/ q)& = \ln( a^{-1})- \delta_0/2=\frac{ 2 \ln(a^{-1})}{2+\delta},
 \frac{\ln(1/\eta^2)}{\ln(1/ q)}       = \frac{2\ln(a^{-1})}{\ln(1/ q)}= 2+ {\delta_0 \over  \frac{  \ln( a^{-1})}{1+\delta/2}}=2+\delta,\\
\textrm{ and~ ~}D&= 1/\ln((q/c)^e)=1/\ln(e^{\delta_0 e/2})={2 \over \delta_0 e}={1+2\delta^{-1}\over e \ln(a^{-1})}.
\end{align*}
Consequently, by the  expression $\ell_i(\eta)$ defined in \eqref{comp-bound},
 we   obtain \eqref{upper-bound}   when $\eta$ is set as $a$.
\hfill $\Box$

\begin{remark} ~Several points need reinforcement in light of
the obtained complexity statements: \\
\noindent {\bf (i)} First, the iteration complexity bound grows slowly
 in the number of players, a desirable feature
of any distributed algorithm employed for games with a
large collection of agents. \\
\noindent {\bf (ii)} Second, when $N = 1$, the complexity bound reduces to ${\cal O}(1/\epsilon^{2+\delta})$  which
can be made arbitrarily close to the bound
	for standard stochastic convex programs suggesting that  the bounds may well be optimal. \\
\noindent {\bf (iii)} Third, in an effort to examine the tightness of the bound,
	we show that our bound on the derived iteration complexity is bounded from
	below by a term that is of the order ${\cal O}(N/\epsilon^2)$.
	Note this is not a lower bound on the complexity itself which we leave as future research.
\end{remark}

 \begin{corollary} ~\label{cor-bound}
Let   Algorithm \ref{inexact-sbr-cont}  be applied to the
$N-$player stochastic Nash game (SNash), where $\alpha_{i,k}=\eta^{k+1}$ for some $ \eta \in (a,1)$,
$\mathbb{E}[\|x_{i,0}-x_i^*\|] \leq C,$
and an  inexact  proximal  BR solution satisfying \eqref{inexact-sub} is computed via   \eqref{sa-inner}.
 Suppose  Assumptions \ref{assump-play-prob} and   \ref{assp-gamma} hold.
 Then the number of projected gradient steps  for player $i$ to compute an   $\epsilon-$NE$_2$  is no greater than
\begin{align}\label{cor-comp-bound}\frac{Q_i}{  \eta^{4}\ln(1/\eta^2)}
\left(\frac{{\sqrt{N}\left(C+ {\eta \over \eta-a}\right)}}{\epsilon}\right)^2
+\left\lceil \frac{\ln\left( {\sqrt{N} \over \epsilon}\left(C+ {\eta \over \eta-a}\right)\right)}{\ln(1/\eta)}\right \rceil.
\end{align}
 \end{corollary}
 { \bf Proof.}
{Since ${a \over \eta} \in (0,1) $ by $ \eta \in (a,1)$, we have that
\begin{align*}
   \sum_{j=1}^{k} a^{k-j} \eta^j &   =\eta^k  \sum_{j=0}^{k-1} (a/\eta)^{j}
		 \leq {1 -(a/\eta)^k  \over 1-a/\eta}\eta^k= {\eta^k -a^k  \over 1-a/\eta}.
\end{align*}
Then by  $u_0 \leq \sqrt{N} C$,  $a<\eta$, and by invoking the second  inequality in \eqref{seq-bd-uk},
 we obtain the following bound:
\begin{align*}
 &   u_k \leq  \sqrt{N} C a^k + {\sqrt{N} \eta  \over \eta-a } (\eta^k-a^k) \leq  \left(\sqrt{N} C  + {\sqrt{N} \eta  \over \eta-a }\right)\eta^k.
\end{align*}
Thus,  by setting  $\bar{k}=\left \lceil \frac{\ln(1/\bar{\epsilon})}{\ln(1/\eta)}\right \rceil $ with  $\bar{\epsilon}={\epsilon \over \sqrt{N} (C  + \eta /(\eta-a ))}$,
we get $u_{k+1} \leq   \epsilon~\forall   k\geq\bar{k}-1 .$
	Suppose  $j_{i,k}=\left \lceil \frac{Q_i}{ \eta^{2(k+1)} }\right \rceil$   steps  of (SA$_{i,k}$)
are taken at major iteration $k$ to obtain   $x_{i,k+1} .$
  Then by Lemma \ref{lem-rate-sa},   we have that
	 $ \mathbb{E}\left[\|x_{i,k+1}-\wh{x}_i(y_k)\|^2\big | \mathcal{F}_k\right] \leq \alpha_{i,k}^2$ for all $i \in \mathcal{N} .$ By invoking   \eqref{exponent-sum}  for $\beta > 1$,
we may bound the  iteration complexity  of player $i$  as  follows:
\begin{align*}   \sum_{k=0}^{\left\lceil
	\frac{\ln(1/\bar{\epsilon})}{\ln(1/\eta)}\right\rceil-1}  j_{i,k}  &=   \sum_{k=1}^{\left \lceil
	\frac{\ln(1/\bar{\epsilon})}{\ln(1/\eta)}\right\rceil}  \frac{Q_i}{\eta^{2k}}  +\left\lceil \frac{\ln(1/\bar{\epsilon})}{\ln(1/\eta)}\right \rceil
  \leq  \frac{Q_i }{\ln(1/\eta^2)}   \eta^{ -2\left(\frac{\ln(1/\bar{\epsilon})}{\ln(1/\eta)} +2\right)} +\left\lceil \frac{\ln(1/\bar{\epsilon})}{\ln(1/\eta)}\right \rceil
\end{align*}
which results in  \eqref{cor-comp-bound}  by  utilizing    \eqref{change-exp} and by  recalling  $\bar{\epsilon}={\epsilon \over \sqrt{N} (C  + \eta /(\eta-a ))}$.}
	\hfill $\Box$

 \begin{remark}~
It seems that Corollary \ref{cor-bound}  gives   better results than those  of Theorem \ref{prop-iter-comp} since it shows  that  the iteration complexity is exactly of  $ {\cal O}(1/\epsilon^2)$  instead of  ${\cal O}(1/\epsilon^{2+\delta})$.   However,   Corollary \ref{cor-bound}   restricts the parameter $\eta$ to be  $\eta \in (a,1)$ while
  Theorem \ref{prop-iter-comp}  provides   more flexibility on  the selection of  $\eta$ only by requiring $\eta \in (0,1).$   As a matter of fact,
the numerical  results demonstrated in Table \ref{TAB1}  indicate that  smaller $\eta \in (0,1)$ may give better  empirical iteration complexity.
In the  following  two sections,  we establish the    iteration complexity for the randomized  and the asynchronous algorithms  through   the same  path of Theorem \ref{prop-iter-comp}.
\end{remark}

 We now provide a related result for a probabilistically $\epsilon$-Nash equilibrium. Specifically, a random strategy profile $x: { \Omega} \to \Real^{n}$ is an
 \red{$\epsilon_{\mathbb{P}\delta}-$NE$_2$} if
\red{							 \begin{align}
    \mathbb{P} \left( \omega:    \left\|   \pmat{ \|x_{1}(\omega) -{x}_{1}^*\| \\
				\vdots \\
 \|x_{N}(\omega)  - {x}_{N}^*\|} \right \| \leq \epsilon\right) \geq (1-\delta). \notag
\end{align}
}

{\begin{corollary} ~Let   Algorithm \ref{inexact-sbr-cont}  be applied to the
$N-$player stochastic Nash game (SNash), where $\alpha_{i,k}=\eta^{k+1}$ for some $ \eta \in (a,1)$,
$\mathbb{E}[\|x_{i,0}-x_i^*\|] \leq C,$
and an  inexact  proximal  BR solution satisfying \eqref{inexact-sub} is computed via   \eqref{sa-inner}.
 Suppose  Assumptions \ref{assump-play-prob} and   \ref{assp-gamma} hold.
 Then the iteration complexity for  player $i$  to compute an  \red{$\epsilon_{\mathbb{P}\delta}-$NE$_2$} is bounded   by \eqref{upper-bound} with  $\epsilon$ replaced by $\epsilon \delta.$
\end{corollary}
{{\bf Note.} This  follows  because $u_k\leq \epsilon \delta $ implies that $\mathbb{P}(v_k\geq \epsilon)
   \leq    {\mathbb{E}[v_k]  / \epsilon} =    {u_k  / \epsilon}  \leq \delta$ by the Markov inequality, \red{
   where  $v_k$ and $u_k$ are defined by \eqref{ineaxt-contractive} and \eqref{def-uk}, respectively.}
   Then the $ \epsilon \delta $-NE$_2$ $x$ is  also   an \red{$\epsilon_{\mathbb{P}\delta}-$NE$_2$} solution. 
}

\section{A Randomized  Inexact Proximal  BR Algorithm } \label{sec:random}
In this  section, we propose a randomized inexact proximal  BR  algorithm, where a subset of players
   is randomly  chosen to update their strategies  in each major  iteration,
 and the inexact    proximal   BR solutions are achieved via   the SA scheme.

 \subsection{Algorithm description}
 The    randomized block coordinate  descent method was   proposed in  \cite{nesterov2012efficiency} to solve large-scale optimization problems,
 where     the coordinates are partitioned    into several  blocks and  at each iteration only one block of variables is randomly chosen to  update
 while  the other  blocks are kept  invariant. Subsequently,   this
convergence was   analyzed  to nonconvex and fixed-point
regimes~\cite{nesterov2012efficiency,richtarik2014iteration,lu2013randomized,combettes2015stochastic}.
  Motivated by  these  research,    we   design a  randomized  inexact  proximal  BR scheme as follows:
 For any $i \in \mathcal{N}$, let $\{\chi_{i,k}\}_{k\geq 0}$ be a sequence of  i.i.d Bernoulli random variables taking values in $\{0,1\}$.
The variable $\chi_{i,k}$ signals whether the $i$th player updates at major iteration $k $:
if  $ \chi_{i,k}=1$,  then player $i$ initiates an update at major iteration $k$ and computes  an   inexact proximal    BR solution satisfying \eqref{rand-inexact-sub2},      where  $\mathcal{F}_k$ is the $\sigma$-field of the entire information  used by  the algorithm up to  (and including)  the update   of  $x_k,$ and $\alpha_{i,k}$ is a nonnegative random variable  adapted  to  $\mathcal{F}_k$. It is worth noticing that $\mathcal{F}_k$  contains  the  sample $\chi_{i,l}$ for any $i \in \mathcal{N} $ and $l:0\leq l\leq k-1,$ and $\alpha_{i,k}$ may also depend on  the  samples $\{\chi_{i,l}\}_{0 \leq l \leq k-1}$.
We   define $\mathcal{F}_k$ in  Section \ref{sec:rand:complexity}, and specify the  selection of the  sequence
  $\{\alpha_{i,k}\}_{k\geq 0}$ when we proceed to   investigate  the convergence properties for which we need  the following condition in the convergence analysis.

\begin{algorithm}[H]
\caption{{Randomized inexact proximal  BR scheme }}
\label{rand-inexact-sbr-cont}
 Let $k:=0$, $ y_{i,0}=x_{i,0} \in X_i$  for $i = 1, \hdots, N$.
\begin{enumerate}
\item[(1)] If  $ \chi_{i,k}=1$,   then    player $i$ updates $x_{i,k+1}$ satisfying the following:  \begin{align} \label{rand-inexact-sub2}
x_{i,k+1}  \in  \left\{ z\in X_i : \mathbb{E}\left[\|z -
\wh{x}_{i}(y_k)\|^2\big |\mathcal{F}_k\right]   \leq \alpha_{i,k}^2~~a.s. \right\},
\end{align}
where $\wh{x}_{i}(y_k)$ is defined  in \eqref{def-wh_x_inft}.  Otherwise,   $x_{i, k+1}=x_{i, k}$. 
\item[(2)] For $i = 1, \hdots, N$,  $y_{i,k+1} :=x_{i,k+1}$;
\item [(3)] $k:=k+1$; If $k < K$, return to (1); else STOP.
\end{enumerate}
\end{algorithm}

\begin{assumption}~ \label{assp-rand}
	 For any $i \in \cal{N}$, $\mathbb{P}(\chi_{i,k}=1)=p_i>0$ and $\chi_{i,k}$ is independent of $\mathcal{F}_k$.
			 \end{assumption}

\begin{remark}~\label{rem-poisson}  We make  the following clarifications on Assumption \ref{assp-rand}:

\noindent {\bf (i)} The condition $\mathbb{P}(\chi_{i,k}=1)=p_i>0$ guarantees that each player initiates an update
with positive  probability at major  iteration $k$ of Algorithm \ref{rand-inexact-sbr-cont}. It  accommodates the
special case where only one player is  randomly chosen with positive probility to update in each iteration.

\noindent {\bf (ii)}  The Poisson model \cite{aysal2009broadcast,boyd2006randomized} is   a special case:
Each player $i \in \cal{N}$ is activated  according to   a local  Poisson clock,  which  ticks according to a Poisson process with rate $\varrho_i>0$.    Suppose that there is a virtual global clock which ticks whenever any of the local Poisson clocks tick.
Suppose that the local  Poisson clocks are independent, then  the global clock ticks according to a
Poisson process with rate $\sum_{i=1}^N\varrho_i$.  Let  $Z_k$ denote the time of the $k$-th tick of the global clock,
and $I_k \in \{1,\cdots,N\}$ denote the set of players   whose clocks tick at time $Z_k.$
Since the  local Poisson clocks  are independent,  with probability one,  $I_k$  contains a single  element and $\mathbb{P}(I_k =i)={\varrho_i \over \sum_{i=1}^N\varrho_i} \triangleq p_i.$
 Besides, the memoryless property of the Poisson process indicates  that  the process  $\{I_k \}_{k\geq 0}$ is  i.i.d.  As a result,  the processes $\{\chi_{i,k}\}~\forall i \in\cal{N} $ are mutually independent,  and  for each $i \in \mathcal{N}$,
 $\{\chi_{i,k}\} $  is an i.i.d  sequence with $\mathbb{P}(\chi_{i,k}=1)=p_i>0$ and $\sum_{i=1}^N p_i=1.$
\end{remark}

\subsection{Convergence analysis}

{We now establish the almost sure convergence and the  geometric rate of convergence of the $\{x_k\} $
to $x^*$ under suitable conditions on the inexactness  sequences $\{\alpha_{i,k}\}, i\in \mathcal{N}.$}
\begin{lemma} ~\label{lem-as-asy}({\bf a.s. convergence})
Let  $\{x_k\}_{k \geq 0}$ be generated by Algorithm~\ref{rand-inexact-sbr-cont}.  Suppose
 Assumptions \ref{assump-play-prob},   \ref{assp-gamma},  and  \ref{assp-rand} hold.
For any $i \in \cal{N}$,  $0 \leq \alpha_{i,k} <1$ and  $\sum_{k=0}^{\infty}\alpha_{i,k}<\infty ~a.s$.
 Then for  any $i \in \mathcal{N} $,
  $ \lim\limits_{ k \to\infty}  x_{i,k} =x_i^*~~a.s.$
\end{lemma}
{\bf Proof.}   See Appendix \ref{app-lem-4}. \hfill $\Box$

Define  $\beta_{i,0}=0 $ and  $\beta_{i,k}=\sum_{p=0}^{k-1}\chi_{i,p} ~{\rm for~ all ~}k \geq 1$.
Thus, $\beta_{i,k}$ is adapted  to  $\mathcal{F}_k$.   Note that for the Poisson model described in Remark \ref{rem-poisson},
  players  may not able to know the exact number of ticks of the global clock. Then   it is impractical  to set $\alpha_{i,k}$ be a function of $k$, instead,  a function of $\beta_{i,k}$ is  appropriate.

	\begin{lemma} [Geometric   Convergence] \label{rand-geometric}
	  ~Let $\{x_k\}_{k \geq 0}$ be generated by Algorithm~\ref{rand-inexact-sbr-cont},
	  where  $\mathbb{E}[\|x_{i,0}-x_i^*\|] \leq C$ and $\alpha_{i,k}=\eta^{\beta_{i,k}+1}$ for some $\eta \in (0,1)$.
	Suppose  Assumptions \ref{assump-play-prob},   \ref{assp-gamma}  and    \ref{assp-rand} hold.
	 Define $\tilde{c} \triangleq\max\{\tilde{a},\tilde{\eta}\}$  with  $\tilde{a}$ and $\tilde{\eta}$ defined by
 \eqref{asy-as-a} and \eqref{asy-bound-inexct-sequnece}, respectively. Then  the following holds for  any  $ \tilde{q} \in( \tilde{c},1)$:
 \begin{align*}
  \mathbb{E}\left[\|x_{k}-x^*\|_P\right]     \leq  \sqrt{N}( \widetilde{C} +\tilde{D})\tilde{q}^k\quad \forall  k\geq 0,
\end{align*}
 where  $\| \cdot\|_P$ is defined in \eqref{weighted-norm},
    $D \triangleq 1/\ln ((\tilde{q}/\tilde{c})^e),$ $ \widetilde{C}= C\left(\sum_{i=1}^N N^{-1}p_i^{-1}\right)^{1/2}$,  and  $ \widetilde{D}= D\eta  \tilde{\eta}^{-1} $.
\end{lemma}
{\bf Proof.}  See Appendix \ref{app-lem-5}. \hfill $\Box$

\begin{remark}\label{rem-u}~
 By  $\| \cdot\|_P$ defined in \eqref{weighted-norm}, we get $\|x \|_P^2  \geq {1 \over p_{\max}} \sum_{i=1}^N  \|x_i\|^2,$
where $p_{\max}=\max_{i \in \cal{N}} p_i.$ Thus,
 $$ \|x_{k}-x^*\|_P \geq {1 \over \sqrt{p_{\max}}} \sqrt{\sum_{i=1}^N  \|x_{i,k}-x_i^*\|^2}= {1 \over \sqrt{p_{\max}}}  v_k,$$
 where $v_k$ is defined in \eqref{ineaxt-contractive}. Then  by Lemma \ref{rand-geometric}   and  $u_k$ is defined in \eqref{def-uk}, we have that
$u_k  \leq ( p_{\max} )^{1/2} \mathbb{E}[\|x_{k}-x^*\|_P]  \leq ( N p_{\max} )^{1/2} ( \widetilde{C} +\tilde{D})\tilde{q}^k.$
\end{remark}

\subsection{Overall iteration  complexity} \label{sec:rand:complexity}

We  proceed to estimate the  iteration complexity, where the inexact proximal BR solution  is computed via a SA scheme.
 At the major iteration $k$:  if $ \chi_{i,k}=1$,    then    player $i$ takes    SG steps  (SA$_{i,k}$) from $t = 1, \hdots, j_{i,k}$.
 We define   $\xi_{i,k}$, $I_k$, $\mathcal{F}_k$
 as  $\xi_{i,k}\triangleq ( \xi_{i,k}^1,\cdots, \xi_{i,k}^{j_{i,k}}),$ $I_k\triangleq \{ i\in \mathcal{N}:\chi_{i,k}=1\},$ and  $\mathcal{F}_k
\triangleq \sigma\{x_0,\{\xi_{i,l}\}_{i \in I_k},  \{ \chi_{i,l}\}_{ i \in \mathcal{N}},0 \leq l \leq k-1\}$, respectively.
  Then,  by Algorithm \ref{rand-inexact-sbr-cont} ,   $x_k$ is adapted to   $\mathcal{F}_k$.
  Similar to  the proof  of  Lemma \ref{lem-rate-sa},  we    obtain    the following result  for the  SA scheme   (SA$_{i,k}$).

\begin{lemma}~\label{lem-rate-sa3} Let Assumption~\ref{assump-play-prob} hold.
Assume  that  for any $  i \in I_k$  the random variables $\{\xi_{i,k}^t\}_{1 \leq t \leq j_{i,k}}$ are i.i.d,
and in addition,   that  the random vector    $\xi_{i,k}$
is independent  of $\mathcal{F}_k$.   Then for any $ i \in I_k$ and any $ t \geq 1$ we have that
 $ \mathbb{E}[\|z_{i,t}-\wh{x}_i(y_k)\|^2\big | \mathcal{F}_k] \leq   {Q}_i/{(t+1)}~~a.s. $
 \end{lemma}

 \begin{theorem}~ \label{rand-thm-complexity} ~Let the sequence $\{x_k\}_{k \geq 0}$ be generated by Algorithm~\ref{rand-inexact-sbr-cont},
where $\alpha_{i,k}=\eta^{\beta_{i,k}+1}$ for some $\eta \in (0,1)$ and  $\mathbb{E}[\|x_{i,0}-x_i^*\|] \leq C$.
  Suppose that  an inexact solution   characterized  by \eqref{rand-inexact-sub2}
is computed via    \eqref{sa-inner}.
Let  Assumptions \ref{assump-play-prob},   \ref{assp-gamma}  and  \ref{assp-rand} hold.
Suppose     $\tilde{a}$, $\tilde{\eta}$   and $\tilde{\eta}_0$ are defined by
 \eqref{asy-as-a}, \eqref{asy-bound-inexct-sequnece} and   \eqref{rand-def-eta}.
    Define    $ \tilde{c}  \triangleq\max\{\tilde{a},\tilde{\eta}\}$
    and let $ \tilde{q} \in ( \tilde{c} ,1)$.   Then expectation of the number of projected gradient steps  for player $i$  to compute an
   $\epsilon-$NE$_2$  is
no greater than
\begin{align}\label{randcomp-bound}
\tilde{\ell}_i(\eta) \triangleq   \frac{p_iQ_i  }{\eta^2 \tilde{\eta}_0^2\ln(1/\tilde{\eta}_0^2)}
	\left(\frac{1}{ \tilde{\epsilon} }\right)^{\frac{\ln(1/\widetilde{\eta}_0^2)}{\ln(1/\tilde{q} )}}
	+ \left \lceil \frac{\ln(1/ \tilde{\epsilon})}{\ln(1/\tilde{q})}\right \rceil
, \end{align}
where
$ \tilde{\epsilon} \triangleq \frac{\epsilon}{(N p_{\max})^{1/2}(\widetilde{C}+\widetilde{D})}$ with $ \widetilde{C}= C\left(\sum_{i=1}^N N^{-1}p_i^{-1}\right)^{1/2}$    and  $ \widetilde{D}= D\eta  \tilde{\eta}^{-1} $ with   $D \triangleq 1/\ln ((\tilde{q}/\tilde{c})^e)$. If $\eta=a$, then
given any  $\delta > 0$, $ \tilde{\ell}_i(\eta)$ satisfies the following upper bound:
 \begin{equation}\label{rand-bb-estimate2}
  \tilde{\ell}_i(\eta)  \leq   \frac{p_iQ_i }{\tilde{\eta}^{2} \eta^2\ln(1/\tilde{\eta}^2)}
 \left(\frac{1}{ \tilde{\epsilon} }\right)^{2 \ln(\tilde{\eta}_0^{-1})/\ln(\tilde{\eta}^{-1})+\delta}+
 \left \lceil  \ln(1/ \tilde{\epsilon}) \left({ 1\over \ln(\tilde{\eta}^{-1})}+{\delta \over 2\ln(\tilde{\eta}_0^{-1})}  \right)\right \rceil.
\end{equation}
\end{theorem}

\noindent {\bf Proof.}
Suppose $j_{i,k}=\left \lceil \frac{Q_i}{\eta^{2(\beta_{i,k}+1)} }\right\rceil$ steps of  the SA scheme (SA$_{i,k}$)
are taken at major iteration $k$  to obtain $x_{i,k+1}$, i.e.,  $x_{i,k+1}=z_{i,j_{i,k}}.$
Then by Lemma \ref{lem-rate-sa3}, it follows that
	 $$ \mathbb{E}\left[\|x_{i,k+1}-\wh{x}_i(y_k)\|^2\big | \mathcal{F}_k\right]
	  = \mathbb{E} \left[\|z_{i,j_{i,k}}-\wh{x}_i(y_k)\|^2\big | \mathcal{F}_k\right] \leq
	\frac{ Q_i }{j_{i,k}}  \leq \alpha_{i,k}^2 \quad \forall i \in \mathcal{N}.$$
By Remark \ref{rem-u}, we have that
 \begin{align}  u_{k+1}  & \leq (N p_{\max})^{1/2}(\widetilde{C}+\widetilde{D}) \tilde{q}^{k+1} \leq \epsilon  \Rightarrow \ k \geq \frac{ \ln( 1/\epsilonbar)}{\ln(1/\tilde{q})} -1 ,\notag
				\end{align}
 where $ \tilde{\epsilon} \triangleq \frac{\epsilon}{(N p_{\max})^{1/2}(\widetilde{C}+\widetilde{D})}.$
 Similar to  \eqref{asy-bound-inexct-sequnece}   we have that
 \begin{equation}\label{rand-def-eta}
 \mathbb{E} [\eta^{-2\beta_{i,k}}]   =
  (p_i\eta^{-2}+1-p_i)^k=  \left(p_i(\eta^{-2}-1)+1\right)^k\leq \left(p_{\max}(\eta^{-2}-1)+1\right)^k \triangleq \tilde{\eta}_0^{-2k}~~\forall i \in \cal{N}.
\end{equation}
Then by invoking    \eqref{exponent-sum}    for  $\beta > 1$,
  the expectation of the overall iteration complexity  of player $i$   is   bounded by the following:
\begin{align*} & \mathbb{E}\left[  \sum_{k=0}^{\left \lceil\frac{\ln(1/ \tilde{\epsilon})}{\ln(1/\tilde{q})}\right\rceil-1}  j_{i,k} \chi_{i,k}  \right]  = \sum_{k=0}^{\left\lceil\frac{\ln(1/ \tilde{\epsilon})}{\ln(1/\tilde{q})}\right\rceil-1}
 \mathbb{E}\left[  j_{i,k}  \right]  \mathbb{E}[\chi_{i,k}]
 \quad {\scriptstyle \left(  \mathrm{since~}  \chi_{i,k} {\rm~is ~independent ~of ~} j_{i,k} \right)}
\\&\leq  \sum_{k=0}^{\lceil \frac{\ln(1/ \tilde{\epsilon})}{\ln(1/\tilde{q})}\rceil-1}
\left( \frac{p_iQ_i}{\eta^2\tilde{\eta}_0^{2k}}+1\right)
 =  \frac{\tilde{\eta}_0^2}{\eta^2}\sum_{k=1}^{\lceil
	\frac{\ln(1/ \tilde{\epsilon})}{\ln(1/\tilde{q})}\rceil}
  \frac{p_iQ_i}{\tilde{\eta}_0^{2k}}+ \left \lceil
	\frac{\ln(1/ \tilde{\epsilon})}{\ln(1/\tilde{q})} \right\rceil	
	  \leq \frac{\tilde{\eta}_0^2}{\eta^2} \frac{p_iQ_i }{\ln(1/\tilde{\eta}_0^2)}
			\left(\frac{1}{\tilde{\eta}_0^{ 2\left(\frac{\ln(1/ \tilde{\epsilon})}{\ln(1/\tilde{q})} +2\right)}}\right) +
			\left \lceil \frac{\ln(1/ \tilde{\epsilon})}{\ln(1/\tilde{q})}\right \rceil,\end{align*}
{which results in   \eqref{randcomp-bound} by \eqref{change-exp}.} If  $\eta=a$, then  by $\tilde{a}$ and $\tilde{\eta}$ defined in
 \eqref{asy-as-a} and \eqref{asy-bound-inexct-sequnece}, we have  that $  \tilde{a}=\tilde{\eta}$ and
$ \tilde{c}=  \tilde{a}.$  Define $\tilde{q }=  \tilde{\eta}e^{\delta_0/2}$, where
 $\delta_0=\frac{\delta \ln(\tilde{\eta}^{-1})}{\ln(\tilde{\eta}_0^{-1})/\ln(\tilde{\eta}^{-1})+\delta/2}$. Then we obtain the following:
\begin{align}
\ln(1/\tilde{q} )        &=
    \ln( \tilde{\eta}^{-1})- \delta_0/2=\frac{\ln(\tilde{\eta}_0^{-1}) }{\ln(\tilde{\eta}_0^{-1})/\ln(\tilde{\eta}^{-1})+\delta/2}=
    \frac{1 }{1/\ln(\tilde{\eta}^{-1})+\delta/2\ln(\tilde{\eta}_0^{-1})}, \notag \\
&  \frac{\ln(1/\widetilde{\eta}_0^2)}{\ln(1/\tilde{q} )}    =
    \frac{2 \ln(\tilde{\eta}_0^{-1})}  {\ln(1/\tilde{q} )}  =  \frac{2 \ln(\tilde{\eta}_0^{-1})}{\ln( \tilde{\eta}^{-1}) }+\delta, \notag
   \\& {\rm and ~}D= 1/\ln((\tilde{q}/\tilde{c})^e)=1/\ln(e^{\delta_0 e/2})={2 \over \delta_0 e}={\ln(\tilde{\eta}_0^{-1}) /\ln(\tilde{\eta}^{-1})+2\delta^{-1}\over e \ln( \tilde{\eta}^{-1})}.
 \label{def-d}
\end{align}
Consequently, the result  \eqref{rand-bb-estimate2} follows by  the
expression $\tilde{\ell}_i(\tilde{\eta})$ defined in \eqref{randcomp-bound}.\hfill $\Box$

\begin{remark} ~\label{rand-remark} We make the following illustrations on Theorem \ref{rand-thm-complexity}:

\noindent {\bf (i) }If $p_i=1~\forall i \in \cal{N}$, then $\widetilde{Q}=Q$ and  $p_{\min}=1.$
 Thus,  from  \eqref{asy-as-a},  \eqref{asy-bound-inexct-sequnece}, \eqref{rand-def-eta}
 we see that  $\tilde{a}=a,\tilde{\eta}=\eta, \tilde{\eta}_0=\eta,$
and hence $ \widetilde{C}= C,\widetilde{D}= D$.
As a result,  the results  of Lemma \ref{rand-geometric} and Theorem  \ref{rand-thm-complexity} for the randomized scheme  reduce  to those given in  Proposition \ref{prp-linear} and  Theorem \ref{prop-iter-comp} for the synchronous  scheme, respectively.

\noindent {\bf (ii) }
Note that the  iteration complexity  in Theorem  \ref{rand-thm-complexity} of the randomized algorithm  is described
via the    expected number  of projected gradient steps,  since the  number of gradient steps to get
an inexact solution   \eqref{rand-inexact-sub2} is a random variable dependent on the realization of updates in the associated trajectory.  

\noindent {\bf (iii) }By  definitions  \eqref{asy-bound-inexct-sequnece} and \eqref{rand-def-eta} we
have that $$(\tilde{\eta}/\tilde{\eta}_0)^2=1+p_{\max}(\eta^{-2}-1)p_{\min}\left({1 \over p_{\min}}-{\eta^2 \over p_{\max}}-(1-\eta^2) \right) \geq 1+p_{\max}(1-p_{\min})\eta^{-2}(1-\eta^2)^2 \geq 1,$$ where the equality holds only  if $p_{\min} =1.$
Thus, for the case $p_{\min} <1,$ the complexity bound is greater  than  the bounds derived   in  Proposition \ref{prp-linear} and  Theorem \ref{prop-iter-comp},  a consequence of the  cost of randomization.  
\end{remark}

\section{An Asynchronous Inexact Proximal  BR Algorithm } \label{sec:asy}
Asynchronous methods  date  back to \cite{chazan1969chaotic} when they were
employed for the solution to systems of linear equations.  Subsequently,  they
were  used in optimization  problems,  in  which a partially asynchronous
gradient projection algorithm is proposed in   \cite{bertsekas1989parallel},
while the convergence rate is analyzed  in \cite{tseng1991rate}.
Here,   we     adapt the scheme developed in \cite{bertsekas1989parallel}
to   stochastic Nash games  and propose an   asynchronous inexact proximal  BR
algorithm.  {Recall  that in  \cite{bertsekas1989parallel},  the
$\infty$-norm is utilized in  the rate analysis of the asynchronous schemes
for problems with maximum norm contraction mappings.  By
assuming  that the proximal BR map is contractive in the $\infty$-norm, we
obtain a geometric  rate of convergence for an appropriately chosen
inexactness sequence, establish    the  overall iteration complexity    in
terms of the number of   projected  gradient steps, and analyze the associated
complexity  bound.}

\subsection{Algorithm design}
The synchronous algorithm  designed  in the previous  section requires that all players  update their strategies simultaneously.
In a network with a large  collection of noncooperative  players,    players might not be able to make simultaneous  updates nor may they  have  access to their rivals' latest information.  In this context, we  propose the following  asynchronous scheme.    Let $T=\{0,1,2,\cdots\}$ be a set of epochs  at  which one or more players update their strategies. Denote by  $I_k \in \mathcal{N}$  the set of players which update  their strategies  at time $k.$
  For  any $i \in I_k$, player $i$  may not obtain  its rivals' latest information, instead,   the outdated data
 $y_{k}^i\triangleq (x_{1,k-\tau_{i1}(k)},\cdots,x_{N,k-\tau_{iN}(k)}  )$ is available  to player $i$, where
$\tau_{ij}(k) \leq k,j=1,\cdots,N$  are random  nonnegative integers    representing  communication delays from player
$j$ to player $i$ at  time $k$.  Set $\tau_{ii}(k)=0$ without loss of generality.

The following  Algorithm  \ref{asy-inexact-sbr-cont} presents the asynchronous inexact proximal  BR scheme,
where   $\mathcal{F}_k=\sigma \big\{ \mathcal{F}_k',  \{\tau_{ij}(k)\}_{i \in I_k,j \in \mathcal{N}} \big\}$ with $\mathcal{F}_k'$  being
the $\sigma$-field of the entire information employed by  the algorithm up to  (and including)  the update  of  $x_k$.
It is worth noticing that $y_{k}^i$ is adapted to $\mathcal{F}_k$ while  is not adapted to  $\mathcal{F}_k'$
since $y_{k}^i$ depends on $ \{\tau_{ij}(k)\}_{ j \in \mathcal{N}}$. We define $\mathcal{F}_k$ in Section \ref{sec:asy:complexity} while  the sequence   $\{\alpha_{i,k}\}$    will be specified
when  analyzing  the convergence properties    of  Algorithm \ref{asy-inexact-sbr-cont}.

\begin{algorithm} [H]
\caption{{Asynchronous inexact proximal  BR scheme}}
\label{asy-inexact-sbr-cont}
 Let $k:=0$, $x_{i,0}  \in X_i$  for $i \in \mathcal{N}$.
\begin{enumerate}
\item[(1)] For  any $i \in I_k$,  set $y_{k}^i=(x_{1,k-\tau_{i1}(k)},\cdots,x_{N,k-\tau_{iN}(k)}  )$.
\item[(2)] If  $i \in I_k$,   then     player $i$ updates $x_{i,k+1}$ that satisfies  the following:
 \begin{align} \label{inexact-sub2}
x_{i,k+1}  \in  \left\{  z\in X_i  : \mathbb{E}\left[\|z -
\wh{x}_{i}(y_k^i)\|^2\big |\mathcal{F}_k\right]   \leq \alpha_{i,k}^2~~a.s. \right\},
\end{align}  where $\wh{x}_{i}(y_k^i)$ is defined in \eqref{def-wh_x_inft} with $y_k$  replaced by $y_k^i$.
  Otherwise,   $x_{i, k+1}=x_{i, k}$.
\item [(3)] $k:=k+1$; If $k < K$, return to (1); else STOP.
\end{enumerate}
\end{algorithm}

  We impose the following assumptions on the asynchronous protocol,
  and on   the  parameters  $\zeta_{i,\min} , \zeta_{ij,\max} ,i,j=1,\cdots,N$ that are defined in \eqref{minmax-twice-diff}.

\begin{assumption}~\label{assp-asy}
(a) The sequence of sets  $\{I_k\}_{k\geq 0}$ is  deterministic.\\
(b)  Each player $i\in {\cal N }$ updates  its strategy   at least once  during any time interval of length $b_i$. In addition,   there exists a   positive integer $B_1$ such that $b_i\leq B_1~\forall i\in {\cal N }$;  \\
(c) There exists a random variable  $\tau_{ij}$ such that  $\tau_{ij}(k) \leq \tau_{ij}~\forall k \geq 0, i \in I_k,j \in \mathcal{N} ~~a.s.$ Furthermore,   there  exists a nonnegative   integer   $B_2$  such that    $\tau_{ij} \leq B_2~\forall i,j\in \mathcal{N}~a.s.$
\end{assumption}
By Assumption \ref{assp-asy}(b) we see that for any $k\geq 0$ and any $i \in I_k,j \in \mathcal{N}  $
 \begin{equation}\label{bd-delay}
\max\{0,k-B_2\}\leq  \max\{0,k-\tau_{ij}\}\leq  k-\tau_{ij}(k) \leq k \quad  a.s.
\end{equation}
\vspace{-0.25in}
 \begin{assumption}~\label{assp-dd}(Strict Diagonal Dominance)
For any $  i=1,\cdots,N$,
$\zeta_{i,\min} > \sum_{j \neq i} \zeta_{ij,\max} .$
\end{assumption}
 Then from this assumption   by $\Gamma$ defined in
  \eqref{matrix-hessian},  we  have that $a_{\infty}\triangleq  \|\Gamma\|_{\infty} <1$.
 Let $y: {\Omega} \to X$ be a
random variable defined on the probability space $({ \Omega}, {\cal F},\mathbb{P})$.
Then by taking expectations on  both sides of \eqref{cont-prox-best-resp},
and   by taking the infinity norm,  we have that
\begin{align}
\left\|  \pmat{\mathbb{E} [\|\wh{x}_1(y) -\wh{x}_1(x^*)\|] \\
				\vdots \\
\mathbb{E} [ \|\wh{x}_N(y) - \wh{x}_N(x^*)\|]}  \right\|_{\infty} \leq a_{\infty} \left\|  \pmat{ \mathbb{E} [ \|y_1- x_1^*\| ]\\
				\vdots \\
\mathbb{E} [ \|y_N - x^*_N\|]}  \right\|_{\infty} .\nonumber
\end{align}
As a result, by  $x_i^* = \wh{x}_i (x^*)$ we see that  for any $i=1,\cdots,N$:
 \begin{align}
\mathbb{E} [\|\wh{x}_i(y) -x_i^*\|]=\mathbb{E} [\|\wh{x}_i(y) -\wh{x}_i(x^*)\|] \leq a_{\infty} \max_ {i\in {\cal N}} \mathbb{E} [ \|y_i- x_i^*\| ]  .\label{cont-prox-resp-2}
\end{align}

 \subsection{Rate of convergence }  Denote by $\beta_{i,k}$  the number of updates player
 $i$ has carried out up to (and including) major iteration $k.$
\begin{lemma}\label{asy-lem-gemo}~
Let the asynchronous inexact proximal  BR scheme (Algorithm \ref{asy-inexact-sbr-cont}) be applied to
the $N$-player stochastic Nash game (SNash),  where  $ \alpha_{i,k}= \eta^{\beta_{i,k}}$ for some scalar $\eta \in (0,1)$ and    $\mathbb{E}[\|x_{i,0}-x_i^*\|] \leq C$.
  Suppose  Assumptions   \ref{assump-play-prob},  \ref{assp-asy},  and \ref{assp-dd} hold.
 Then we have the following for any $  k \geq 0$:
\begin{align} \label{asy-geom-err}  \max_{i \in {\cal N}}\mathbb{E}[\|x_{i,k} -x_i^*\|]  &  \leq (C+k)\rho^{\left \lfloor  {k \over B_1} \right\rfloor},
\end{align}
 where $\rho = \left(\max\{a_{\infty},\eta\}\right)^{1/(n_0+1)}$   with $ n_0 \triangleq  \left \lceil \frac{ B_2}{B_1} \right \rceil.   $   Furthermore,
if $q >c\triangleq  \rho^{  {1 \over B_1}}$, and   $D \red{\geq}1/\ln((q/c)^e)$,
\begin{align} \label{upper-asy-geom-err} \max_{i \in {\cal N}} ~\mathbb{E}[\|x_{i,k} -x_i^*\|]  & \leq
\rho^{- {B_1-1 \over B_1}} (C+D) q^k, \quad  \forall k \geq 0.
\end{align}
\end{lemma}
{\bf Proof.} See Appendix \ref{app-lem7}. \hfill $\Box$

 \subsection{Overall iteration  complexity analysis}\label{sec:asy:complexity}
We  proceed to derive a bound on the overall iteration complexity, when the inexact proximal BR
 solution  is computed via  SA. At   major iteration $k$,   if $i \in I_k$,    then  given $y_{k}^i,$ the $i$th  player takes the following sequence of {\em stochastic} gradient steps from $t = 1, \hdots, j_{i,k}$
\begin{align} \label{SA2}
z_{i,t+1} & := \Pi_{X_i} \left[ z_{i,t} -  		 \gamma_{t}  { \left(\nabla_{x_i}
 \psi_i(z_{i,t},y_{-i,k}^i;\xi_{i,k}^t) + \mu (z_{i,t}-x_{i,k})\right) } \right],
\end{align}
where  $z_{i,1}=x_{i,k}$, and $ \gamma_{t}=1/\mu(t+1)$. 	
Define   $\mathcal{F}_k'\triangleq\sigma\{x_0,\xi_{i,l},  \{\tau_{i}(l)\}_{i \in \textrm{I}_l},0 \leq l \leq k-1\}$ and  $\mathcal{F}_k\triangleq \sigma\{ \mathcal{F}_k',  \{ \tau_i(k)\}_{ i\in I_k} \}$,
where $\tau_i(k)\triangleq (\tau_{i1}(k),\cdots,\tau_{iN}(k))$.
  Then by Algorithm \ref{asy-inexact-sbr-cont},  we see that $x_k$ is adapted to   $\mathcal{F}_k'$ while
  $\{y_k^i\}_{ i \in I_k}$  is adapted to $\mathcal{F}_k.$
  {Since  $\mathcal{F}_k'\subset \mathcal{F}_k$,
  similar to Lemma \ref{lem-rate-sa}, we   also have   the following result  for  the scheme  \eqref{SA2}.}

\begin{lemma}~\label{lem-rate-sa2} Let Assumption~\ref{assump-play-prob} hold.
\red{Consider the asynchronous inexact proximal BR scheme given by   Algorithm \ref{asy-inexact-sbr-cont}.
Assume  that  for any $  i \in I_k$ and  any $t=1,\cdots, j_{i,k},$  $  \mathbb{E}\left[  \psi_i(z_{i,t},y^i_{-i,k};\xi_{i,k}^t)  \big | \sigma\{\mathcal{F}_k, \xi_{i,k}^{[t-1]}\}\right]  =\nabla_{x_i}f_i(z_{i,t},y^i_{-i,k})   ~a.s.,$ and   $  \mathbb{E}\left[ \|  \psi_i(z_{i,t},y^i_{-i,k};\xi_{i,k}^t) \|^2 \big | \sigma\{\mathcal{F}_k, \xi_{i,k}^{[t-1]}\}\right]  =
\| \nabla_{x_i}f_i(z_{i,t},y^i_{-i,k})  \|^2 ~a.s.$     Then for any $ i \in I_k$ and    any $t=1,\cdots, j_{i,k} $, we obtain  that
 $ \mathbb{E}[\|z_{i,t}-\wh{x}_i(y_k^i)\|^2\big | \mathcal{F}_k] \leq   {Q}_i/{(t+1)}~~a.s.  $}
\end{lemma}

 Distinct from  the definition of $\epsilon-$Nash equilibrium given in \eqref{def-epsilon-NE1}, for the asynchronous  algorithm  a
random strategy profile $x: { \Omega} \to \Real^{n}$ is called an  {$\epsilon-$NE$_{\infty}$} when
\begin{align} \label{def-epsilon-NE2} \max_{i\in \mathcal{N}}~\mathbb{E}  [\left\| x_{i} -{x}_{i}^* \right \|]  \leq \epsilon.
\end{align}
 \begin{theorem}~ \label{thm-complexity} Let Algorithm \ref{asy-inexact-sbr-cont}    be applied to the
  stochastic Nash game (SNash), where   $ \alpha_{i,k}= \eta^{\beta_{i,k}}$ for some scalar $\eta \in (0,1)$
   and $ \mathbb{E}[\|x_{i,0}-x_i^*\|^2] \leq C^2$.    Suppose that  an inexact solution   characterized  by \eqref{inexact-sub2}
is computed via   \eqref{sa-inner}. Let    Assumptions \ref{assump-play-prob},  \ref{assp-asy}  and \ref{assp-dd} hold.
Suppose    $\rho=\left(\max\{a_{\infty},\eta\}\right)^{1/(n_0+1)}$ with  $n_0=\left \lceil {B_2 \over B_1}\right\rceil $.  Define $c\triangleq  \rho^{  {1 \over B_1}}$ and let  $q \in (c,1)$.  Then the number of projected gradient steps  for player $i$ to compute an  {$\epsilon-$NE$_{\infty}$}  is no greater than
\begin{align}\label{asy-comp-bound}
	\ell_i^{(1)}( \eta)\triangleq \frac{Q_i}{ \eta^4\ln(1/\eta^2)}
		\left(\frac{ 1}{\hat{\epsilon} }\right)^{ \frac{\ln(1/\eta^2)}{\ln(1/q )}}+\left \lceil\frac{\ln(1/\hat{\epsilon})}{\ln(1/q)} \right\rceil,
\end{align}
where $\hat{\epsilon} $ is defined in  \eqref{asy-def-epsilonbar} \red{with $D \geq /\ln((q/c)^e)$.}
 In particular,  if $B_1=1,B_2=0$, then  \eqref{asy-comp-bound} reduces to
$\frac{Q_i}{ \eta^4\ln(1/\eta^2)}
			\left(\frac{C+D}{\epsilon }\right)^{ \frac{\ln(1/\eta^2)}{\ln(1/q )}}$ with $q >c\triangleq  \max\{a_{\infty},\eta\}  $  and   $D \red{\geq}1/\ln((q/c)^e)$.
\end{theorem}
\noindent {\bf Proof.}     For any $i \in  I_k$,   let  $j_{i,k}=\lceil \frac{Q_i}{  \eta^{2(k+1)}}\rceil$ be the number of steps of   the scheme \eqref{SA2}  taken by player $i$ in major iteration $k$ to obtain $x_{i,k+1}$, i.e.,   $x_{i,k+1}:=z_{i,j_{i,k}}$.  Since $\beta_{i,k}\leq k+1$ and $\eta \in(0,1)$,
by Lemma \ref{lem-rate-sa2},  we have the following for any  $i \in  I_k$:
$$  \mathbb{E}\left[\|x_{i,k+1}-\wh{x}_i(y_k^i)\|^2\big| \mathcal{F}_k\right]  = \mathbb{E}\left[\|z_{i,j_{i,k}}-\wh{x}_i(y_k^i)\|^2\big| \mathcal{F}_k\right]  \leq	\frac{ Q_i }{j_{i,k}} \leq  \eta^{2(k+1)} \leq  \eta^{2\beta_{i,k}}  = \alpha_{i,k}^2 ~~a.s.$$ By \eqref{upper-asy-geom-err},  we may obtain a lower bound on $k$ as follows:
\begin{align}
  &\max_{i \in {\cal N}}\mathbb{E}[\|x_{i,k+1}-x_i^*\|]     \leq \rho^{- {B_1-1 \over B_1}} (C+D) q^{k+1}  \leq \epsilon   \Rightarrow
q^{k+1} \leq  {\epsilon  \over C+D }  \rho^{  {B_1-1 \over B_1}}\triangleq \hat{\epsilon}  \label{asy-def-epsilonbar}
  \Rightarrow  k \geq \frac{\ln(1/\hat{\epsilon})}{\ln(1/q)} -1.  \end{align}
{Then by  invoking \eqref{exponent-sum}
 for  $\beta > 1$,  we derive the following    complexity bound for player $i$:
\begin{align*}
\sum_{k=0}^{\left\lceil\frac{\ln(1/\hat{\epsilon})}{\ln(1/q)}\right\rceil-1}  j_{i,k}  &\leq \sum_{k=0}^{\left\lceil\frac{\ln(1/\hat{\epsilon})}{\ln(1/q)}\right\rceil-1} \left( \frac{Q_i}{  \eta^{2(k+1)}}	+1\right)
= \sum_{k=1}^{\left\lceil\frac{\ln(1/\hat{\epsilon})}{\ln(1/q)}\right\rceil}  \frac{Q_i}{  \eta^{2k}}	 +\left \lceil\frac{\ln(1/\hat{\epsilon})}{\ln(1/q)} \right\rceil
\\ & \leq \frac{Q_i}{ \ln(1/\eta^2)}   \eta^{-2\left(\frac{\ln(1/\hat{\epsilon})}{\ln(1/q)}+2 \right)}
+\left \lceil\frac{\ln(1/\hat{\epsilon})}{\ln(1/q)} \right\rceil,\end{align*}
which results in  \eqref{asy-comp-bound}  by invoking \eqref{change-exp}.} \hfill $\Box$

 Coordinate  descent   methods  advocate that  only a small  block  of variables are updated in each iteration while others are  kept  fixed \cite{xu2015block}.  The blocks of variables can be updated in  a {\it cyclic, random,} or {\it greedy} fashion.  In Section \ref{sec:random},   we examined  randomized    update schemes, whereas here we consider a cyclic
  update scheme  to  broader problem settings, and analyze  the iteration complexity.

\begin{corollary}   \label{prp-cyclic}~ Let Algorithm \ref{asy-inexact-sbr-cont}    be applied to the
  stochastic Nash game (SNash),  where    $ \alpha_{i,k}= \eta^{\beta_{i,k}}$ for some $\eta \in(0,1)$.
  Suppose that  an inexact solution   characterized  by \eqref{inexact-sub2}
is computed via  \eqref{sa-inner}.  Let   Assumptions \ref{assump-play-prob},  \ref{assp-asy}  and \ref{assp-dd} hold.  Suppose that the players are chosen to update  cyclicly,   and  {  $\rho =\left(\max\{a_{\infty},\eta\}\right)^{1/(n_0+1)}$  with   $n_0=\left \lceil {B_2 \over N}\right\rceil$.
Define $c\triangleq  \rho^{  {1 \over  N}} $ and let  $q \in (c,1)$. } If $\tilde{  \epsilon}\triangleq \frac{\epsilon }{C+D}\rho^{N-1 \over N}$
\red{with $D =1/\ln((q/c)^e)$},   $\tilde{\eta}=\eta^{1/N},$ then the number of projected gradient steps for player $i$ to compute an  {$\epsilon-$NE$_{\infty}$}  is no greater than
$$ 	 \ell_i^{(2)}( \eta)\triangleq \frac{Q_i}{ \tilde{\eta}^2 \eta^2\ln(1/\tilde{\eta}^2)}
			\left( { 1\over \tilde{  \epsilon}}\right)^{ \frac{\ln(1/\tilde{\eta}^2)}{\ln(1/q )}}
			+\left \lceil\frac{\ln\left( 1/\tilde{  \epsilon}  \right)}{\ln(1/q)} \right\rceil.$$		
\end{corollary}

\noindent {\bf Proof. }   For any $i \in  I_k$,    let  $j_{i,k} = \left\lceil \frac{Q_i}{ \eta^{2+2k /N}}\right\rceil$
 be the steps of the SA scheme   \eqref{SA2}  computed for player $i$ at major iteration  $k$. Set  $x_{i,k+1}:=z_{i,j_{i,k}}$.  By the cyclic rule,  we see that $B_1=N$, and each  player $i\in \cal{N}$   is picked
 exactly once among $N$ steps.  As a result,  $\beta_{i,k}=\left \lceil\frac{k+1}{N}\right \rceil$ when   $i \in  I_k$,  and hence
  $\beta_{i,k} \leq {k \over N}+1.$  Then by $\eta \in(0,1)$ and Lemma \ref{lem-rate-sa},   we obtain that
 \begin{align*} \mathbb{E}\left[\|x_{i,k+1}-\wh{x}_i(y_k^i)\|^2\big| \mathcal{F}_k\right]
  &= \mathbb{E}\left[\|z_{i,j_{i,k}}-\wh{x}_i(y_k^i)\|^2\big| \mathcal{F}_k\right] \leq
	\frac{ Q_i }{j_{i,k}} \leq  \eta^{2+2k /N} \leq  \eta^{2\beta_{i,k}}  = \alpha_{i,k}^2,\quad\forall i \in  I_k.\end{align*}
{Then by similar procedure for deriving \eqref{asy-comp-bound},   we may obtain the following  iteration complexity  for player $i$:
 \begin{align*}
 \sum_{k=0}^{\left\lceil\frac{\ln(1/\hat{\epsilon})}{\ln(1/q)}\right\rceil -1}   j_{i,k}
 & \leq   \frac{1}{ \eta^{2(N-1) /N}}  \sum_{k=0}^{\left\lceil\frac{\ln(1/\hat{\epsilon})}{\ln(1/q)}\right\rceil -1}
  \frac{Q_i}{ \tilde{\eta}^{2(k+1) }}+\left \lceil\frac{\ln(1/\hat{\epsilon})}{\ln(1/q)} \right\rceil
\\& \leq \frac{1}{ \eta^{2(N-1) /N}}   \frac{Q_i}{ \tilde{\eta}^{4}\ln(1/\tilde{\eta}^2)}
			\left(\frac{1}{ \hat{\epsilon}}\right)^{\frac{\ln(1/\tilde{\eta}^2)}{\ln(1/q )}} +\left \lceil\frac{\ln(1/\hat{\epsilon})}{\ln(1/q)} \right\rceil.
\end{align*}
Thus,  the result follows by $\tilde{\eta}=\eta^{1/N}$
  and by   $\hat{\epsilon} $ defined in \eqref{asy-def-epsilonbar} with $B_1=N$.}
\hfill $\Box$

\subsection{Complexity bound }
In the following  theorem, we  drive  the upper bound for    $\ell_i^{(1)}( \eta)$ and  $\ell_i^{(2)}( \eta)$  defined in Theorem \ref{thm-complexity}
  and Corollary \ref{prp-cyclic}, respectively.

\begin{theorem} \label{asy-upper-bound} ~
 Set $\eta \triangleq a_{\infty}$.  Then  we have
  the following  bounds on the expression $ \ell_i^{(1)}(\eta)$ and $ \ell_i^{(2)}(\eta)$.

a) {Define $n_0= \left\lceil {B_2 \over B_1}\right\rceil$ and   $n'=B_1(1+n_0)$. Then the following holds for any  given $\delta > 0$:}
$$ \ell_i^{(1)}(\eta)   \leq \frac{Q_i}{\eta^4\ln(1/\eta^2)}\left(\frac{1}{\hat{\epsilon} }\right)^{2n'+\delta}
+\left\lceil \ln\left({ 1 \over \hat{\epsilon }} \right) { 2n'+\delta \over 2\ln( \eta^{-1}) } \right\rceil,$$
where { $ \hat{\epsilon}={  \epsilon \over C+D } \eta^{  -{B_1-1 \over B_1(n_0+1)}} $
with  $ D $ defined in \eqref{def-D}.} In particular,  if    $B_1=1$,
  then $$ \ell_i^{(1)}(\eta)   \leq \frac{Q_i}{\eta^4\ln(1/\eta^2)} \left(\frac{C+D}{\epsilon}\right)^{2(1+B_2)+\delta} +\left\lceil \ln\left({ C+D \over \epsilon}  \right) { 2(B_2+1)+\delta \over 2\ln( \eta^{-1}) } \right\rceil,$$
 where $ D= {(1+B_2)+2(1+B_2)^2\delta^{-1}\over e \ln(\eta^{-1})}$.
If $ D= {1+2\delta^{-1}\over e \ln(\eta^{-1})}$ and $B_2=0$, then given any $\delta>0$, we have that
 \begin{align}\label{asy-syn}
   \ell_i^{(1)}(\eta)   \leq \frac{Q_i}{\eta^4\ln(1/\eta^2)}
\left(\frac{C+D}{\epsilon}\right)^{2 +\delta}+\left\lceil \ln\left({ C+D \over \epsilon}  \right) { 2+\delta \over 2\ln( \eta^{-1}) } \right\rceil,
\end{align}
	 		
b)  Given any  $\delta > 0$, $\tilde{  \epsilon}\triangleq \frac{\epsilon }{C+D}\rho^{N-1 \over
N(1+n_0)}$, $n_0=\left \lceil \frac{ B_2}{N}\right \rceil$, $D $ given
by \eqref{asy-parameter-D}, we have that
$$\ell_i^{(2)}(\eta)  \leq \frac{Q_i}{\eta^{2+{2\over N}}\ln(1/\eta^{2 \over N})}
\left( \frac{1}{\tilde{\epsilon} } \right)^{2(n_0+1)+\delta} + \left \lceil \ln\left( {1 \over \tilde{  \epsilon} } \right)  {  2(n_0+1)+ \delta  \over 2\ln( \eta^{-1})/N}\right\rceil .$$
 In particular,  if     $B_2=0$, then given any  $\delta > 0$,  $\tilde{  \epsilon}\triangleq \frac{\epsilon }{C+D}\rho^{N-1 \over N }$ and $D={ N\left(1+2  \delta^{-1}\right) \over e \ln(\eta^{-1})}$, we have that
	$$\ell_i^{(2)}(\eta)  \leq \frac{Q_i}{\eta^{2+{2\over N}}\ln(1/\eta^{2 \over N})}	\left( \frac{1}{\tilde{  \epsilon} } \right)^{2 +\delta}+ \left \lceil \ln\left( {1 \over \tilde{  \epsilon} } \right)  {  (2 + \delta)N  \over 2\ln( \eta^{-1}) }\right\rceil .$$
\end{theorem}
\noindent {\bf Proof.} a)    Define $ q= c e^{\delta_0/2n'}$ with $\delta_0=\frac{\delta \ln(\eta^{-1})}{  n'+\delta/2}$.
{By $\eta=a_{\infty}$,  we  get  $\rho^{1+n_0}=\eta$ and $ c=  \rho^{  {1 \over B_1}}=\eta^{1 \over n'}$.    Then
 \begin{align}
 \ln(1/q )&=  \ln \left( \eta^{- 1/ n' }e^{-\delta_0/2n'} \right)=
 {\ln( \eta^{-1})\over n'}-{ \delta_0 \over 2n' } = {2\ln( \eta^{-1})\over
2n'+\delta} , \frac{\ln(1/\eta^2)}{\ln(1/q )} =   \frac{2
\ln(\eta^{-1})}{\ln(1/q )} =2n'+ \delta \notag
   \\{\rm and ~~}D&= 1/\ln((q/c)^e)=1/\ln(e^{\delta_0 e/2n'})={2n' \over \delta_0 e}={n'+2(n')^2\delta^{-1}\over e \ln(\eta^{-1})} \label{def-D}. \end{align}}
Consequently, by the  expression $\ell_i^{(1)}(\eta)$ defined in \eqref{asy-comp-bound},
 we  show the first assertion of  part (a). The  remaining  two assertions can be easily  followed  by
 the definitions of $n_0,n', $ and $D.$

b)   Define $ q= c e^{\delta_0/2N(1+n_0)}$ with $\delta_0=\frac{\delta \ln(\eta^{-1})}{  1+n_0+\delta/2}$.
 By $\eta=a_{\infty}$,  we  get $\rho^{n_0+1}=\eta$ and   $ c= \rho^{  {1 \over N}}=\eta^{1 \over N(1+n_0)}$.
  Then by $\tilde{\eta}=\eta^{1/N} $ we have the following:
 \begin{align}
 \ln(1/q )&=  \ln \left(  \eta^{-1/ N(1+n_0)} e^{-\delta_0/2N(1+n_0)} \right)=
  {\ln( \eta^{-1})- \delta_0/2\over N(1+n_0)} ={2\ln( \eta^{-1})/N \over  2(n_0+1)+ \delta  }, \notag  \\
 \frac{\ln(1/\tilde{\eta}^2)}{\ln(1/q )}&= \frac{\ln(1/\eta^{2/N})}{\ln(1/q )}=
   \frac{2 \ln(\eta^{-1})/N}{\ln(1/q )}   =2(n_0+1)+ \delta , \notag \\
{\rm and~~}    D&= 1/\ln((q/c)^e)=1/\ln(e^{\delta_0 e/2N(1+n_0)})={2N(1+n_0) \over \delta_0 e}=N(1+n_0){1+2(1+n_0) \delta^{-1}\over e \ln(\eta^{-1})}.
\label{asy-parameter-D}
   \end{align}
 Then by     $\ell_i^{(2)}( \eta)$  defined in  Corollary \ref{prp-cyclic},    we   obtain  the first assertion of  part  (b).
 If $B_2=0$, then by definitions of $n_0$ and $D$, we get  $n_0=0$ and $D={ N\left(1+2  \delta^{-1}\right) \over e \ln(\eta^{-1})}$. Thus,
 the remaining  assertion  of  part (b) holds.
\hfill $\Box$

The complexity bounds given in the part (a) and part (b) of  Theorem \ref{asy-upper-bound}  are summarized as follows:
\vspace{-0.2in}
\begin{table}[H]
 \scriptsize
    \begin{minipage}[c]{0.52\textwidth}
\centering
\subtable{
       \begin{tabular}{|c|c|c|}
        \hline
 \multicolumn{3}{|c|}{(a)-almost cyclic rule}   \\
 \hline
update frequency  &  delay & complexity bound \\ \hline
$B_1$ & $B_2$ &  $\mathcal{O}\left( (1/\epsilon)^{2B_1 \left(1+\left \lceil {B_2 \over B_1}\right \rceil \right)+\delta}\right)$     \\ \hline
1 & $B_2$ &  $\mathcal{O}\left( (1/\epsilon)^{2(1+B_2)+\delta}\right)$  \\ \hline
1 & $0$ &  $\mathcal{O}\left ((1/\epsilon)^{2+\delta}\right)$   \\ \hline
 \end{tabular}
}
\end{minipage}
 \begin{minipage}[c]{0.4 \textwidth}
\centering
\subtable{   \begin{tabular}{|c|c|c|}
 \hline
 \multicolumn{3}{|c|}{(b)-cyclic rule}  \\
 \hline
update frequency  &  delay & complexity bound \\ \hline
   $N$ & $B_2$ & $\mathcal{O} \left( (1/\epsilon)^{2\left ( \left \lceil {B_2 \over N} \right \rceil+1\right)+\delta}\right)$ \\ \hline
  $N$ & $0$ & $\mathcal{O}\left ((1/\epsilon)^{2+\delta}\right)$ \\ \hline
 \end{tabular}
}
\end{minipage}
      \caption{Summary of  Complexity Bounds  \label{TAB1}}
\vspace{-0.2in}
\end{table}

\begin{remark} ~Based on established results,  we  make the   following clarifications:  \\
	 \noindent {\bf (i) } If $B_1=1,~B_2=0$, then  the asynchronous  algorithm degenerates to the   synchronous  algorithm.
 Notably, we achieve  the complexity bound   ${\cal O} \left( 1/\epsilon^{2 +\delta} \right)$ given by  \eqref{asy-syn}  in Theorem \ref{asy-upper-bound}(a) compared  with    ${\cal O} \left( \sqrt{N}/\epsilon \right)^{2 +\delta}$  given
	  in  Theorem \ref{prop-iter-comp}(b). This is because the definition  of  $\epsilon$-NE$_{2}$
 by \eqref{def-epsilon-NE1} (resp.    $\epsilon$-NE$_{\infty}$   by
\eqref{def-epsilon-NE2}) as well as the analysis of
Theorem \ref{prop-iter-comp}(b) and  Theorem \ref{asy-upper-bound} (a) are  based on  two-norm and    infinity norm, respectively; \\
	 \noindent {\bf (ii)} Assume that   the players update in a cyclic  manner and $B_2=0$.
 Then   the upper bound of iteration complexity is  shown to be ${\cal O} \left( 1/\epsilon^{2 +\delta} \right)$ in Theorem \ref{asy-upper-bound}(b),
 which is of the same order  as that of the synchronous case given by  \eqref{asy-syn}  in Theorem \ref{asy-upper-bound}(a).     \\
 \noindent {\bf (iii) }From Corollary \ref{prp-cyclic}  and Theorem \ref{asy-upper-bound}(b)   we infer that  for some specific
selection of update indices  in each iteration,  we can get a  much  better upper bound  than
that given in  Theorem \ref{asy-upper-bound}(a).  So, the  update schemes of players are  critical in improving  the   complexity bound.  Nevertheless, these  findings imply that the complexity  bounds established   in Theorem \ref{asy-upper-bound}  for Algorithm \ref{asy-inexact-sbr-cont}    are optimal to some extent.\\
 \noindent {\bf (iv) }  From  the analysis of Algorithm \ref{asy-inexact-sbr-cont},   the geometric  rate of convergence and  the iteration complexity   can also be established for Algorithm \ref{inexact-sbr-cont}   and  Algorithm \ref{rand-inexact-sbr-cont}  by replacing the assumption $\| \Gamma\|<1$ with  $\| \Gamma\|_{\infty}<1.$
\end{remark}

\begin{remark} ~ The assumption $\| \Gamma\|_{\infty}<1$
 can be further  weakened to $\rho(\Gamma)<1$  based on  the following:
By  \cite[Corollary 2.6.1]{bertsekas1989parallel}, there exists a nonnegative vector $w \in \mathbb{R}^{N}$ such that
$a_{\infty} \triangleq \max_{i \in \cal{N}} \frac{1}{w_i} \sum_{j=1}^N \gamma_{ij} w_{j}<1.$ Then by \eqref{cont-prox-best-resp}, we have that
\begin{align*}
 diag\{ w\}^{-1}  \pmat{\|\wh{x}_1(y') -\wh{x}_1(y)\| \\
				\vdots \\
\|\wh{x}_N(y') - \wh{x}_N(y)\|}  \leq    diag\{ w\}^{-1}   \Gamma   diag\{ w\}  diag\{ w\}^{-1}    \pmat{\|y_1'- y_1\| \\
				\vdots \\
\|y'_N - y_N\|}  ,
\end{align*}
where  $diag\{ w\} $ is a diagonal matrix with the diagonal entries being $w_i~i\in \mathcal{N}$.
Similar to that given by the proof \cite[Proposition 5]{scutari2014real} we have  the following inequality:
\begin{align*}
  \left \|  diag\{ w\}^{-1}  \pmat{\|\wh{x}_1(y') -\wh{x}_1(y)\| \\
				\vdots \\
\|\wh{x}_N(y') - \wh{x}_N(y)  \|}\right  \|_{\infty}  \leq a_{\infty}  \left \|  diag\{ w\}^{-1}   \pmat{\|y_1'- y_1\| \\
				\vdots \\
\|y'_N - y_N\|} \right \|_{\infty} ,
\end{align*}
by which we  are able to establish the geometric   rate of convergence and  the same order   complexity bound for Algorithm  \ref{asy-inexact-sbr-cont} but with different constants.
\end{remark}

\section{Two-stage Recourse  }\label{sec:recourse}

In the previous   sections,  we assumed that the  functions
$f_i(x)$ are   twice   differentiable, a requirement that is
essential in deriving the contractive properties of the  proximal  BR map.
{We   weaken this assumption by  considering a stochastic Nash game in   the following form:}
\begin{align*}
\min_{x_i \in X_i} \ f_i(x_i,x_{-i}) + \red{g_i(x_i)},  
\end{align*}
where $f_i(x_i,x_{-i}) \triangleq \mathbb{E} \left[\psi_i(x_i,x_{-i};\xi({\omega})) \right]$ and $ \red{g_i(x_i)}$ can be merely
convex, possibly nonsmooth, and  expectation-valued. It is relatively simple to show that if the proximal BR map,
	constructed using the functions $f_1, \hdots, f_N$, admits a
		contractive property, then this modified problem  also  admits such a  property.  Consequently, the only challenge that emerges is the need to get inexact
solutions to the proximal  BR problem with an objective $f_i(x) +  \red{g_i(x_i)}$ instead of merely $f_i(x)$.
In this section, we  consider a  two-stage recourse-based model to  introduce
much needed flexibility into this framework.

  Two-stage recourse-based stochastic programs  originate from the work by
Dantzig~\cite{ dantzig1955linear} and Beale~\cite{beale55minimizing} in the 50s.
  In the first stage, the
decision-maker makes a decision prior to the revelation of the
uncertainty,  while in the second-stage, a scenario-specific   decision  is made contingent on the first-stage
decision and on the realization of uncertainty. This second-stage decision is
referred to as a {\em recourse} decision, which originates  from the
notion that in the second-stage, the decision-maker takes {\em
	recourse} based on the realization of uncertainty and on the
	first-stage decision. The expected cost of second-stage decisions is
	incorporated into the first-stage problem through the {\em recourse}
	function. This model has tremendous utility and finds applicability
	in financial planning, inventory control, power systems operation, etc
	  (cf.~\cite{birge97introduction,shapiro09lectures}).
	 Nevertheless,  such models and algorithms are not designed to
accommodate multiple decision makers who do not cooperate in a multi-agent system.
Recently, in \cite{pang2017two},  the authors consider two-stage non-cooperative games, where each agent is risk-averse and solves a rival-parameterized stochastic program with quadratic
recourse.  To handle  the nonsmooth recourse functions, the authors    propose smoothing
  schemes leading to differentiable approximations, then  design an iterative BR scheme
	 for  the smoothed problem and  show its convergence.  We propose a different avenue that leverages stochastic approximation and allows for deriving rate and complexity statements, that are unavailable in~\cite{pang2017two}.  Next, we  incorporate two-stage  linear
and  quadratic recourse into this framework  and derive iteration complexity statements for the proposed BR schemes.
 \subsection{Linear   recourse}\label{sec-linear-recourse}
 Consider the following  two-stage stochastic Nash game:
\begin{align}  \label{SNash_rec} \tag{SNash$_{\rm rec}(x_{-i})$}
\min_{x_i \in X_i} ~ \mathbb{E}\left[\psi_i(x_i,x_{-i}; \red{\xi(\omega)}\right] +
c_i(x_i) + \mathbb{E}[{\cal Q}_i(x_i, \omega  ) ],
\end{align}
where $c_i(x_i)$ is the continuously differentiable  convex cost of the first-stage decision $x_i$ and
\begin{align}
\tag{Rec$_{\rm LP}(x_i)$}
{\cal Q}_i(x_i,\omega) \triangleq
\min_{\red{q_{i,\omega}} } \quad    \left\{d_{i,\omega}^T \red{q_{i,\omega}}  \mid    T_{i,\omega} x_i +W_i \red{q_{i,\omega}}  = h_{i,\omega} ,\, \red{q_{i,\omega}}  \geq 0\right\} \label{two-stage}.
\end{align}
The function ${\cal Q}_i(x_{i},\omega)$,  being the cost of taking recourse given decision $x_i$ and scenario $\omega$,
 is shown to be convex in   $x_i \in X_i$ for any given $\omega \in \Omega$ (See~\cite[Proposition 2.1]{shapiro09lectures}).

Given a  strategy $y \in X$, player $i$ can {compute} a proximal   BR $\wh{x} _i(y)$ by solving the following  stochastic optimization problem
\begin{align}\label{eq-prox-recourse}
\wh{x} _i(y) \triangleq 
{\operatornamewithlimits{\mbox{argmin}}_{x_i \in X_i}}
 \left[   \mathbb{E}
[\psi_i(x_i,y_{-i};\omega)] +c_i(x_i) + \mathbb{E}[{\cal Q}_i(x_i,\omega)]+ {\mu \over 2} \|x_i-y_i\|^2 \right].
\end{align}
It is easily seen  that $\wh{x}(\cdot) $   defined by \eqref{eq-prox-recourse}  is a contractive map since $ \red{g_i(x_i)}:=c_i(x_i) + \mathbb{E}[{\cal Q}_i(x_i,\omega)]$ is convex.  We can also use  Algorithms \ref{inexact-sbr-cont},  \ref{rand-inexact-sbr-cont},  and
  \ref{asy-inexact-sbr-cont} to  find the Nash equilibrium of  the problem (SNash$_{\rm rec}(x_{-i})$)  but with  $\wh{x} _i(\bullet) $ defined  by \eqref{eq-prox-recourse} instead of  \eqref{eq-br}. We intend to use the SG method to obtain  an inexact   proximal   BR  solution.  So, we impose the following condition on the problem (Rec$_{\rm LP}(x_i)$)   to guarantee  the existence and boundedness  of the stochastic
 subgradient for  function $ \mathbb{E}\left[{\cal Q}_i(x_i,\omega)\right]$.
\begin{assumption}~\label{assp-exist}
 For almost all $\omega \in \Omega$ and  all $x_i \in X_i$, \\
 a) there exists a strictly positive vector $\red{q_{i,\omega}}$ such that
 $ T_{i,\omega} x_i +W_i \red{q_{i,\omega}}  = h_{i,\omega} $,     \\
b)  for any  $v_{i,\omega}  \geq 0$ with $W_i v_{i,\omega}  = 0$, it holds that $d_{i,\omega}^T v_{i,\omega} \geq 0 $,\\
c) $\mathbb{E}[\|T_{i,\omega}\|]<\infty$ for all $i \in \cal{N}$.
  \end{assumption}
\uvs{We adapt the following result from~\cite{shapiro09lectures} to
	analyze the subgradient of ${\cal Q}(x_i,\omega)$.}
\begin{lemma} ~Let Assumption \ref{assp-exist}  hold. Then for almost all $\omega \in \Omega$ and  all $x_i \in X_i$, \\
(i) ${\cal Q}_i(x_i,\omega)$ is finite  and its  subdifferential at $x_i$ is given by
\begin{equation}\label{subg-recourse}
\partial_{x_i} {\cal Q}_i(x_i,\omega)  = -T_{i,\omega}^T{\cal D}_i(x_i,\omega),
 \mbox{ where~ } {\cal D}_i(x_i,\omega)  \triangleq \argmax_{
		\pi_{i,\omega}:
		W_i^T\pi_{i,\omega} \leq d_{i,\omega} }  \
			(h_{i,\omega}-T_{i,\omega} x_i)^T\pi_{i,\omega};
\end{equation}
(ii) ${\cal D}_i(x_i,\omega) $ defined by  \eqref{subg-recourse} is  bounded;\\
(iii)  there exists a positive constant $M_s>0$ such that
 \begin{equation}
\label{subg-bound}
E[\| s_i \|^2]\leq M_s^2,\quad\forall s_i \in \partial {\cal Q}_i(x_i,\omega) .
\end{equation}
\end{lemma}
{\bf Proof.}  (i)    Assumption \ref{assp-exist}(a) implies   the existence of  a feasible solution  to problem  (Rec$_{\rm LP}(x_i)$).
  Then by Assumption \ref{assp-exist}(b), we see that   for almost all $\omega \in \Omega$ and  all $x_i \in X_i$,
  (Rec$_{\rm LP}(x_i)$) has an optimal solution and
  ${\cal Q}_i(x_i,\omega)$ is finite.  Thus,  by  \cite[Proposition 2.2]{shapiro09lectures}  we obtain the result.

(ii)  Note that $ h_{i,\omega} -T_{i,\omega} x_i $ is an interior point of the positive hull of matrix $W_i$ by   Assumption \ref{assp-exist}(a).  It is shown in \cite[p. 29]{shapiro09lectures} that the optimal solution set ${\cal D}_i(x_i,\omega) $ of \eqref{subg-recourse} is  bounded.

(iii)  By   Assumption \ref{assp-exist}(c) and the boundedness of ${\cal D}_i(x_i,\omega) $   we have the result.
\hfill $\Box$

Since  ${\cal Q}_i(x_i,\omega)$ is   subdifferentiable,  the following SA scheme is utilized  to approximate
 the proximal   BR $\wh{x} _i(y_k) $:
\begin{align} \tag{SA$^{rec}_{i,k}$}
z_{i,t+1} & := \Pi_{X_i} \left[ z_{i,t} 
		 -  \gamma_{t} \left(\nabla_{i}
c_i(z_{i,t}) + \nabla_{x_i} \psi_i(z_{i,t},y_{-i,k};\red{\xi_{i,k}^t}) +\mu (z_{i,t}-y_{i,k})
+s_{i,t}\right)  \right],  \label{sa-recourse}
\end{align}
where $z_{i,1}=x_{i,k}$,  $ \gamma_{t} = 1/\mu (t+1)$ and $s_{i,t} \in \partial  {\cal Q}_i(z_{i,t},\omega_{i,k}^t)$. In fact,   the SA scheme (SA$^{rec}_{i,k}$) requires a solution of the  second-stage dual problem  \eqref{subg-recourse}  to obtain a stochastic  subgradient.  Since $\nabla_{i}c_i(\cdot)$ is continuous and $X_i$ is compact, there exists $M_c>0$ such that
	$\|\nabla_{i}c_i(x_i)\|\leq M_c~~\forall x_i \in X_i.$
Combined  with Assumption  \ref{assump-play-prob}(d),    \eqref{diameter} and  \eqref{subg-bound} yields that
for any $x_i \in X_i,y \in X$
$$ \mathbb{E}\left[\|\nabla_{i}
c_i(x_i) + \nabla_{x_i} \psi_i(x_i,x_{-i};\xi) +\mu (x_i-y_i)
+s_i||^2  \right] \leq 4( M_i^2+M_s^2+M_c^2+\mu^2 D_{X_i}^2). $$
By \cite[Theorem 7.47]{shapiro09lectures},  under suitable regularity conditions, it is seen  that $\mathbb{E}\left[\partial  {\cal Q}_i(x_i,\omega)\right] = \partial \mathbb{E} \left[ {\cal Q}_i(x_i,\omega)\right]$,   and hence  $s_{i,t} \in \partial  {\cal Q}_i(z_{i,t},\omega_{i,k}^t)$
 is an unbiased estimate for some  subgradient of   $ \mathbb{E}\left[{\cal Q}_i(x_i,\omega)\right]$ at  point $z_{i,t}.$

Consequently,    Algorithms \ref{inexact-sbr-cont},    \ref{rand-inexact-sbr-cont}, and   \ref{asy-inexact-sbr-cont}  are applicable to the  stochastic Nash game with linear-recourse.  Further,  the convergence results   all   hold for the problem \eqref{SNash_rec} with  ${\cal Q}_i(x_i,\omega)$ defined  by (Rec$_{\rm LP}(x_i)$), except that  the constant   $Q_i \triangleq  \frac{2 M_i^2 }{ \mu^2}  +2 D_{X_i}^2 $ is replaced by $Q_i \triangleq   \frac{4( M_s^2+M_c^2+M_i^2 )}{\mu^2} +4  D_{X_i}^2 $.

 \subsection{Quadratic  recourse}
In the past, two-stage  stochastic quadratic  programs with fixed recourse have
been considered  in   \cite{rockafellar1986lagrangian,chen1995newton}, where
the first-stage term is a quadratic  function and the second-stage term  is the
expectation of   the minimum value of a quadratic program.  In
\cite{rockafellar1986lagrangian} and \cite{chen1995newton}, the authors  work
with the recourse subproblems'  dual that is approximated by a sequence of
quadratic program  subproblems, and propose a  Lagrangian finite generation
technique and a Newton's method to solve the problem, respectively.
In~\cite{kulkarni12recourse,shanbhag06decomposition}, Dorn duality is employed
in developing a Benders scheme.  Motivated by this, we incorporate   quadratic
recourse in the   stochastic Nash game and   consider the two-stage  problem
(SNash$_{\rm rec}(x_{-i})$), for which    $c_i(x_i)$  is  a continuously
differentiable  convex cost of   $x_i$, and   the recourse function ${\cal
Q}_i(x_i,\omega)$ is   the optimal value of a convex quadratic program
parameterized by  the decision $x_i$ and scenario $\omega$:
\begin{align}
\tag{Rec$_{\rm QP}(x_i)$}
\min_{\red{q_{i,\omega}}} \quad    &  \left\{ d_{i,\omega}^T \red{q_{i,\omega}} +{1 \over 2} \red{q_{i,\omega}} ^T H_{i,\omega}  \red{q_{i,\omega}} \mid  \red{q_{i,\omega}} \in   \Upsilon_{i,\omega}\triangleq \{ \red{q_{i,\omega}}  \geq 0: T_{i,\omega} x_i +W_i \red{q_{i,\omega}}  = h_{i,\omega} \}\right\}, \label{two-stage-QP}
\end{align}
where $H_{i,\omega} $ is a symmetric positive semidefinite matrix. By Dorn duality \cite{dorn1960duality},
a    dual problem to \eqref{two-stage-QP}  has the following form:
\begin{align} \label{quar-rec-dual}
\begin{aligned}
\max_{u_{i,\omega},\pi_{i,\omega}}  \quad   &\quad  \phi_i(x_i,u_{i,\omega},\pi_{i,\omega}) \triangleq (    h_{i,\omega}-T_{i,\omega} x_i  )^T \pi_{i,\omega}-{1 \over 2} u_{i,\omega}^T H_{i,\omega}   u_{i,\omega} \\
 \mbox{subject to }  & \quad   (u_{i,\omega},\pi_{i,\omega}) \in S_{i,\omega} \triangleq \{u_{i,\omega},\pi_{i,\omega}: W_i^T\pi_{i,\omega}- H_{i,\omega}  u_{i,\omega} \ \leq d_{i,\omega}\}.
\end{aligned}
\end{align}\begin{assumption}~\label{assp-exist-quar}
 For any $x_i \in X_i$ and almost all $\omega \in \Omega$,

  a) there exists a strictly positive vector $\red{q_{i,\omega}}$ such that
 $ T_{i,\omega} x_i +W_i \red{q_{i,\omega}}  = h_{i,\omega} $;

 b) for any  $  \red{q_{i,\omega}} \in   \Upsilon_{i,\omega}, v_{i,\omega} \geq 0 $ with $ W_i v_{i,\omega}=0, v_{i,\omega}^TH_{i,\omega} v_{i,\omega}=0  $,  it holds that $ (H_{i,\omega} \red{q_{i,\omega}} +d_{i,\omega} )^Tv_{i,\omega}  \geq 0;$

 c)  $\mathbb{E}[\|T_{i,\omega}\|]<\infty$ for all $i \in \cal{N}$.

   \end{assumption}

Before   proceeding to show the  finiteness of $  {\cal Q}_i(x_i,\omega)$, we  first provide a lemma concerning the existence of optimal solutions for  quadratic programs. Consider a quadratic  program in the following form:
\begin{equation} \label{quadratic}
\min_{x \in \mathbb{R}^n} \left\{ \quad{1 \over 2} x^TDx+c^Tx \mid \quad  x \in \Delta (A,b) \triangleq \{x\geq 0:Ax=b\} \right\},
\end{equation}
where $D\in \mathbb{R}^{n \times n}$ is   positive semidefinite.
Then by \cite[Corollary 2.6]{lee2006quadratic}, we have the following lemma.
\begin{lemma}\label{lem-existence}~
The problem  \eqref{quadratic}  has a solution  if and only if $\Delta (A,b) $ is nonempty
and the following condition is satisfied:
$$ \rm{ If ~} x,v\in \mathbb{R}^n \rm{~are ~such~ that~}x \in \Delta (A,b) , v \geq 0, Av=0,
\rm{~and~}v^TDv=0, \rm{~then~}(Dx+c)^Tv\geq 0.$$
\end{lemma}
Before providing the result, we state a simple result about the quadratic
programming \uvs{which may be proved directly or by
	invoking~\cite[Cor.~2.3.7]{facchinei02finite}}.
\begin{lemma}\label{lem-quad-sing}~
Assume that problem  \eqref{quadratic}  has a solution set $S^*$.
Then the set $DS^*$ contains a singleton.
\end{lemma}
\uvs{We now analyze the subdifferential of the recourse function.}
 \begin{lemma}\label{lem-quar-rec}~
 Let  Assumption  \ref{assp-exist-quar}  hold.  Then we may claim  the following:

\noindent (a)
The optimal solution set  $ S_{i,\omega}^{\rm opt}(x_i)$ of the dual problem \eqref{quar-rec-dual} is nonempty,
and ${\cal Q}_i(x_i,\omega)=\max\limits_{ (u_{i,\omega},\pi_{i,\omega}) \in S_{i,\omega}}
  \phi_i(x_i,u_{i,\omega},\pi_{i,\omega}).$  Additionally,  ${\cal D}_i(x_i,\omega) = \left \{\pi_{i,\omega}^{\rm opt}(x_i): \left(u_{i,\omega}^{\rm opt}(x_i), \pi_{i,\omega}^{\rm opt}(x_i)\right) \in S_{i,\omega}^{\rm opt}(x_i) \right\}$  is a bounded set;

\noindent (b) ${\cal Q}_i(x_i,\omega)$  is  convex in   $x_i \in X_i$ for any given $\omega \in \Omega$;

\noindent (c) ${\cal Q}_i(x_i,\omega)$ is subdifferentiable at $x_i$ with
$\partial _{x_i}{\cal Q}_i(x_i,\omega)   = -T_{i,\omega}^T{\cal D}_i(x_i,\omega).$
 In particular, if ${\cal D}_i(x_i,\omega)$  consists of a   unique  point $ \pi_{i,\omega}^{\rm opt}(x_i)$, then
 ${\cal Q}_i(x_i,\omega)$ is  differentiable at $x_i$ with
$ \nabla_{x_i} {\cal Q}_i(x_i,\omega)   = -T_{i,\omega}^T \pi_{i,\omega}^{\rm opt}(x_i);$

\noindent (d)  there exists a positive constant $M_s>0$ such that
$ E[\| s_i \|^2]\leq M_s^2~~\forall s_i \in \partial {\cal Q}_i(x_i,\omega).$
 \end{lemma}

\noindent   {\bf Proof.}
(a) By invoking Assumption  \ref{assp-exist-quar} together   with Lemma \ref{lem-existence}, we have that ${\cal Q}_i(x_i,\omega) $
defined by problem   \eqref{two-stage-QP} has at least one optimal solution $\red{q_{i,\omega}^*}.$
Note that  $$T_{i,\omega} x_i +W_i \red{q_{i,\omega}}  = h_{i,\omega}  \Leftrightarrow \begin{pmatrix}
      &  W_i \\
      &  -W_i
\end{pmatrix} \red{q_{i,\omega}} \geq \begin{pmatrix}
      &  h_{i,\omega}- T_{i,\omega} x_i  \\
      &  -\left( h_{i,\omega}- T_{i,\omega} x_i\right)
\end{pmatrix}.$$
Then by    \cite[Theorem (Dual)-i]{dorn1960duality}  it is easily  seen   that
a solution $\left(u_{i,\omega}^{\rm opt}(x_i), \pi_{i,\omega}^{\rm opt}(x_i)\right)=\left(\red{q_{i,\omega}^*}, \pi_{i,\omega}^{\rm opt}(x_i)\right)$  exists to the dual problem \eqref{quar-rec-dual}    and
 ${\cal Q}_i(x_i,\omega)=\max_{ (u_{i,\omega},\pi_{i,\omega}) \in S_{i,\omega}}    \phi_i(x_i,u_{i,\omega},\pi_{i,\omega}).$

By  \cite[Theorem (Dual)-ii]{dorn1960duality} we see that
  for any  $ \left(u_{i,\omega}^{\rm opt}(x_i), \pi_{i,\omega}^{\rm opt}(x_i)\right) \in S_{i,\omega}^{\rm opt}(x_i)$,
  a solution $ \red{q^*_{i,\omega}}$  satisfying  $H_{i,\omega} \red{q^*_{i,\omega}}=H_{i,\omega}  u_{i,\omega}^{\rm opt}(x_i)$ to problem   \eqref{two-stage-QP}  also exists.  We denote by $ O_{i,\omega}^{opt}$ the optimal solution set of   problem \eqref{two-stage-QP}.  Then by Lemma \ref{lem-quad-sing}, the set  $ H_{i,\omega}  O_{i,\omega}^{opt}$ contains a    singleton denoted by  $H_{i,\omega} \red{q_{i,\omega}^*}$.  Thus,  for any  $ \left(u_{i,\omega}^{\rm opt}(x_i), \pi_{i,\omega}^{\rm opt}(x_i)\right) \in S_{i,\omega}^{\rm opt}(x_i)$,  $H_{i,\omega}  u_{i,\omega}^{\rm opt}(x_i)$ equals $H_{i,\omega} \red{q_{i,\omega}^*}$.
As a result,    by  the optimality condition for problem \eqref{quar-rec-dual} we see  that
 for any  $ \left(u_{i,\omega}^{\rm opt}(x_i), \pi_{i,\omega}^{\rm opt}(x_i)\right) \in S_{i,\omega}^{\rm opt}(x_i)$, $\pi_{i,\omega}^{\rm opt}(x_i)$ is a solution of the following linear program
\begin{align} \label{dual-max}
\max_{ \pi_{i,\omega}}   \quad   \left\{  (    h_{i,\omega}-T_{i,\omega} x_i  )^T \pi_{i,\omega}  \mid   W_i^T\pi_{i,\omega} \leq   H_{i,\omega} \red{q_{i,\omega}^*}+d_{i,\omega}\right\}.
\end{align}
Note that  since  $ h_{i,\omega} -T_{i,\omega} x_i $ is an interior point of the positive hull of matrix $W_i$ by  Assumption  \ref{assp-exist-quar}(a),  from
 \cite[p. 29]{shapiro09lectures}, the  optimal solution set
    of problem \eqref{dual-max} must be bounded. Hence  ${\cal D}_i(x_i,\omega) $ is bounded.

  (b)  Since    $ \phi_i(x_i,u_{i,\omega},\pi_{i,\omega})$ is convex in $x_i$ for all $(u_{i,\omega},\pi_{i,\omega}) \in S_{i,\omega}$,   for any $\alpha \in [0,1]$ and any $x_i,x_i' \in X_i$
\begin{align*}
\phi_i(\alpha x_i+(1-\alpha)x_i',u_{i,\omega},\pi_{i,\omega}) &\leq
\alpha \phi_i(x_i ,u_{i,\omega},\pi_{i,\omega})+ (1-\alpha)
\phi_i( x_i',u_{i,\omega},\pi_{i,\omega}) \\ & \leq \alpha {\cal Q}_i(x_i,\omega) +(1-\alpha) {\cal Q}_i(x_i',\omega),
\end{align*}
where the last inequality  follows from part (a).
Taking the maximum with respect to $(u_{i,\omega},\pi_{i,\omega}) \in S_{i,\omega}$ we obtain
${\cal Q}_i (\alpha x_i+(1-\alpha)x_i',\omega)  \leq \alpha {\cal Q}_i(x_i,\omega) +(1-\alpha) {\cal Q}_i(x_i',\omega)$,
and hence the recourse function ${\cal Q}_i(x_{i},\omega)$ is  convex in $x_i\in X_i$ for every
$\omega \in \Omega$.

(c) Note that  $\frac{\partial  \phi_i(x_i,u_{i,\omega},\pi_{i,\omega})}{ \partial x_i}=-T_{i,\omega}^T\pi_{i,\omega} $.
Since    ${\cal D}_i(x_i,\omega) $  is a compact convex set, by the proof of   Danskin's theorem \cite[Proposition A.22]{bertsekas1971control} we obtain
\begin{align*}\partial {\cal Q}_i(x_i,\omega)  &={\rm conv} \{ -T_{i,\omega}^Tu_{i,\omega}  :(u_{i,\omega},\pi_{i,\omega}) \in S_{i,\omega}^{\rm opt}(x_i)\} \\&    = {\rm conv} \{ -T_{i,\omega}^Tu_{i,\omega}  :u_{i,\omega} \in  {\cal D}_i(x_i,\omega)\}
  =-T_{i,\omega}^T {\cal D}_i(x_i,\omega).
\end{align*}
Hence  the  result  follows.

(d)  By   Assumption \ref{assp-exist-quar}(c) and the boundedness of ${\cal D}_i(x_i,\omega) $  we achieve   the result.
\hfill $\Box$

From Lemma \ref{lem-quar-rec}(b), $ \red{g_i(x_i)}:=c_i(x_i) + \mathbb{E}[{\cal Q}_i(x_i,\omega)]$ is convex. Then by recalling the fact that $f_1,\cdots, f_N$ admit  a contractive property,  $\wh{x}(\cdot) $   defined by \eqref{eq-prox-recourse}  is a contractive map.
 As a result, the  algorithm  and   convergence results  shown in Section \ref{sec-linear-recourse}   for  the problem  \eqref{SNash_rec}   with ${\cal Q}_i(x_i,\omega)$ defined by (Rec$_{\rm LP}(x_i)$)  still  carry over to the problem \eqref{SNash_rec}   with ${\cal Q}_i(x_i,\omega)$ defined by (Rec$_{\rm quar}(x_i)$).

 \section{Numerics}\label{sec:numerics}
 In this section, we demonstrate the performance of the proposed algorithms through   \red{an asset management problem and  a two-stage capacity expansion game.}  In Section \ref{sec-numerics-formulation}, we  formulate a
	portfolio  investment  problem and detailed descriptions of our
	simulations are provided in Section \ref{sec-numerics-simulation}.
	 \red{In  Section \ref{sec:num-2stage},  we  formulate a two-stage  competitive capacity expansion problem
	  and provide numerical simulations for the synchronous algorithm.}
 \subsection{Competitive portfolio selection} \label{sec-numerics-formulation}

The  portfolio selection problem has a long history in  the field of
modern financial theory, dating back to the seminal work of
Markowitz~\cite{markowitz1952portfolio}. We consider the case  with $N$ investors where
investor $i \in \mathcal{N}=\{1,\cdots,N\}$  may  invest an amount $x_{ij}$ in asset $j \in \{1,\cdots, n\}.$
Suppose that asset $j=1,\cdots,n$ has a return rate $r_j$, a random variable with expectation
$\nu_j=\mathbb{E}[r_j]$. Denote by $R=\mathbb{E}[(r-\nu)(r-\nu)^T]  \in
\mathbb{R}^{n\times n}$  the  covariance  matrix, where $r=(r_1,\cdots, r_n)^T  \in \mathbb{R}^n$ and $\nu=(\nu_1,\cdots, \nu_n)^T \in \mathbb{R}^n$. Then for investor $i \in \mathcal{N}$, the expected return is $\sum_{j=1}^n  \nu_j x_{ij}=\nu^T x_i$,   while the  variance (risk) of the return is $ \mathbb{E}[ \left( r^Tx_i-\nu^Tx_i\right)^2] =x_i^T Rx_i$,  where $x_i=(x_{i1},\cdots, x_{in})^T \in \mathbb{R}^n$. The expected return   and variance of
each asset $i  $ can be estimated from the history.  We consider a setting where each investor trades off between return and risk through the following utility:
    $U_i(x_i) \triangleq  \rho_i x_i^T Rx_i-\nu^T x_i$,
  where   $\rho_i>0$ is the \uvs{risk-aversion parameter} of investor
  $i$, and   larger  $\rho_i$ results in a
\uvs{larger emphasis on risk.}

 In a practical market  incorporating multiple investors, trades of
 diverse investors are usually pooled and  executed together~ \cite{o2006pooling,yang2013multi}. The transaction cost for a single investor may depend on the overall  trading levels  in the market and not just its own trading. Each investor  $i $ has an initial holding $x_i^0,$ the transaction  price is a function of  trade  size $x_j-x_j^0~\forall j \in \mathcal{N}$ from all investors, e.g., the price takes the form
 $ \phi(\omega) \sum_{i=1}^N (x_i-x_i^0)$,  where $\phi:\Omega \rightarrow \mathbb{R}^{n \times n}$ is a diagonal random   matrix with  diagonal entries being positive.   Besides,  the transaction
  cost of each investor is proportional  to its trade size. Based on the above considerations,  the portfolio selection  problem of  investor $i$  can be modeled as
   the following  stochastic program:
  \begin{align}\label{quad-payoff}
   \min_{x_i \in X_i} \quad \quad  & f_i(x_i,x_{-i})=U_i(x_i) + \mathbb{E}\left[(x_i-x_i^0)^T\phi(\omega) \sum_{j=1}^N (x_j-x_j^0) \right],    \end{align}
  where $X_i$  denotes  the set of feasible portfolios of player $i$.  Denote by  $\Phi \triangleq \mathbb{E}[\phi(\omega)]$.
    Then  for any $i \in \mathcal{N}$,
 $$\nabla_{x_i}^2 f_i(x)=2\rho_i R+2\Phi,~~\nabla_{x_ix_j}^2 f_i(x)= \Phi~~\forall j \neq i.$$
 Assume that the parameters $\rho_i, R, \Phi$ satisfy  the following condition
 \begin{equation}\label{exmp-cond}
\lambda_{\min} (2\rho_i R+2\Phi) > (N-1) \| \Phi\|~~\forall i \in \mathcal{N}.
\end{equation}
 \uvs{If condition  \eqref{exmp-cond} holds, then   $\|\Gamma\|_{\infty}<1$.
 Since both   $R $ and $\Phi$  are symmetric and positive definite
 matrices,  all of their eigenvalues  are strictly positive}. Then for sufficiently large $\rho_i$,  condition
 \eqref{exmp-cond}  certainly holds.   Notably,    large  $\rho_i$ means
 that  investor $i$ is more risk-averse.

  \subsection{Numerical results}\label{sec-numerics-simulation}
We consider a market of 6 investors and
  4 assets, i.e.,  $N=6, n=4.$  Set $\nu=(0.5,0.35,0.4,0.3)$ and let $R$ be a diagonal matrix with
  $R_{11}=0.16, R_{22}=0.1,R_{33}=0.12, R_{44}=0.09$.  Let each diagonal  element  of $\phi(\omega)$ be a random variable satisfying a uniform distribution $[0.12,0.18]$, and hence $\Phi=0.15\mathbf{I}_4.$
 Set  $\rho_i=3+i/N  $.   Let $X_i=\{x_i \in \mathbb{R}^n:x_{ij}\geq 0, x_{ij} \leq \textrm{cap}_{ij}, j=1,\cdots,n\}$.
  We set $ \textrm{cap}_{ij}=0.5$ and   the initial holdings $x_{ij}^0=0$ for all $i \in \mathcal{N}$ and $j =1,\cdots,n. $
  Then  Assumptions 1, 2 and 5 hold.
 Throughout this section, we assume that  the  empirical mean of the   error is calculated by averaging across  50 trajectories.

\subsubsection{Convergence of the Synchronous Algorithm}
We now present simulations for Algorithm \ref{inexact-sbr-cont} and examine its empirical
convergence rate and iteration complexity.
{Suppose  $j_{i,k}=\left \lceil \frac{1}{\eta^{2k}} \right \rceil$  steps of   (SA$_{i,k}$)
are taken at major iteration $k$ to get an inexact solution to    \eqref{def-wh_x_inft},
where  $\eta=a^{\kappa/2} $ for some   $\kappa>0$.
 Then $\alpha_{i,k}  \leq \eta^k$ by Lemma \ref{lem-rate-sa},   and hence the sequence $\{\alpha_{i,k}\}$ is summable by noting that $\eta \in (0,1)$. } We carry out  simulations for different selections of $ \mu,\kappa,$  for each case the smallest number of
 projected SG steps each player has carried out   to  obtain  $  u_k \leq $  2.5e-03  is  listed  in Table \ref{TAB1}, where $u_k$ is defined in \eqref{def-uk}. From this table, it is seen that  the parameters $ \mu,\kappa$   play  an  important rule in the overall iteration complexity.
   \begin{table} [!htb]
  \scriptsize
 \begin{center}
     \begin{tabular}{|c|c|c|c|c|c|c|c|c|}
        \hline
           \diagbox{$\mu$} { $\kappa$} &  0.8& 1.2 &1.6 &2&  2.4 & 2.8 &3.2  &3.6 \\ \hline
0.75 & 10240&7694&5908&  5076 & $4125$ & 4176 & 3237 &2941 \\ \hline
1 &8354& 5684 &4351&  $4192$ & $3018$  & 2428& 2782 &2740\\ \hline
1.25 &7266&5117&4051&  $3068$ & $3010$ &2724   & 1876&2635\\ \hline
1.5 &6178&4325 &3407 &  $2862$ & $2501$ &1910& 1813 &1928\\ \hline
1.75 & 6168&4173 &2993 &  $2732$ & $2178$ &1783 &1784&1514 \\ \hline
2 & 5572&  3765& 2709& 2652&  2300 &2054&1769 &1594 \\ \hline
2.25 & 5153& 3477 & 2774 &  $2608$ & $2090$ & 1956 & 1455 &1339\\ \hline
2.5&  4637& 3755& 2848& 2582& 1692& 1609& 1491&1433\\ \hline
      \end{tabular}
        \caption{Empirical iteration complexity for the synchronous algorithm with different  $\mu$ and $\eta$ \label{TAB1}}
\end{center}
\end{table}

 The rate of convergence of $u_k$ is   shown in Figure \ref{GOne},
 which demonstrates   that the iterates converge  in mean to the unique equilibrium at a linear rate.
{ We further  examine  the convergence behavior of Algorithm \ref{inexact-sbr-cont}  with $\mu=2,\eta=a$
 by plotting the  error bar  of   $x_{11,k}$,
where the  blue line denotes the trajectory of the mean  sequence computed via the sample average, the black dashed line denotes the equilibrium strategy $x_{11}^*$, while the red bars capture the variance of $x_{11,k}$.
   Figure \ref{GTwo}  shows that the  equilibrium estimate  approaches zero both with probability one  and in mean   with
   the    variance across the  samples  decaying  to zero.
   This corresponds well with the theoretical findings that the iterates  converge a.s.  to the NE  and that the variance
    converges to zero   by Proposition \ref{prp-as}  and  Proposition \ref{prp-mean-convergence}(b), respectively.}
The empirical relation between $\epsilon$ and $K(\epsilon)$ is shown in Figures \ref{fig:comp:a} and \ref{fig:comp:b}, where  $K(\epsilon)$ denotes  the  smallest number of    projected SG  steps  the  player has carried out to make  $ u_k<\epsilon.$ The  red  solid  curve   represents  the empirical  data, while the blue dashed curve demonstrates   the corresponding  quadratic  fit, which  \uvs{aligns} well with  the empirical  data.  It indicates that the  empirical  iteration complexity is of order $(1/\epsilon^2)$ in  this simulation setting.

 \begin{figure}[!htb]
 \centering
 \begin{minipage}{0.3\linewidth}
       \centering
  \includegraphics[width=2.1in]{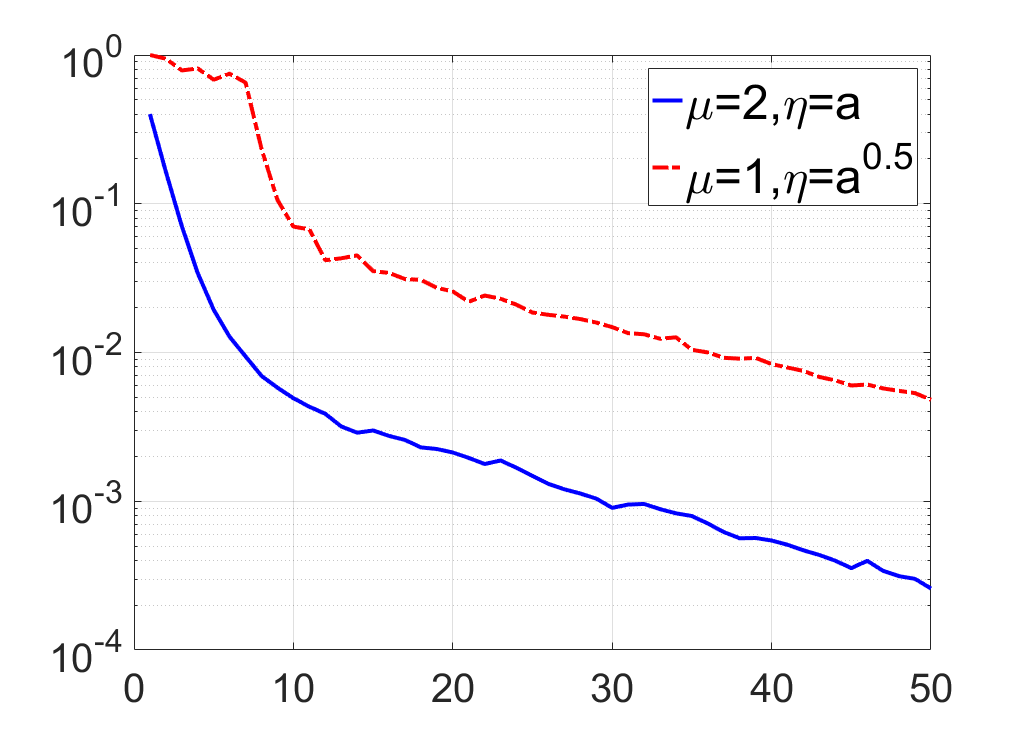}
  \caption{Linear Convergence of  Synchronous Algorithm}
   \label{GOne}
    \end{minipage}%
      \hspace{0.2in}
     \begin{minipage}{0.3\linewidth}
       \centering
 \includegraphics[width=2.1in]{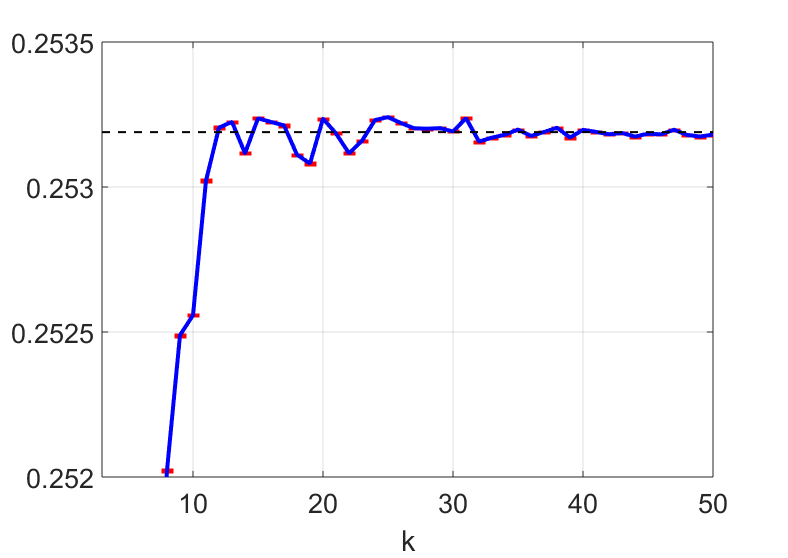}
 \caption{ Almost Sure Convergence  of Synchronous Algorithm with  Summable $\{\alpha_{i,k}\}$}     \label{GTwo}
    \end{minipage}%
      \hspace{0.2in}
     \begin{minipage}{0.3\linewidth}
       \centering
 \includegraphics[width=2.1in]{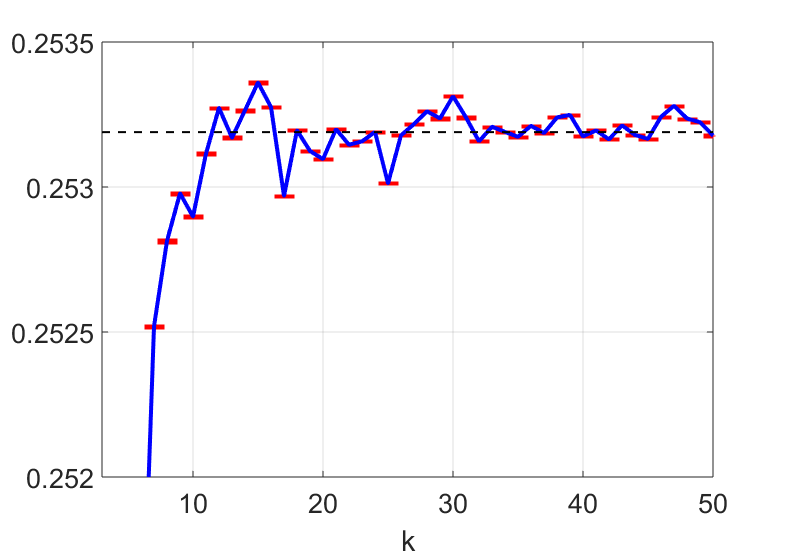}
 \caption{ Almost Sure Convergence   of Synchronous Algorithm with Unsummable $\{\alpha_{i,k}\}$ }     \label{GUnsum0}
    \end{minipage}%
\end{figure}

 While the previous numerical simulations    are all  for summable inexactness sequence  $\{\alpha_{i,k}\}$, we now investigate the performance   for   non-summable $\{\alpha_{i,k}\}.$
{Let $\mu=2$ and  suppose $j_{i,k}=k^2$  steps of (SA$_{i,k}$)
are taken at major iteration $k$ to get an inexact solution to    \eqref{def-wh_x_inft}.    Then
$\alpha_{i,k}=1/k$  by Lemma \ref{lem-rate-sa},    and hence  $\{\alpha_{i,k}\}$ is not summable.}
We plot  the  error bar  of   $x_{11,k}$ in Figure \ref{GUnsum0}, where the  blue line denotes the trajectory of the mean  sequence, the black dashed line denotes the equilibrium strategy $x_{11}^*$, and the red bars denotes the variance of $x_{11,k}$.
 It  shows that  the    variance across the  samples  decays  to zero, which is consistent with   Proposition \ref{prp-mean-convergence}(b)  by noting that $\alpha_{i,k}\to 0$.  Figure \ref{GUnsum0} also displays  that the estimates appear to converge a.s. to the NE though we cannot theoretically claim this from Proposition \ref{prp-as} since $\sum_{k=1}^{\infty }\alpha_{i,k}=\infty$. Further, the empirical relation between $\epsilon$ and $K(\epsilon)$ is shown in Figure  \ref{fig:comp:c} with the  red  solid  curve   representing  the empirical  data  and  the blue dashed curve demonstrating  its  quadratic  fit.  It is easily seen that   the complexity  bound of    non-summable $\{\alpha_{i,k}\}$ is worse than that of  summable   $\{\alpha_{i,k}\}$ as shown in Figure \ref{GThree}.
 \vspace{-0.2in}
\begin{figure}[!htb]
  \centering
  \subfigure[~Summable $\alpha_{i,k}$]{
     \label{fig:comp:a}
      \includegraphics[width=2in]{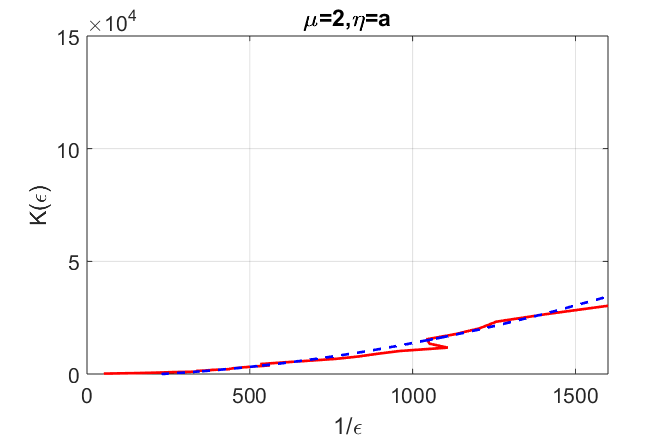} }
  \subfigure[~Summable $\alpha_{i,k}$]{
  \label{fig:comp:b} 
    \includegraphics[width=2in]{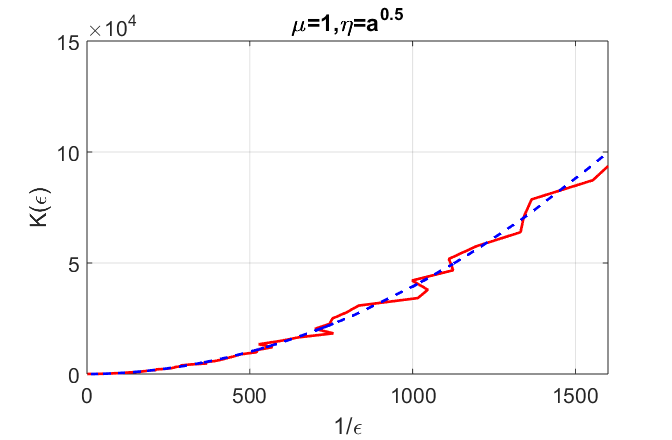} }
       \subfigure[~Non-summable $\alpha_{i,k}$]{
    \label{fig:comp:c} 
    \includegraphics[width=2in]{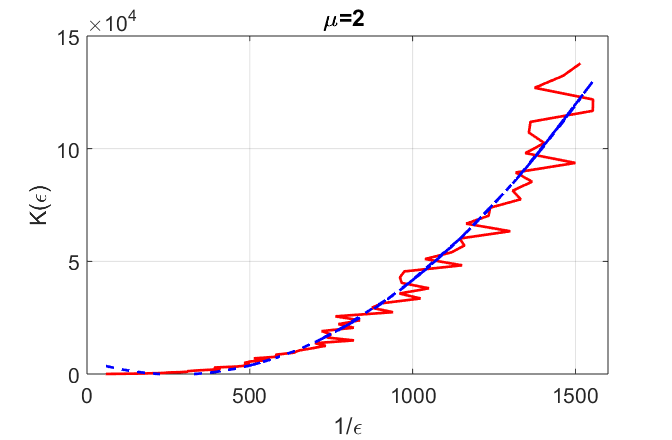}}
  \caption{Empirical Iteration Complexity and  Quadratic Fitting of Synchronous Algorithm}   \label{GThree}
\end{figure}

 \noindent  {\bf Comparisons with stochastic gradient method:}  Set    $\mu=2.5$ and
 suppose   $j_{i,k}=\left \lceil \frac{1}{a^{3k}} \right \rceil$  steps of   (SA$_{i,k}$)
are taken at major iteration $k$ to get an inexact solution to    \eqref{def-wh_x_inft}.
 We compare
both  Algorithm \ref{inexact-sbr-cont} and the standard  SG method for computing an NE}
in terms of the  iteration complexity   and communication overhead  for achieving  the same accuracy.
 The  empirically observed relationship between $\epsilon$ and $K(\epsilon)$ for both methods are shown in  Figure \ref{Gcompare}. From the figure,
  it can be observed that the  iteration complexity  are of the same orders while the constant of the SG method is superior to that of Algorithm \ref{inexact-sbr-cont}.
{Note that in  the delay-free SG method,
each player performs a single  projected gradient step by invoking a
communication  with  its rivals, consequently,  the  resulting  communication
overhead   is proportional  to the total number of  projected gradient steps.
In contrast,   the synchronous inexact proximal  algorithm  carries out
 an increasing number of player-specific projected gradient steps  after a single round of communication with its
rivals.  Note that this communication overhead is expected to be less in stochastic gradient schemes with delay~\cite{agarwal2011distributed}.
    The communication overhead  of   the SG method   and Algorithm \ref{inexact-sbr-cont} is  compared in Figure \ref{Gdata}.
From the results demonstrated in  Figure \ref{Gcompare} and  Figure  \ref{Gdata}, we  conclude that  Algorithm \ref{inexact-sbr-cont}
 compares  well with the SG method  in terms of overall projected gradient  steps while
  showing markedly less  communication overhead. In fact,  in certain applications, high communication overhead tends to render a scheme impractical.

\vspace{-0.2in}
\begin{figure}[H]
     \begin{minipage}[]{0.3\linewidth}
\includegraphics[width=2in]{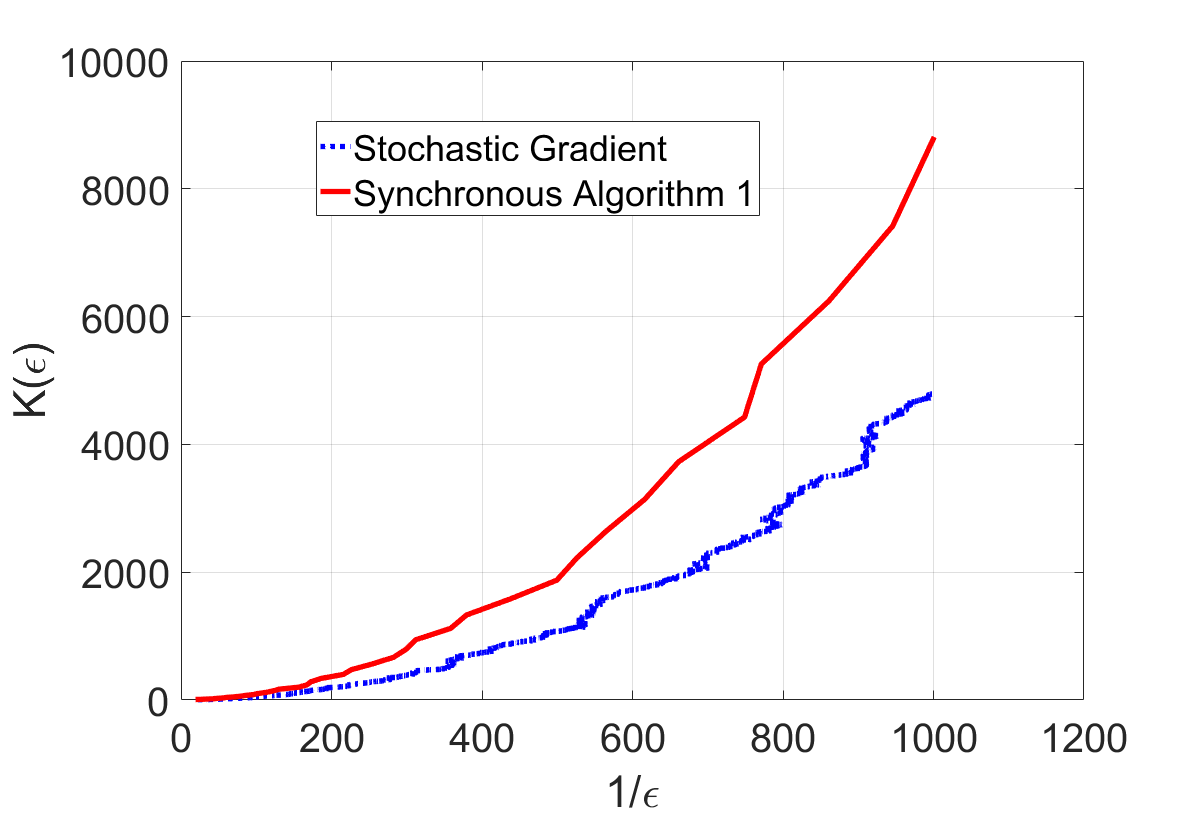}
  \caption{ Comparison of   Iteration Complexity} \label{Gcompare}
    \end{minipage}
       \hspace{0.2in}
   \begin{minipage}{0.3\linewidth}
       \centering
 \includegraphics[width=2in]{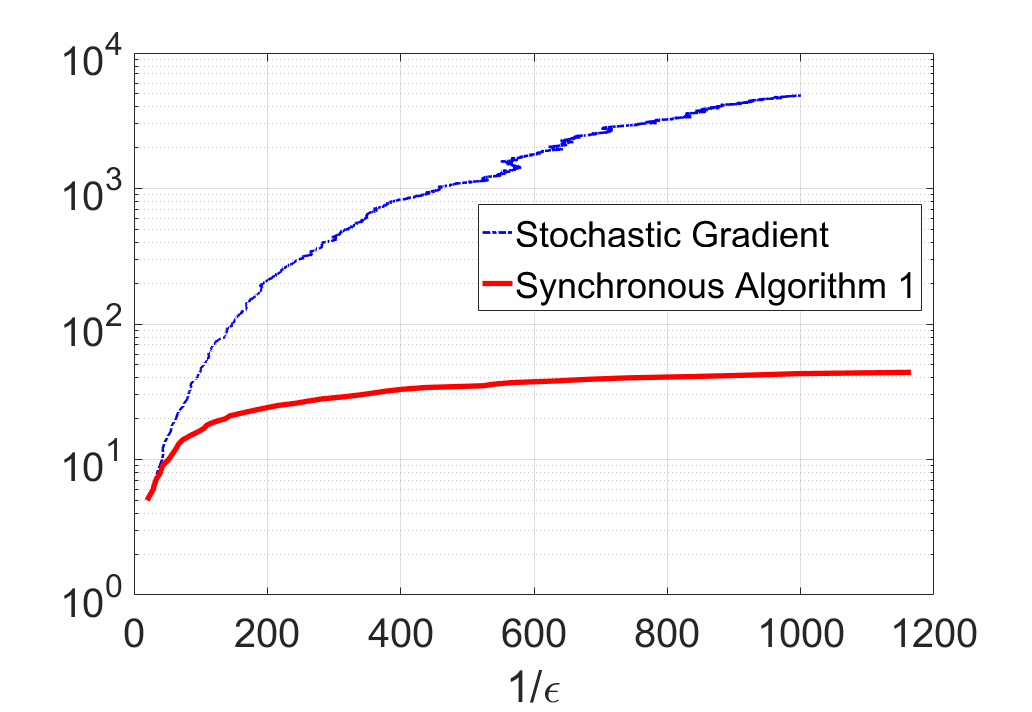}
   \caption{ Comparison of  Communication Overhead}   \label{Gdata}
    \end{minipage}%
       \hspace{0.2in}
          \begin{minipage}{0.3\linewidth}
\includegraphics[width=2in]{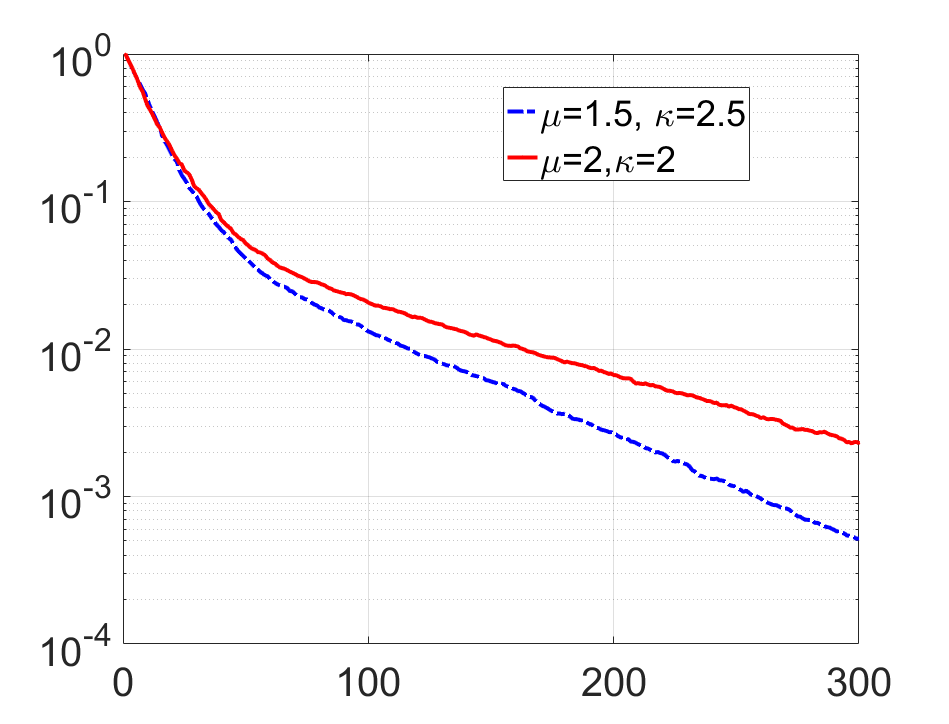}
  \caption{Linear Convergence of Randomized Algorithm } \label{rand_linear_convergence}
    \end{minipage}
 \end{figure}

\subsubsection{Convergence of the Randomized  and Asynchronous Schemes}
We  now run simulations for Algorithm \ref{rand-inexact-sbr-cont},  {where  $p_i=1/N~\forall i \in \cal{N}.$
Suppose  $j_{i,k}=\left \lceil \frac{1}{\eta^{2(\beta_{i,k}+1)}}\right \rceil$ steps of
 \eqref{sa-inner} are taken at major iteration $k$ to get an inexact solution   satisfying \eqref{rand-inexact-sub2}, where
   $\eta=a ^{\kappa/2} $ for some   $\kappa>0$.}
The trajectories  of  $u_k$  are shown in Figure \ref{rand_linear_convergence},  while   the empirical   mean of the  number of
projected SG steps    is  displayed in  Figure  \ref{rand_interation_complexity}.
It is seen that   randomized algorithm still  displays  linear convergence but its  empirical  iteration complexity of randomized algorithm  is  larger than that of the synchronous algorithm, a less surprising observation.

 \begin{figure}[!htb]
 \begin{minipage}{0.3\linewidth}
       \centering
 \includegraphics[width=2in]{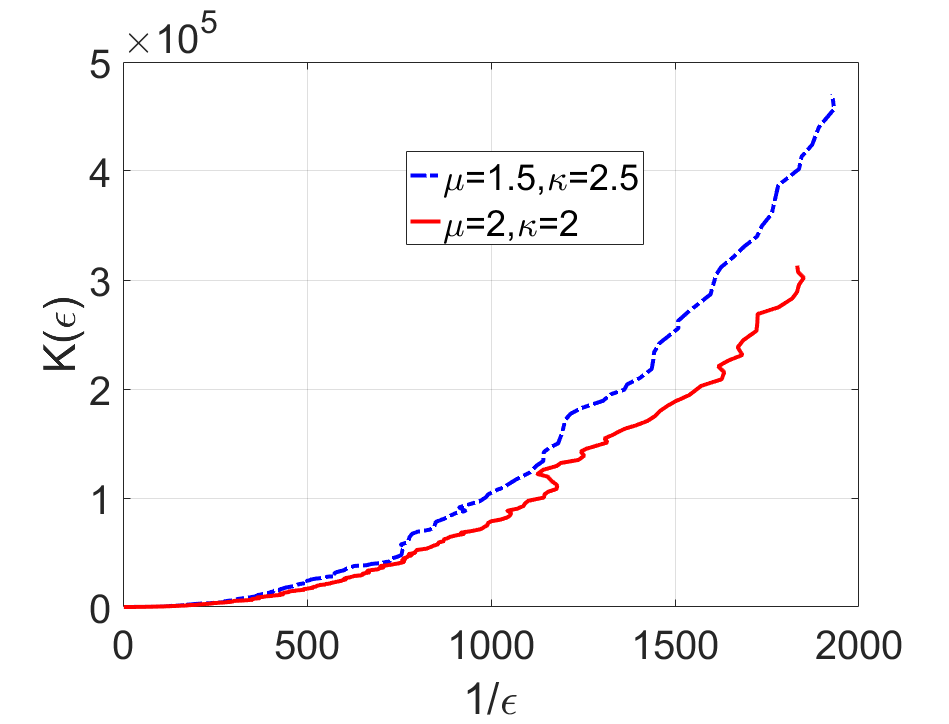}
   \caption{Iteration Complexity   of Randomized Algorithm}   \label{rand_interation_complexity}
    \end{minipage}%
       \hspace{0.2in}
     \begin{minipage}{0.3\linewidth}
\includegraphics[width=2in]{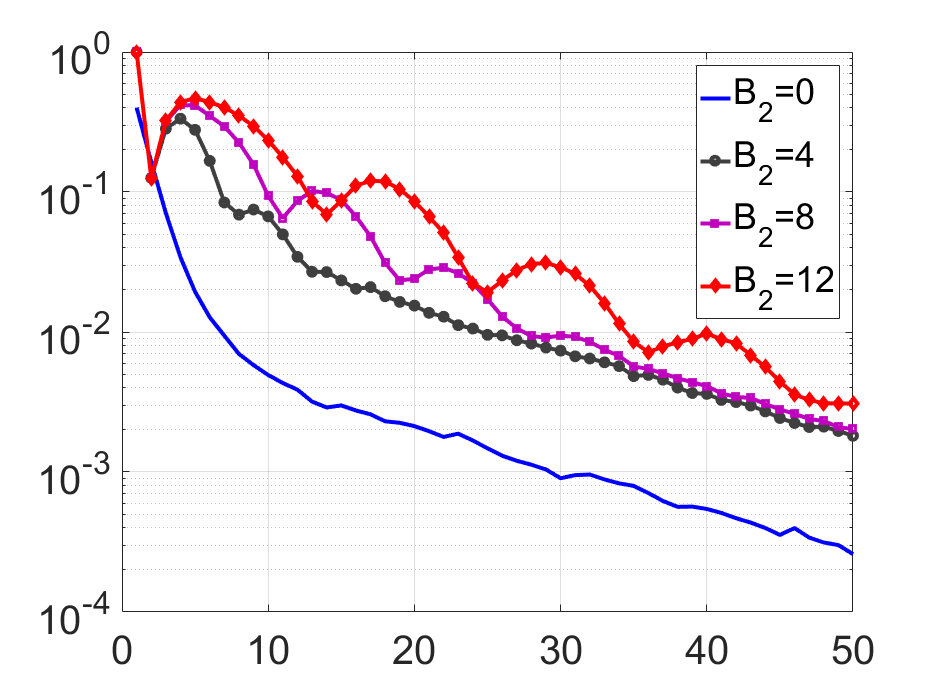}
  \caption{Linear Convergence  of Asynchronous  Algorithm  } \label{Asy_linear_convergence}
    \end{minipage}
       \hspace{0.2in}
   \begin{minipage}{0.3\linewidth}
       \centering
 \includegraphics[width=2in]{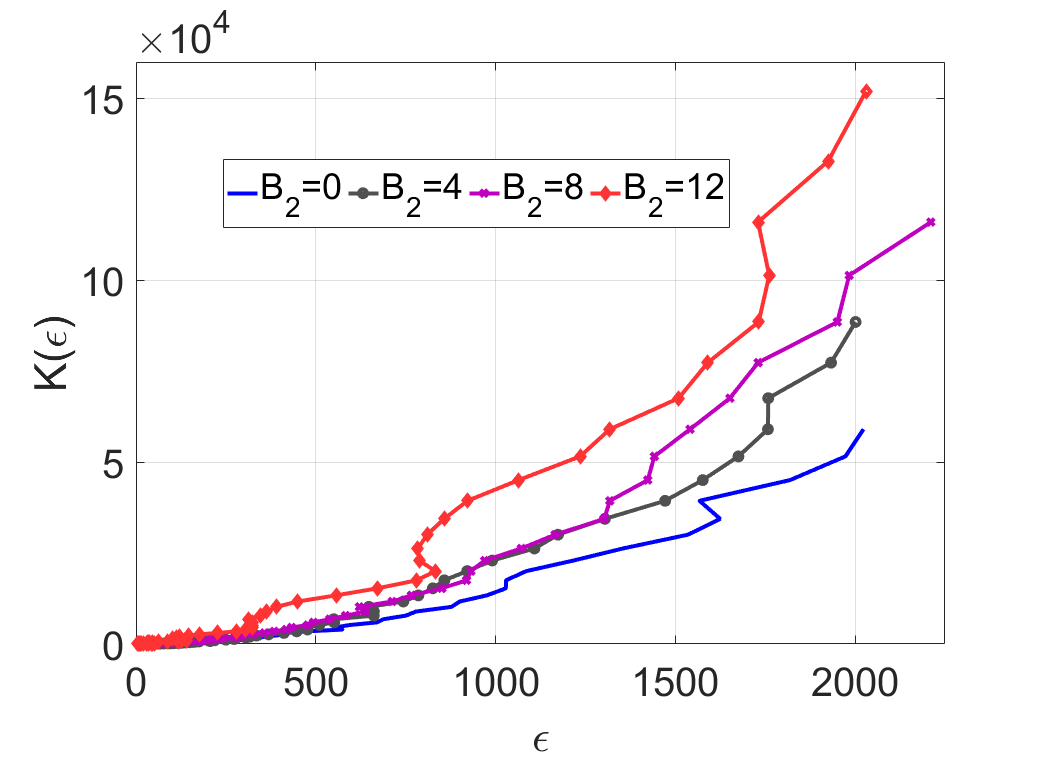}
   \caption{ Iteration Complexity  of Asynchronous  Algorithm }   \label{Asy_interation_complexity}
    \end{minipage}%
 \end{figure}
   Next, we simulate the performance of   Algorithm \ref{asy-inexact-sbr-cont},  where  for any
$k\geq0, i\in I_k, j \in \mathcal{N}$,  the communication delays
$ \tau_{ij}(k)$ are independently generated from a uniform distribution on the set  $\{0,1,\cdots, B_2\}$.
 Set $\mu=2, \eta=a_{\infty}, B_1=1$. Suppose that  $j_{i,k}=\left\lceil \frac{1}{\eta^{2k}} \right\rceil$ steps of    \eqref{SA2} are taken at major iteration $k$ to get an inexact solution to  problem \eqref{inexact-sub2}.
We carry out simulations  for  different communication delays    $B_2=0,4,8,12$.
  The trajectories  of  $u_k$  are shown in Figure \ref {Asy_linear_convergence},  \uvs{from which   it is observed that the iterates  produced by  Algorithm \ref{asy-inexact-sbr-cont}    converge in mean to the unique equilibrium  at a linear rate; however, the trajectories tend to display less of the monotonically decreasing behavior that characterizes the synchronous schemes.}
The empirically observed  iteration complexity is plotted  in  Figure  \ref{Asy_interation_complexity} \uvs{and we note that this worsens as $B_2$ grows.}

\subsubsection{Comparisons of Empirical and  Theoretical  Results}
 Set  $\rho_i=4~\forall i \in \mathcal{N}$.  Let each player take 40 BR steps  and  the inexact proximal   BR solution be computed via a SA scheme.  The  theoretical   and   empirical  rate of convergence   are shown in Table  \ref{TAB4},
  from which it is seen that the theoretical error bound is relatively conservative while the empirical error  is seen to be far smaller in practice.

  \begin{table} [htbp]
 \centering
 \scriptsize
 \begin{tabular}{|c|c|c|c|c|c|c|c|c|}
 \hline   \multicolumn{2}{|c|}{\multirow{2}{*}{parameters}} &  \multicolumn{2}{c|}{  synchronous } &
 \multicolumn{2}{c|}{randomized} & \multicolumn{2}{c|}{  asynchronous } \\  \cline{3-8}
\multicolumn{2}{|c|}{}  &  empirical  & theoretical &  empirical  & theoretical & empirical  & theoretical   \\ \hline
{\multirow{3}{*}{$\mu=1$}}  &  $\eta= a^{0.5}$     & 2e-03& 1.89&2.64e-03 &1.98e+01&  1.43e-03 & 1.36e+01   \\   \cline{2-8}
 &  $ \eta=a^{0.75} $     &  4.76e-04& 7.18e-01 &7.42e-04&1.73e+01&  3.24e-04 & 1.36e+01    \\ \cline{2-8}
 &  $ \eta=a$    &  1.08e-04&3.18e-01 &2.27e-04 &1.53e+01&  7.94e-05 & 1.37e+01   \\ \hline
{\multirow{3}{*}{$\mu=2$}} &
  $\eta=a^{0.5}$   &6.1e-03& 6.98e+00 &7.88e-03 & 2.49e+01& 4.33e-03 &  3.69e+01 \\ \cline{2-8}
 &  $\eta= a^{0.75}$     &  2.26e-03& 3.89e+00 &3.3e-03 &2.27e+01&   1.72e-03&3.69e+01   \\  \cline{2-8}
 &  $ \eta=a $     &  9.39e-04&2.33e+00 &1.39e-03 &2.09e+01&   6.9e-04&3.69e+01  \\ \hline
 {\multirow{3}{*}{$\mu=5$}} &   $\eta=a^{0.5}$    & 1.3e-02& 2.6e+01 &2.06e-02& 3.16e+01& 1.11e-02& 9.62e+01 \\ \cline{2-8}
 &  $ \eta=a^{0.75}$     &  7.5e-03& 2.01e+01&1.24e-02&3.01e+01& 6.55e-03 &9.62e+01  \\  \cline{2-8}
 &  $\eta= a $     &  4.8e-03&1.58e+01&8.4e-03 &2.89e+01& 4.2e-03 & 9.62e+01   \\ \hline
 \end{tabular}
 \vskip 2mm
     \caption{Comparison of theoretical and empirical error \label{TAB4}}
\end{table}

\subsection{Two-stage stochastic Nash games.}\label{sec:num-2stage}
{Consider a set of $N$ players denoted by $\mathcal{N} \triangleq \{1,\cdots,N\},$ where the $i$-th player solves the following two-stage problem
\begin{align}\label{example-2stage}
 \min_{x_i: 0\leq x_i\leq cap_i } \quad  & \underbrace{ C_i(x_i)-P(x) x_i}_{\triangleq f_i(x_i,x_{-i})}+ \mathbb{E}[Q_i(x_i,\omega)], \mbox{ where }
Q_i(x_i,\omega)\triangleq \red{\max_{ q_{i,\omega}: 0\leq  q_{i,\omega} \leq x_i}~ d_{i,\omega}  q_{i,\omega}-{h_{i,\omega}  \over 2} q_{i,\omega}^2} .
\end{align}
where $P(x)=a-b\sum_{i=1}^Nx_i,$ $C_i(\cdot)$ is a twice continuously  differentiable and $\eta_i$-strongly convex function of $x_i$.
 This can be viewed as a game where $N$ firms  compete  in Cournot in a  capacity market defined by an inverse demand function
 $P(x)$ and subsequently make production decisions  subject to a demand's capacity constraints while  faced by random prices and costs,
where   $P(x)$, the  market price, is a decreasing  function of total production, and $C_i(x_i)$  is the  cost function of  firm $i$.   Capacity markets are utilized to price generation capacity in power markets~\cite{smeers11generation,abada17multiplicity}. Note  that
\begin{align*} \bigtriangledown ^2 f=  \begin{pmatrix}
      & \bigtriangledown_{11}^2 f &\cdots &  \bigtriangledown_{1N}^2 f   \\
      & \vdots  & \ddots  &\vdots \\
        & \bigtriangledown_{N1}^2 f &\cdots &  \bigtriangledown_{NN}^2 f
\end{pmatrix} = \begin{pmatrix}
      & \bigtriangledown_{11}^2 C_1   \\
      &    & \ddots    \\
        &   & &  \bigtriangledown_{NN}^2 C_N
\end{pmatrix} +b\left(\mathbf{I}_N+\mathbf{1}_N\mathbf{1}_N^T \right).
\end{align*}
Since  $C_i(\cdot)$ is a twice continuously  differentiable and $\eta_i$-strongly convex function of $x_i$,
by definition \eqref{minmax-twice-diff} we have that for any $i\in \mathcal{N}$, $
\zeta_{i,\min} =\eta_i+2b \mbox{ ~and ~} \zeta_{ij,\max} =b  ~~      \forall j \neq i.$
Then  by \eqref{matrix-hessian}  the following holds:
\begin{align*}
&  \Gamma= \pmat{\frac{\mu}{\mu+\eta_1+2b} &
	\frac{b }{\mu+\eta_1+2b} & \hdots &
		\frac{b}{\mu+\eta_1+2b}\\
\frac{b}{\mu+\eta_2+2b} &
	\frac{\mu}{\mu+\eta_2+2b} & \hdots &
		\frac{b}{\mu+\eta_2+2b}\\
	\vdots & & \ddots & \\
\frac{b}{\mu+\eta_N+2b} &
	\frac{b}{\mu+\eta_N+2b} & \hdots &
		\frac{\mu}{\mu+\eta_N+2b}}
\end{align*}
Since for any $i\in \mathcal{N}$ and for any given $\omega$,  the function $Q_i(\cdot,\omega)$ is convex. Then $Q_i(\cdot)$ is convex in $x_i$,
and hence $\widehat{x}(\bullet)$ is contractive when    the spectral radius    $ \rho(\Gamma)< 1$,
for which one of  the sufficient conditions is $\| \Gamma\|_{\infty}<1.$  As such, the  assumption
$\min\limits_{i\in \mathcal{N}} \eta_i >(N-3)b$ ensure  that  the proximal BR map is contractive.

 By Dorn duality \cite{dorn1960duality},
a    dual problem to  the second stage problem \eqref{example-2stage} is as follows:
 \begin{equation}\label{dual-2stage}
\min_{u_{i,\omega},v_{i,\omega}} \left\{ \red{ \phi_i}(x_i,u_{i,\omega},v_{i,\omega}) = {\red{h_{i,\omega}}  \over 2}u_{i,\omega}^2 + x_{i} v_{i,\omega} \mid  \red{h_{i,\omega}}  u_{i,\omega}+v_{i,\omega}\geq \red{d_{i,\omega}}, {~\rm~} \red{v_{i,\omega}}\geq 0 \right\}.
\end{equation}
Similar to Lemma \ref{lem-quar-rec}, we can also show that
 the optimal solution set $S_{i,\omega}^{\rm opt}(x_i) $ of \eqref{dual-2stage} is bounded.
Note that $\frac{\partial  \phi_i(x_i,u_{i,\omega},\pi_{i,\omega})}{ \partial x_i}=v_{i,\omega} $,
  by invoking Danskin's theorem,  we obtain that
$\partial {\cal Q}_i(x_i,\omega)  ={\rm conv} \{ v_{i,\omega}   :(u_{i,\omega},v_{i,\omega}) \in S_{i,\omega}^{\rm opt}(x_i)\}.$}
In the simulations, we set $N=5,$  $\mu=1,$ $a= 2 , b=0.5, $ $cap_i=0.3+0.1\sqrt{i}  $,
and $\eta_i=(N-2.5)b.$
Suppose $C_i(x_i)={\eta_i x_i^2 \over 2}$. Let  \red{$h_{i,\omega}$ and $d_{i,\omega}$} be random variables satisfying the uniform distributions $[0.45,0.55]$ and $[0.3,0.4]$, respectively.
We now present simulations for  the synchronous algorithm and examine its empirical
convergence rate and iteration complexity, where the   empirical results are obtained  by averaging across  50 trajectories.
 Suppose  $j_{i,k}=\left \lceil \frac{1}{a^k} \right \rceil$  steps of  \red{\eqref{sa-recourse}}
are taken at major iteration $k$ to get an inexact solution to    \eqref{def-wh_x_inft}.
 The rate of convergence of $u_k$ is   shown in Figure \ref{2slinear},
 which demonstrates   that the iterates converge  in mean to the unique equilibrium at a linear rate.
 The empirical relation between $\epsilon$ and $K(\epsilon)$ is shown in Figure \ref{2scomplexity},
 \red{from which it is seen that    the empirical  data  aligns well with  the corresponding  quadratic  fit.
   It demonstrates  that the  empirical  iteration complexity is still of order $\mathcal{O}(1/\epsilon^2)$ in the  settings of a two-stage stochastic Nash game.}

 \begin{figure}[!htb]
 \centering
 \begin{minipage}{0.45\linewidth}
       \centering
  \includegraphics[width=2.1in]{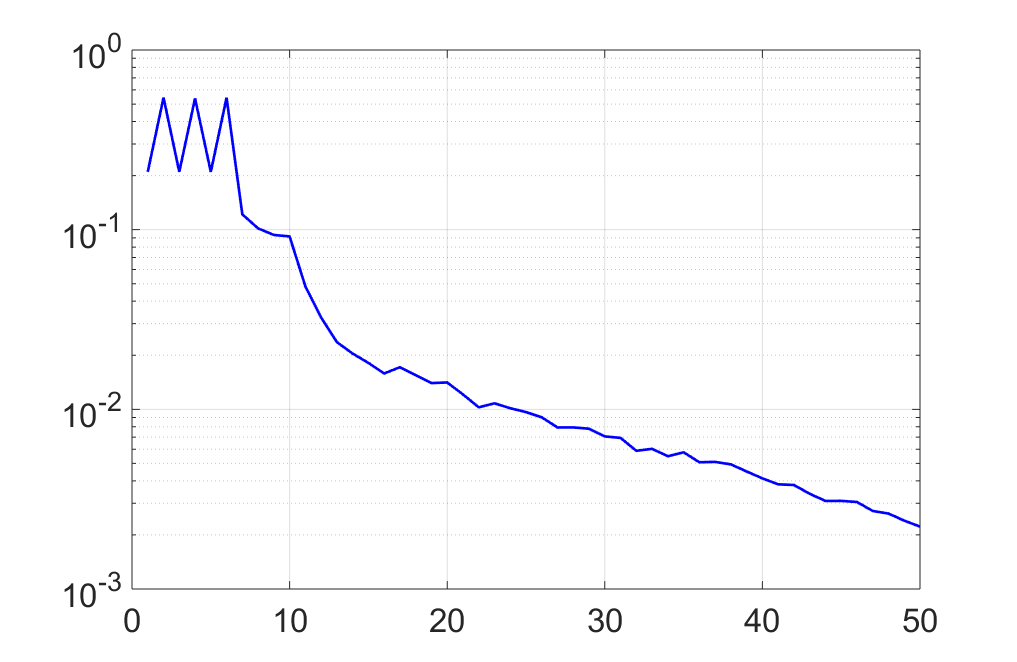}
  \caption{Linear Convergence }
   \label{2slinear}
    \end{minipage}%
      \hspace{0.2in}
     \begin{minipage}{0.45\linewidth}
       \centering
 \includegraphics[width=2.1in]{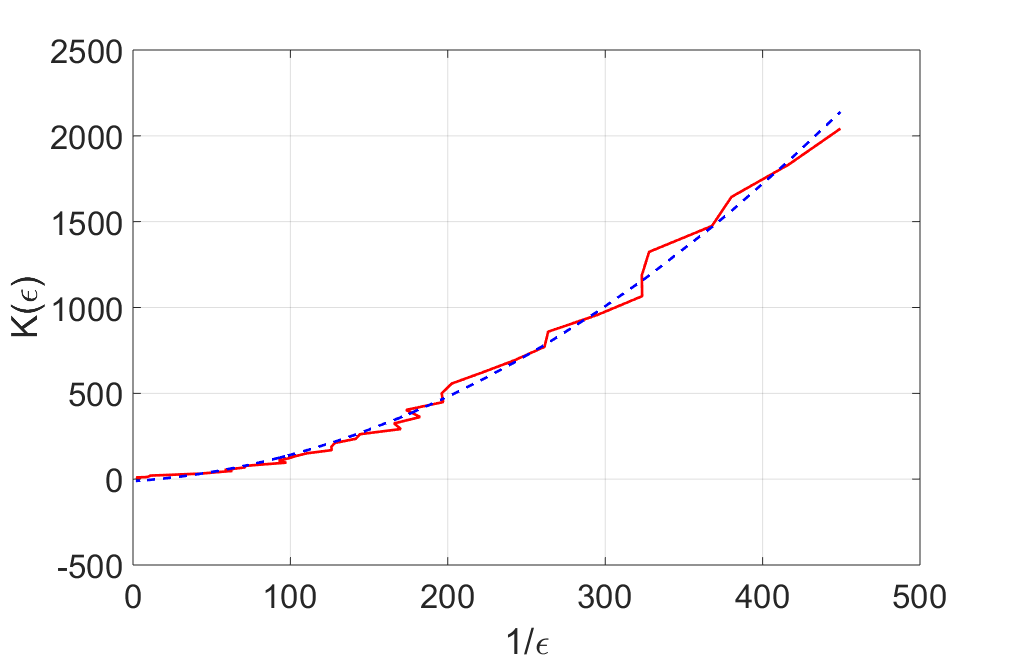}
 \caption{ Iteration Complexity }     \label{2scomplexity}
    \end{minipage}%
\end{figure}

\section{Concluding Remarks}\label{sec:conclusion}
This paper considers a class of  Nash games where each player's payoff function is expectation-valued.
We propose a synchronous inexact proximal   BR scheme  to solve the problem as well as a randomized and an asynchronous variant. Under  suitable contractive properties  on the proximal BR map,   we  separately  prove that
all proposed   schemes produce  iterates that converge in mean to the unique equilibrium  at a
linear rate.  In addition,  we derive  the overall iteration complexity for computing an $\epsilon$-Nash equilibrium  in
terms of projected gradient  steps.  Furthermore,   we consider generalizations that allow for private recourse by allowing for each player to solve a two-stage  stochastic  Nash game with either linear or quadratic second-stage problems.
Finally,  preliminary numerical studies support the theoretical findings in terms of asymptotic behavior and rate statements.

\appendix 
\numberwithin{equation}{section}
\section{Proof of Lemma  \ref{lem-as-asy}}\label{app-lem-4}
Define $t_{i,k}\triangleq x_{i,k}+\chi_{i,k} (\wh{x}_i (x_k)-x_{i,k}),  $ and
$ e_{i,k} \triangleq \chi_{i,k}\left(x_{i,k+1}-\wh{x}_i (x_k)\right)$.
Then $x_{i,k+1}=t_{i,k}+e_{i,k}.$
By  Assumption  \ref{assp-rand},  we have that
\begin{align}\label{cod-exp-t}
&\mathbb{E}\left[ \|t_{i,k}-x_i^*\|^2 \big | \mathcal{F}_k\right]=(1-p_i) \|x_{i,k}-x_i^*\|^2+p_i\| \wh{x}_i (x_k)-x_i^*\|^2.
\end{align}
{For $x=\{x_i\}_{ i=1}^N \in X $
and $P= {\rm diag} \left \{ {1\over p_1} \otimes \mathbf{I}_{n_1}, \cdots, {1\over p_N} \otimes \mathbf{I}_{n_N}  \right  \},$
 we  define the weighted norm
 \begin{equation}\label{weighted-norm}
\| x \|_P^2 =x^T Px=\sum_{i=1}^N    \|x_i\|^2/ p_i.
\end{equation}
Then by denoting $t_k=\{t_{i,k}\}_{ i=1}^N $,  from  \eqref{cod-exp-t} and \eqref{weighted-norm} it follows that
\begin{align}
\mathbb{E}\left[\|t_{k}-x^*\|_P^2 \big | \mathcal{F}_k\right]&=\sum_{i=1}^N {1 \over p_i}\mathbb{E}\left[ \|t_{i,k}-x_i^*\|^2 \big | \mathcal{F}_k\right] \notag
=\|x_{k}-x^*\|_P^2 +\sum_{i=1}^N \| \wh{x}_i (x_k)-x_i^*\|^2 -\sum_{i=1}^N  \|x_{i,k}-x_i^*\|^2 \notag
\\&\leq \|x_{k}-x^*\|_P^2  - (1-a^2)\  \|x_{k}-x^*\|^2, \label{asy-as-1}
\end{align}
where the last inequality follows by $\wh{x}(x^*)=x^*$ and  the inequality \eqref{cont-prox-resp}.  By \eqref{rand-inexact-sub2} and Assumption  \ref{assp-rand},   we see that for  any $ i \in \mathcal{N}$
 \begin{equation} \label{asy-as-02} \mathbb{E}\left[ \|e_{i,k} \|^2 \big | \mathcal{F}_k\right]= p_i \mathbb{E}\left[ \| x_{i,k+1}-\wh{x}_i (x_k) \|^2\big | \mathcal{F}_k\right]
 \leq p_i \alpha_{i,k}^2 ~~a.s. \end{equation}
{Then by the definition \eqref{weighted-norm},    we obtain   the following:
 \begin{align}\label{asy-as-2}
\mathbb{E}\left[ \|e_k\|_P^2 \big | \mathcal{F}_k\right] \leq \sum_{i=1}^N \alpha_{i,k}^2 ~~a.s.
\end{align}}
Note that by the Cauchy-Schwarz inequality,
$  \|x_{i,k+1}-x_i^*\|^2  \leq   \|t_{i,k}-x_i^*\|^2  + \|e_{i,k} \|^2 +2 \| t_{i,k}-x_i^*\| \|e_{i,k}\|.$}
Then  by  taking expectations conditioned on $\mathcal{F}_k$,  and by the condition Jensen's inequality, we have that
 \begin{equation}   \begin{split}
  \mathbb{E}\left[ \|x_{i,k+1}-x_i^*\|^2  \big | \mathcal{F}_k\right] &\leq  \mathbb{E}\left[ \|t_{i,k}-x_i^*\|^2 \big | \mathcal{F}_k\right]  +
\mathbb{E}\left[  \|e_{i,k} \|^2 \big | \mathcal{F}_k\right] +2  \sqrt{  \mathbb{E}\left[ \|t_{i,k}-x_i^*\|^2 \big | \mathcal{F}_k\right] \mathbb{E}\left[  \|e_{i,k} \|^2 \big | \mathcal{F}_k\right] } .\nonumber
\end{split}
\end{equation}
 {Since  $t_{i,k} \in X_i$,  by \eqref{diameter} and invoking the definition of weighted norm \eqref{weighted-norm} we have that
  \begin{equation}\label{asy-as-3}
  \begin{split}
  \mathbb{E}\left[ \|x_{k+1}-x^*\|_P^2  \big | \mathcal{F}_k\right] &\leq \mathbb{E}\left[\|t_{k}-x^*\|_P^2 \big | \mathcal{F}_k\right]
  +  \mathbb{E}\left[ \|e_k\|_P^2 \big | \mathcal{F}_k\right] +2 \sum_{i=1}^N  {   D_{X_i}   \over p_i}   \sqrt{\mathbb{E}\left[  \|e_{i,k} \|^2 \big | \mathcal{F}_k\right] }
\\& \leq  \|x_{k}-x^*\|_P^2  - (1-a^2)\  \|x_{k}-x^*\|^2 +\sum_{i=1}^N  \alpha_{i,k}^2 +2 \sum_{i=1}^N {  D_{X_i}    \over \sqrt{p_i}} \alpha_{i,k},
\end{split}
\end{equation}
where  the last inequality  is derived by \eqref{asy-as-1}, \eqref{asy-as-02}, and \eqref{asy-as-2}.
Since  $ 0 \leq \alpha_{i,k} <1$ and $\sum_{k=0}^{\infty}\alpha_{i,k}<\infty ~a.s.$,
  we have  that  $\sum_{k=0}^{\infty} \alpha_{i,k}^2<\infty ~~a.s.  $}
   Then  by the  Robbins-Siegmund theorem (\cite[Theorem 1]{robbins1985convergence}),
$\|x_{k}-x^*\|_P^2$ converges almost surely and
$\sum_{k=0}^{\infty} (1-a^2)\  \|x_{k}-x^*\|^2<\infty ~~a.s.$
Consequently, $\|x_{k}-x^*\|^2 $ converges to zero almost surely, and   hence we obtain the  result.
	 \hfill $\Box$

\section{Proof of Lemma  \ref{rand-geometric}} \label{app-lem-5}
Note that
   $\|x_{k}-x^*\|^2    \geq p_{\min} \sum_{i=1}^N   \|x_{i,k}-x_i^*\|^2 /p_i=p_{\min} \|x_{k}-x^*\|_P^2$,
where $p_{\min}=\min_{i \in \cal{N}} p_i.$  Then by \eqref{asy-as-1}, we  get
\begin{equation}\label{asy-as-a}
  \begin{split}
\mathbb{E}\left[\|t_{k}-x^*\|_P^2 \big | \mathcal{F}_k\right]
 \leq  \underbrace{\left(1-p_{\min}(1-a^2)\right)}_{\triangleq \tilde{a}^2}\|x_{k}-x^*\|_P^2 ~~a.s,\end{split}
\end{equation}
and hence by  the conditional  Jensen's inequality,  we obtain the following bound:
 $\mathbb{E}\left[\|t_{k}-x^*\|_P \big | \mathcal{F}_k\right] \leq  \tilde{a}\|x_{k}-x^*\|_P ~~a.s .$
Then by invoking that $x_{k+1}=t_k+e_k$,   the triangle inequality  and  Jensen's inequality, we have that
\begin{equation}\label{rand-bound1}
  \begin{split}
 \mathbb{E}\left[ \|x_{k+1}-x^*\|_P \right] &\leq \mathbb{E}\left[\|t_{k}-x^*\|_P  \right]  +\mathbb{E}\left[ \|e_k\|_P  \right]
 \leq \tilde{a} \mathbb{E}\left[\|x_{k}-x^*\|_P  \right] +\sqrt{\mathbb{E}\left[ \|e_k\|_P ^2 \right]}   . \end{split}
\end{equation}
Since $\mathbb{P}(\beta_{i,k}=m)= \left({k \atop m}\right) p_i^m(1-p_i)^{k-m}$ for all $  k\geq 1.$ Then for any $k \geq 1$,
\begin{equation}\label{asy-bound-inexct-sequnece}
  \begin{split}
  \mathbb{E} [\eta^{2\beta_{i,k}}] & =\sum_{m=0}^k \eta^{2m} \mathbb{P}(\beta_{i,k}=m)
  =\sum_{m=0}^k \left({k \atop m}\right) p_i^m(1-p_i)^{k-m} \eta^{2m}
  \\& =\sum_{m=0}^k \left({k \atop m}\right) (p_i\eta^2)^m(1-p_i)^{k-m} =
  (p_i\eta^2+1-p_i)^k\\& =(1-p_i(1-\eta^2))^k\leq(1-p_{\min}(1-\eta^2))^k   \triangleq \tilde{\eta}^{2k}~~\forall i \in \cal{N} .
\end{split}
\end{equation}
Then by $\beta_{i,0}=0$ and  $\alpha_{i,k}=\eta^{\beta_{i,k}+1}$, we obtain $\mathbb{E} [\alpha_{i,k}^2]   \leq \eta^2  \tilde{\eta}^{2k}~\forall k \geq 0 , $
and hence by \eqref{asy-as-2}  we have  that
 $$ \mathbb{E}\left[ \|e_k\|_P^2  \right] \leq \sum_{i=1}^N \mathbb{E}\left[ \alpha_{i,k}^2  \right] \leq N\eta^2\tilde{\eta}^{2k}~~\forall k \geq 0 .$$
Thus, by \eqref{rand-bound1}, we have the  following inequality for any $k \geq 1$:
   \begin{align*}
& \mathbb{E}\left[\|x_{k}-x^*\|_P \right]  \leq \tilde{a} \mathbb{E}\left[ \|x_{k-1}-x^*\|_P  \right] +\sqrt{N}\eta  \tilde{\eta}^{-1}  \tilde{\eta}^{k} .
\\&\leq    \tilde{a}^k |||x_{0}-x^*|||+ \sqrt{N}\eta  \tilde{\eta}^{-1}   \sum_{j=1}^k \tilde{a}^{k-j} \tilde{\eta}^j
   \leq   C\left(\sum_{i=1}^N p_i^{-1}\right)^{1/2} \tilde{c}^k +\sqrt{N}\eta  \tilde{\eta}^{-1}  k \tilde{c}^k,
\end{align*}
where  $\tilde{c}  \triangleq\max\{\tilde{a},\tilde{\eta}\}$.
Then by Lemma \ref{linear-bd},  we have that
 \begin{align*}  \mathbb{E}\left[\|x_{k}-x^*\|_P\right]
 & \leq   C\left(\sum_{i=1}^N p_i^{-1}\right)^{1/2} \tilde{q}^k +\sqrt{N}\eta  \tilde{\eta}^{-1}  D\tilde{q}^k
 \leq  \sqrt{N}( \widetilde{C} +\tilde{D})\tilde{q}^k,
\end{align*}
 where   $ \tilde{q} > \tilde{c}$,    $D \triangleq 1/\ln ((\tilde{q}/\tilde{c})^e),$ $ \widetilde{C}\triangleq C\left(\sum_{i=1}^N N^{-1}p_i^{-1}\right)^{1/2}$,  and  $ \widetilde{D} \triangleq  D\eta  \tilde{\eta}^{-1} $. \hfill $\Box$

 \section{Proof of Lemma \ref{asy-lem-gemo}.}\label{app-lem7}
For any $i \in I_k$,  by  the triangle  inequality  we have that
\begin{align}
   \quad \|x_{i,k+1}- x_i^*\|
  \leq  \| x_{i,k+1} -\wh{x}_i (y_k^i)\|
			+\|  \wh{x}_{i}(y_k^i) -
			\wh{x}_i (x^*)\|. \nonumber
\end{align}
Then by taking  expectations  conditioned on $\mathcal{F}_k$,
by  \eqref{inexact-sub2} and   the conditional Jensen's inequality,  we obtain:
\begin{align}\label{asy-condx-bd}
\mathbb{E}\left[ \|x_{i,k+1}- x_i^*\| \big | \mathcal{F}_k\right]
			 \leq   \alpha_{i,k}+ \mathbb{E}\left[\|\wh{x}_i (y_k^i) - \wh{x}_i (x^*)\|  \big | \mathcal{F}_k\right]\quad a.s.
			 \end{align}
Since $ \alpha_{i,k}$ is deterministic   by Assumption \ref{assp-asy}(a),
by taking   unconditional expectations on both sides of \eqref{asy-condx-bd},
and by invoking $y_{k}^i=(x_{1,k-\tau_{i1}(k)},\cdots,x_{N,k-\tau_{iN}(k)}  )$  and  \eqref{cont-prox-resp-2},
 we have that for any $i \in I_k$:
\begin{align}\label{update-expt-bd}
\mathbb{E}[ \|x_{i,k+1}- x_i^*\| ]
			 & \leq  \alpha_{i,k}+ a_{\infty} \max_ {j\in {\cal N}} \mathbb{E} [ \|x_{j,k-\tau_{ij}(k) }- x_j^*\| ]
			 . \end{align}			
  We now prove inequality \eqref{asy-geom-err}  by induction.
  It is obvious that  \eqref{asy-geom-err}  holds for  $k=0$ by  $\mathbb{E}[\|x_{i,0}-x_i^*\|] \leq C$ for all $i  \in \mathcal{ N}$.
 Inductively,  we assume that  \eqref{asy-geom-err}  holds for all $k$ up to some nonnegative integer $   \bar{k}$.

   We first prove     the following inequality:
 \begin{equation}\label{bound-delay}
\max_ {j\in {\cal N}} \mathbb{E} [ \|x_{j,k-\tau_{ij}(k) }- x_j^*\| ] \leq  (C+k) \rho^{ \max\{0,p-n_0\}}~~\forall k \in [p B_1 ,\bar{k}]  .
\end{equation}
Notice that $\lfloor \bar{k}/B_1 \rfloor=p.$
Then  by the induction  that    \eqref{asy-geom-err}  holds for all $k \leq   \bar{k}$, we have
 \begin{equation}\label{induct-bound}
\max_ {j\in {\cal N}} \mathbb{E} [ \|x_{j,k}- x_j^*\| ] \leq  (C+k) \rho^{ p}~~\forall k \in [p B_1 ,\bar{k}] .
\end{equation}
   If $B_2=0$, then $k-\tau_{ij}(k) =k$ and $n_0=0$   by  its definition.
   Then  by  \eqref{induct-bound},  we obtain the following:
   $$\max\limits_ {j\in {\cal N}}  \mathbb{E} [ \|x_{j,k-\tau_{ij}(k) }- x_j^*\| ] =\max\limits_ {j\in {\cal N}}  \mathbb{E} [ \|x_{j,k}- x_j^*\| ] \leq    (C+k) \rho^{ p} =(C+k) \rho^{ \max\{0,p-n_0\}}~~\forall k \in [p B_1 ,\bar{k}] .$$
   Thus, \eqref{bound-delay} holds  for $B_2=0$.
   If $B_2$ satisfies $ (n_0-1) B_1+1 \leq B_2 \leq  n_0B_1$ for some  positive integer $n_0  \geq 1$ and
   $k \in [p B_1 ,\bar{k}] ,$  then by  \eqref{bd-delay},  we obtain
        $$\max\{0 ,(p-n_0)B_1\} \leq  \max\{0 , k-n_0B_1 \}\leq \max\{0 , k-B_2\} \leq  k-\tau_{ij}(k) \leq k~~a.s.$$
      Then we have the following  for  all $k \in [p B_1 ,\bar{k}] $:
          \begin{align}\label{asy-bound-exp}
  \max_ {j\in {\cal N}}  \mathbb{E} [ \|x_{j,k-\tau_{ij}(k) }- x_j^*\| ] &\leq \max_{ \max\{0, (p-n_0)B_1\} \leq t \leq k}  ~\max_ {j\in {\cal N}} \mathbb{E} [ \|x_{j,t}- x_j^*\| ]  .
  \end{align}
  We  consider the following two possible cases:

 \noindent i) If  $\max\{0, (p-n_0)B_1\}=0$, then  by  the inductive assumption, we have that   for any $ k \in [p B_1 ,\bar{k}] $:
     \begin{align*}
  \max_{ \max\{0, (p-n_0)B_1\} \leq t \leq k}  ~\max_ {j\in {\cal N}} \mathbb{E} [ \|x_{j,t}- x_j^*\| ]
    &  = \max_{ 0 \leq t \leq k}  ~\max_ {j\in {\cal N}} \mathbb{E} [ \|x_{j,t}- x_j^*\| ]
     \\&  \leq \max_{ 0 \leq t \leq k}   (C+t)\rho^{\lfloor  {t \over B_1} \rfloor} \leq C+k\quad {\scriptstyle \left(  \mathrm{since} \ \rho \in (0,1) \right) } .
  \end{align*}
  ii) If  $\max\{0, (p-n_0)B_1\}=(p-n_0)B_1$, then  by  the inductive assumption, we have that  for any $ k \in [p B_1 ,\bar{k}] $:
     \begin{align*}
  \max_{\max\{0, (p-n_0)B_1\} \leq t \leq k}  ~\max_ {j\in {\cal N}} \mathbb{E} [ \|x_{j,t}- x_j^*\| ]
   &  = \max_{ (p-n_0)B_1 \leq t \leq k}  ~\max_ {j\in {\cal N}} \mathbb{E} [ \|x_{j,t}- x_j^*\| ]
     \\&  \leq \max_{ (p-n_0)B_1 \leq t \leq k}   (C+t)\rho^{\lfloor  {t \over B_1} \rfloor} \leq (C+k)\rho^{p-n_0}\quad {\scriptstyle \left(  \mathrm{since} \ \rho \in (0,1) \right) } .
  \end{align*}
  Combining cases i) and   ii), by \eqref{asy-bound-exp},  we  have  the following bound:
  $\max_ {j\in {\cal N}}  \mathbb{E} [ \|x_{j,k-\tau_{ij}(k) }- x_j^*\| ]   \leq  (C+k) \rho^{\max\{0,  p-n_0\}}$ ,
  and hence  \eqref{bound-delay} holds  for $B_2\geq 1$.
       Consequently,  we have shown   \eqref{bound-delay}.

We will  validate  that    \eqref{asy-geom-err}   holds  for  $k=\bar{k}+1$ by considering two cases:
  $\bar{k}=( p +1)B_1-1$  or $\bar{k} \in [ p B_1 ,( p +1)B_1-1)$ for some nonnegative integer $ p$.\\
       \noindent {\bf Case 1:  $\bar{k}=( p +1)B_1-1$.}
 By the inductive assumption, it is seen that  \eqref{asy-geom-err}  holds for any $k \in  [ p B_1 ,( p +1)B_1)$.
Note that each player updates at least once during any time interval of length $B_1 $ by Assumption \ref{assp-asy}.
 Then player $i$ updates its strategy  at least $p$ times in the  time interval $[0, pB_1), $ and  there exists at least one integer  $k_i \in [ p B_1 ,( p +1)B_1)$ such that $i \in I_{k_i}  $.    Set  $k_i$ to be the largest  integer  in  the set $ [ p B_1 ,( p +1)B_1)$  such that $i \in I_{k_i}  $. Then  $\beta_{i,k_i} \geq   p+1 $, hence by   \eqref{update-expt-bd}  and  \eqref{bound-delay},   we  derive the following:
\begin{equation}\label{asy-bpund-xk}
\begin{split}
\mathbb{E}[ \|x_{i, k_i+1} - x_i^*\| ]
			 & \leq   \alpha_{i,k_i} + a_{\infty} \max_ {j\in {\cal N}} \mathbb{E} [ \|x_{j,k_i-\tau_{ij}(k_i)}- x_j^*\| ] \\&
			  \leq   \eta^{ p+1}+ a_{\infty}  (C+k_i) \rho^{\max\{  p-n_0,0\}}
			  \quad {\scriptstyle \left(\mathrm{since ~}\alpha_{i,k_i}=\eta^{ \beta_{i,k_i}} \rm{~and~}  \eta\in(0,1) \right)}
	 \\& \leq   \rho^{ p+1 }+  (C+k_i)\rho^{  p+1}
	 \quad {\scriptstyle \left(\mathrm{since ~} 1>\rho \geq \max\{ a_{\infty} \rho^{-n_0}, \eta\}   \right)}
\\& \leq  (C+k_i+1)\rho^{  p+1} .\end{split}
\end{equation}	
 Therefore,  for any $i \in {\cal N}$, by the selection of $k_i$ and  the definition of Algorithm  \ref{asy-inexact-sbr-cont},  we  know that
 \begin{align*}
 \mathbb{E}[ \|x_{i, k+1} - x_i^*\| ]  &=\mathbb{E}[ \|x_{i, k_i+1} - x_i^*\| ]
 \leq    (C+k_i+1)\rho^{  p+1}   \\&
 \leq   (C+k+1)\rho^{  p+1} ~~\forall k: k_i \leq k<( p +1)B_1.
 \end{align*}
This implies  that
 $\max\limits_ {i\in {\cal N}} \mathbb{E}[ \|x_{i, ( p+1)B_1} - x_i^*\| ]
  \leq  (C+( p+1)B_1)\rho^{  p+1} ,$
 and hence \eqref{asy-geom-err}  holds for $k= \bar{k}+1$ when  $\bar{k}=( p +1)B_1-1$.

 \noindent {\bf Case 2: $\bar{k} \in [ p B_1 ,( p +1)B_1-1)$.}
Since  $\beta_{i,k_i} \geq   p+1 $ for  any   $i \in I_{\bar{k}}  $,  from  \eqref{update-expt-bd}  and  \eqref{bound-delay} it follows that
\begin{align*}
\mathbb{E}[ \|\wh{x}_{i, \bar{k}+1} - x_i^*\| ]
			 & \leq    \alpha_{i,\bar{k}} + a_{\infty} \max_ {j\in {\cal N}} \mathbb{E} [ \|x_{j,k-\tau_{ij}(\bar{k})}- x_j^*\| ] \\&
			 \leq  \eta^{  p+1}+a_{\infty}  (C+\bar{k}) \rho^{\max\{  p-n_0,0\}}
			 \quad {\scriptstyle \left(\mathrm{since ~} \alpha_{i,\bar{k}}=\eta^{ \beta_{i,\bar{k}}} \rm{~and~}  \eta\in(0,1) \right)}
			 \\& \leq \rho^{  p+1}+ (C+\bar{k})\rho^{  p+1}
	 \quad {\scriptstyle \left(\mathrm{since ~} 1>\rho \geq \max\{ a_{\infty} \rho^{-n_0}, \eta\}   \right)}
  \\&   \leq (C+\bar{k}+1)\rho^{  p+1} <(C+\bar{k}+1)\rho^{  p }~~\forall i \in I_{\bar{k}}   .
			 \end{align*}
 For any   $i \notin I_{\bar{k}}  $, by \eqref{induct-bound}, we  have
$\mathbb{E}[ \|x_{i, \bar{k}+1} - x_i^*\| ]  = \mathbb{E}[ \| x_{i, \bar{k}} - x_i^*\| ]   \leq
   (C+\bar{k} )\rho^{  p} <(C+\bar{k}+1)\rho^{  p}$   ,
and hence    \eqref{asy-geom-err}  holds  for $k=\bar{k}+1.$

By combing Cases 1 and   2,    \eqref{asy-geom-err}  holds  for $k=\bar{k}+1.$
Thus, by  induction,  we obtain   \eqref{asy-geom-err}  for any $k \geq 0$.

Since   $\left\lfloor  {k \over B_1} \right\rfloor \geq  {k \over B_1} -{B_1-1 \over B_1} $ and $~0<\rho<1$,
by \eqref{asy-geom-err}  we derive
\begin{equation}\label{upper-bound-estimate}
\max_{i\in \mathcal{N}}\mathbb{E}[\|x_{i,k} -x_i^*\|]
 \leq (C+k)\rho^{   {k \over B_1} -{B_1-1 \over B_1}}
 =\rho^{ -{B_1-1 \over B_1}} (C+k)\rho^{k \over B_1}=\rho^{ -{B_1-1 \over B_1}} (C+k)c^k,~~  \forall k \geq 0.
\end{equation}
 From Lemma~\ref{linear-bd}, there exist scalars
$q$ and $D$  satisfying  $q \in (c,1)$ and
$D \red{\geq} 1/\ln((q/c)^e)$ such that
$ kc^k  \red{\leq} D q^k~~\forall k \geq 0, $
which incorporating with \eqref{upper-bound-estimate}  and $c^k <q^k$ yields \eqref{upper-asy-geom-err}.   \hfill $\Box$

\bibliographystyle{IEEEtran}\def\cprime{$'$}

\end{document}